\documentclass[secthm,amsthm]{elsart}

\usepackage{amssymb, amsmath}

\makeatletter
 \def\itemize{%
  \ifnum \@itemdepth >\@maxlistdepth
    \@toodeep
  \else
    \advance\@itemdepth \@ne
    \edef\@itemitem{labelitem\romannumeral\the\@itemdepth}%
     \setleftmargin{i}{\qquad}%
     \setleftmargin{ii}{\qquad}%
    \list{\csname\@itemitem\endcsname}%
       {\let\makelabel\right@label}
  \fi}
 \makeatother

 \makeatletter \thm@notefont{\rm} \makeatother

\newcommand{\iso}{\cong}
\newcommand{\integer}{\mathord{\mathbf{Z}}}
\newcommand{\real}{\mathord{\mathbf{R}}}
\newcommand{\torus}{\mathbf{T}}

\newcommand{\quot}{G/\langle s \rangle}

\newcommand{\script}{\mathcal}
\newcommand{\hamilton}{\script H}
\newcommand{\flow}{\script F}
\newcommand{\basics}{\script Q}
\newcommand{\even}{\script E}
\newcommand{\boxprod}{\mathbin{\mbox{\tiny$\square$}}}

\DeclareMathOperator{\Cay}{Cay}
\DeclareMathOperator{\imb}{imb}
\DeclareMathOperator{\kw}{kw}
\DeclareMathOperator{\len}{len}
\DeclareMathOperator{\wt}{wt}
\DeclareMathOperator{\dist}{dist}
\DeclareMathOperator{\knot}{knot}
\DeclareMathOperator{\wli}{wli}

\newcommand{\Xdouble}{\widehat{X}}
\newcommand{\Cdouble}{\widehat{C}}
\newcommand{\Pdouble}{\widehat{P}}
\newcommand{\Qdouble}{\widehat{Q}}
\newcommand{\vdouble}{\widehat{v}}
\newcommand{\wdouble}{\widehat{w}}
\newcommand{\ldouble}{\widehat{I}}

\newcommand{\Xtilde}{\widetilde{X}}
\newcommand{\Ctilde}{\widetilde{C}}
\newcommand{\Ptilde}{\widetilde{P}}

\newcommand{\Itilde}{\widetilde{I}}

\newcommand{\Ho}{H_{\mathrm{o}}}
\newcommand{\He}{H_{\mathrm{e}}}

 \newcommand{\bigset}[2]{\left\{\, #1
 \mathrel{\left| \vphantom {\left\{ #1 \mid #2 \right\} } \right.}
 #2 \,\right\} }

\newbox\GRID \setbox\GRID\hbox{$P_3 \boxprod P_3 \boxprod P_2$}

\newcommand{\pref}[1]{\mbox{\upshape(\ref{#1})}}
\renewcommand{\eqref}[1]{Eq.~\pref{#1}}
\newcommand{\fullref}[2]{\ref{#1}\pref{#1-#2}}
\newcommand{\see}[1]{(see~\ref{#1})}
\newcommand{\seeSect}[1]{(see Section~\ref{#1})}
\newcommand{\seeDefn}[1]{(see Definition~\ref{#1})}
\newcommand{\fullsee}[2]{(see~\ref{#1}\pref{#1-#2})}
\newcommand{\seeand}[2]{(see~\ref{#1} and~\ref{#2})}
\newcommand{\cf}[1]{(cf.~\ref{#1})}

\theoremstyle{definition}
 \newtheorem{assump}[thm]{Assumption}
 \newtheorem{notation}[thm]{Notation}
 \newtheorem{obs}[thm]{Observation}
 \newtheorem{eg}[thm]{Example}

\theoremstyle{remark}
\newtheorem{rmk}{Remark\ignorespaces}
\newtheorem{warn}{Warning\ignorespaces}
\newtheorem{sugg}{Suggestion\ignorespaces}

\renewcommand{\theequation}{{\rm E}\arabic{equation}}

 \newenvironment{proofclaim}[1][\unskip]{\em
 \medskip \noindent Claim.\ }{\unskip\upshape}

 \newcounter{step}
 \newenvironment{step}[1][\unskip]{\refstepcounter{step}\em
 \medskip \noindent Step \thestep\ #1.\ }{\unskip\upshape}
 \renewcommand{\thestep}{\arabic{step}}

 \newcounter{proofcase}
 \newenvironment{proofcase}[1][\unskip]{\refstepcounter{proofcase}\em
 \medskip \noindent Case \thecase\ #1.\ }{\unskip\upshape}
 \renewcommand{\thecase}{\arabic{proofcase}}

 \newcounter{subcase}
 \newenvironment{subcase}[1][\unskip]{\refstepcounter{subcase}\em
 \medskip \noindent Subcase \thesubcase\ #1.\ }{\unskip\upshape}
 \numberwithin{subcase}{proofcase}

 \newcounter{subsubcase}
 \newenvironment{subsubcase}[1][\unskip]{\refstepcounter{subsubcase}\em
 \medskip \noindent Subsubcase \thesubsubcase\ #1.\ }{\unskip\upshape}
 \numberwithin{subsubcase}{subcase}

\begin{document}

\begin{frontmatter}

\title
 {Flows that are sums of hamiltonian cycles 
 in Cayley graphs on abelian groups}

{To Brian Alspach on his sixty-fifth birthday}

\author[osu,leth]{Dave Morris\thanksref{DaveWitte}},
 \thanks[DaveWitte]{Former name: Dave Witte.}
 \ead{dmorris@cs.uleth.ca}
 \ead[url]{http://www.cs.uleth.ca/$\sim$morrisd/}
 \author[leth]{Joy Morris},
 \ead{morris@cs.uleth.ca}
 \ead[url]{http://www.cs.uleth.ca/$\sim$morris/}
 \author[ccr]{David Petrie Moulton}
 \ead{moulton@idaccr.org}

 \address[osu]{Department of Mathematics, Oklahoma State University, Stillwater,
OK 74078, USA}
 \address[leth]{Department of Mathematics and Computer Science, University of
Lethbridge, Lethbridge, AB, T1K~3M4, Canada}
 \address[ccr]{Center for Communications Research,
 805 Bunn Drive, Princeton, NJ 08540, USA}

\begin{abstract}
 If $X$ is any connected Cayley graph on any finite abelian group, we  determine
precisely which flows on~$X$ can be written as a sum of hamiltonian cycles. (This
answers a question of B.~Alspach.) In particular, if the degree of~$X$ is at
least~$5$, and $X$~has an even number of vertices, then the flows that can be so
written are precisely the even flows, that is, the flows~$f$, such that
$\sum_{\alpha \in E(X)} f(\alpha)$ is divisible by~$2$. On the other hand, there
are examples of degree~$4$ in which not all even flows can be written as a sum of
hamiltonian cycles. Analogous results were already known, from work of
B.~Alspach, S.~C.~Locke, and D.~Witte, for the case where $X$~is cubic, or has an
odd number of vertices.
 \end{abstract}


\end{frontmatter}

\section{Introduction}

If $C$ is any cycle in a graph~$X$, then providing~$C$ with an orientation
naturally defines a flow on~$X$. (See~\S\ref{PrelimSect} for definitions and
notation used here in the introduction.) Conversely, it is well known that every
flow can be written as a sum of cycles. Brian Alspach (personal communication)
has asked which flows can be written as a sum of \emph{hamiltonian} cycles.

\begin{notation}
 Suppose $X$ is a graph. Then 
 \begin{itemize}
 \item $\flow = \flow(X)$ denotes the space of all integral flows on~$X$, that
is, the $\integer$-valued flows on~$X$, 
 \item $\even = \even(X)$ denotes the additive subgroup of~$\flow$ consisting of
the even flows, that is, the flows~$f$ such that the sum of the edge-flows of~$f$
is even, and 
 \item $\hamilton = \hamilton(X)$ denotes the additive subgroup of~$\flow$
generated by the oriented hamiltonian cycles.
 \end{itemize}
 Note that $\flow = \even$ if and only if $X$~is bipartite. On the other hand,
$\hamilton \subseteq \even$ whenever $X$ has even order.
 \end{notation}

Locke and Witte \cite{LockeWitte} showed that if $X$~is a connected Cayley graph
on a finite abelian group of odd order, then every flow on~$X$ can be expressed
as a sum of hamiltonian cycles, except for flows on one particular graph, the
cartesian product $K_3 \boxprod K_3$ of two cycles of length~$3$.

\begin{thm}[{Locke-Witte \cite[Thm.~4.1]{LockeWitte}}] \label{oddorder}
 If $X = \Cay(G;S)$~is a connected Cayley graph on a finite, abelian group~$G$ of
odd order, then $\hamilton=\flow$, unless $X \iso K_3 \boxprod K_3$, in which
case, $\flow/\hamilton \iso \integer_3$.
 \end{thm}

They also determined $\hamilton$ in the case with $X$ cubic.

\begin{obs} \label{CubicCayley}
 A connected, cubic Cayley graph on a finite, abelian group is of one of two
types:\ a M\"obius ladder, or a prism over a cycle.
 \end{obs}

\begin{thm}[{Locke-Witte \cite[Prop.~3.3]{LockeWitte}}] \label{cubic}
 ~ 
   \begin{enumerate}
 \item If $X$~is a M\"obius ladder, then
 \begin{enumerate}
 \item $\even/\hamilton \iso \integer_{n/2}$ if $X$~is bipartite, where $n$~is
the number of vertices of~$X$; 
 \item \label{cubic-Mobius-nonbip}
 $\hamilton = \even$ if $X$~is not bipartite. 
 \end{enumerate}
 \item \label{cubic-prism}
 If $X$~is a prism over a cycle of length~$n$, then
 \begin{enumerate}
 \item
 $\even/\hamilton \iso \integer_{n-1}$ if $X$~is bipartite;
 \item $\flow/\hamilton \iso \integer \oplus \integer_{n-1}$ if $X$~is not
bipartite.
 \end{enumerate}
 \end{enumerate}
 \end{thm}

We now complete the classification they began, by calculating~$\hamilton$ for all
the remaining Cayley graphs on finite abelian groups.

\begin{thm} \label{mainthm}
 If $X = \Cay(G;S)$~is a connected, non-cubic Cayley graph on a finite, abelian
group~$G$ of even order, then $\hamilton = \even$, unless
 \begin{enumerate}
 \item \label{mainthm-sqcyc}
 $X$ is the square of a cycle, in which case $\flow/\hamilton \iso
\integer_{n-2}$, where $n = |G|$; or
 \item $X$ has degree~$4$, $X$~is not bipartite, and $|G|$~is not divisible
by~$4$, in which case $\flow/\hamilton \iso \integer_{4}$, unless $X$
is the square of a cycle, in which case, \pref{mainthm-sqcyc} applies.
 \end{enumerate}
 \end{thm}

In the exceptional cases of Theorem~\ref{mainthm} (that is, in those cases with
$\hamilton \neq \even$), Propositions~\ref{SquareCycle} and~\ref{WeirdCase}
determine precisely which flows are in~$\hamilton$. The analogous results of
Locke and Witte \cite{LockeWitte} for the exceptional cases of
Theorems~\ref{oddorder} and~\ref{cubic} are recalled in Section~\ref{CubicSect}.

In the special case of bipartite graphs, the theorem can be restated as follows.

\begin{cor} \label{H=Z}
 Let $X = \Cay(G;S)$~be a connected Cayley graph on a finite, abelian group~$G$.
If $X$~is bipartite and not cubic, then every flow on~$X$ can be written
as a sum of hamiltonian cycles.
 \end{cor}

Analogous results for $\integer_2$-flows (in which the coefficients are taken
modulo~$2$) were obtained by Alspach, Locke, and Witte \cite{AlspachLockeWitte}
\see{ALWthm}.
 We rely heavily on the results of \cite{AlspachLockeWitte,LockeWitte}, and, to a
large extent, we also use the same techniques. Thus, the reader may find
it helpful to look at the proofs in those papers, especially because they provide
drawings of many of the hamiltonian cycles that appear here.

Here is an outline of the paper.
 Section~\ref{PrelimSect} presents notation and definitions. It also states our
standing assumption, which holds everywhere except here in the introduction, that
$G$ has even order.
 Section~\ref{InvolSect} presents some useful observations on involutions.
 Section~\ref{CubicSect} briefly recalls the results of \cite{LockeWitte} that
calculate~$\hamilton$ in the exceptional cases of Theorems~\pref{oddorder}
and~\pref{cubic}.
 Section~\ref{HFSect} develops the main tools to be used in an inductive proof of
our main theorem.
 Section~\ref{4cycleSect} shows that $\hamilton$ often contains certain basic
$4$-cycles.
 Section~\ref{Degree4Sect} proves that if $X$~has degree~$4$, and is not one of
the exceptional cases, then $\hamilton = \even$.
 Section~\ref{Except4Sect} treats the exceptional graphs of degree~$4$.
 Section~\ref{WeirdVanishSect} presents a somewhat lengthy proof that was omitted
from Section~\ref{Except4Sect}.
 Section~\ref{GridSect} shows, for many graphs of degree at least~$5$, that
$\hamilton$ contains all of the basic $4$-cycles.
 Section~\ref{RedundantSect} provides an induction step for the proof of the main
theorem, under the assumption that the generating set~$S$ contains a redundant
generator.
 Section~\ref{TroubleSect} deals with two cases that are not covered by our other
results.
 Section~\ref{Degree5Sect} proves that if the degree of~$X$ is at least~$5$, then
$\hamilton = \even$.
 Combining \pref{degree4}, \pref{SquareCycle}, \pref{WeirdCase}, \pref{degree7},
and the trivial observation that we have $\hamilton = \flow$ for $X$ of degree at
most~$2$ proves Theorem~\ref{mainthm}.

\begin{rmk}
 Although we discuss only integer flows, it is explained in
\cite[\S5]{LockeWitte} that these results are universal. They determine which
$A$-valued flows are linear combinations of hamiltonian cycles, for any abelian
group~$A$.
 \end{rmk}

\begin{ack}
 The authors would like to thank the Department of Mathematics and Statistics of
the University of Minnesota, Duluth for its hospitality.  Almost all of this
research was carried out during several summer visits there. The work was begun
while D.P.M.\ was a participant in an Undergraduate Research Participation
program under the direction of Joseph A.~Gallian, and all of the authors would
like to thank Joe for his encouragement and helpful suggestions.
 They would also like to thank an anonymous referee for numerous comments on the
exposition.
 Some of the authors visited the Tata Institute of Fundamental Research (Mumbai,
India) or the Department of Mathematics and Computing of the Faculty of
Education at the University of Ljubljana (Slovenia).  These authors would like to
thank their hosts for their hospitality.
 The research was partially supported by grants from the National Science
Foundation, NSERC, and the Ministry of Science of Slovenia.
 \end{ack}

\section{Preliminaries} \label{PrelimSect}

\begin{defn}
 Suppose $S$ is a subset of a finite group~$G$.
 \begin{itemize}
 \item $S$ is a \emph{symmetric generating set} for~$G$ if
 \begin{itemize}
 \item $S$~generates~$G$, that is, no proper subgroup of~$G$ contains~$S$, and
 \item we have $s^{-1} \in S$, for every $s \in S$.
 \end{itemize}
 \item The \emph{Cayley graph $\Cay(G;S)$ of~$S$ on~$G$} is the graph defined as
follows:
 \begin{itemize}
 \item the vertices of $\Cay(G;S)$ are the elements of~$G$, and
 \item there is an edge from~$g$ to~$gs$, for every $g \in G$ and $s \in S$.
 \end{itemize} 
 \end{itemize}
 \end{defn}

\begin{notation}
 Throughout this paper,
 \begin{itemize}
 \item $G$ is a finite abelian group (usually written multiplicatively);
 \item $e$ is the identity element of~$G$;
 \item $S$ is a symmetric generating set for~$G$, such that $e \notin S$; and
 \item $X = \Cay(G;S)$.
 \end{itemize}
 \end{notation}

\begin{assump} \label{|G|even}
 Throughout the remainder of this paper, 
 $$ \mbox{$|G|$ is even.} $$
 \end{assump}

\begin{notation}
 When $X = \Cay(G;S)$, and some element~$s$ of~$S$ has been chosen, we let
 \begin{itemize}
 \item $S' = S\setminus\{s,s^{-1}\}$;
 \item $G' = \langle S' \rangle$; 
 \item $X' = \Cay(G';S')$;
 \item $\flow' = \flow(X')$;
 \item $\even' = \even(X')$; and
 \item $\hamilton' = \hamilton(X')$.
 \end{itemize}
 \end{notation}

\begin{notation}
 We use 
 $$[v](t_1,t_2,\ldots,t_n) ,$$
 where $v \in G$ and $t_i \in S$, to denote the path (or cycle) in $\Cay(G;S)$
that visits the vertices 
 $$v, v t_1, v t_1 t_2, \ldots, v t_1 t_2 \cdots t_n .$$
 (When $v = e$, we usually write simply $(t_1,t_2,\ldots,t_n)$.)
 
We use a superscript to denote the concatenation of copies of the same sequence,
and the symbol $\sharp$ denotes truncation of the last term of a sequence. For
example, 
 $$\bigl( (s^2,t)^3\sharp,u \bigr) = (s,s,t,s,s,t,s,s,u) $$
 and
 $$\bigl( (s^2,t)^0,u \bigr) = (u) .$$
 Note that the notation $s^m$ in such a sequence always denotes repetitions of
the generator~$s$, not a single occurrence of the group element~$s^m$. We will
always give a new name (such as $t = s^m$) if we wish to use the group
element~$s^m$ in such a path or cycle.

The following illustrates another useful notation:
 $$ \Bigl( \bigl( s^2,t_i \bigr)_{i=1}^3, u,
 \bigl( s^2,t_i \bigr)_{i=1}^0, u \Bigr)
  = (s,s,t_1,s,s,t_2,s,s,t_3, u, u) . $$
 \end{notation}

\begin{rmk} \ 
 \begin{itemize}
 \item  We do not usually distinguish between a cycle $C = [v](s_1,\ldots,s_m)$
and the corresponding element of~$\flow$.
 \item $[vs](s^{-1})$ and $-[v](s)$ each represent the same edge as $[v](s)$, but
with the opposite orientation. Orientations serve two purposes:\ they
arise in the definition of a flow \seeDefn{FlowDefn}, and they may be
used to indicate that a path traverses a certain edge in a certain direction.
 \item If $(s_1,\ldots,s_m)$ is a cycle, then $-(s_1,\ldots,s_m) = (s_m^{-1},
\ldots, s_1^{-1})$. In particular, 
 $$-(s,t,s^{-1},t^{-1}) = (t,s,t^{-1},s^{-1}) .$$
 \end{itemize}
 \end{rmk}

\begin{sugg}
 Some hamiltonian cycles in~$X$, such as that in \eqref{E=H+2F'pf-H} on
p.~\pageref{E=H+2F'pf-H}, depend on a path $(t_i)_{i=1}^m$ in a subgraph or
quotient graph of~$X$. For simplicity, the reader may find it helpful to assume
that the path is simply~$(t^m)$, so that the subscripts can be ignored. For
example, the above-mentioned hamiltonian cycle \pref{E=H+2F'pf-H} simplifies to
 $$H = \bigl( ( s^{m-1}, t, s^{-(m-1)}, t )^{r/2},
  (t^{n-r-1}, s, t^{-(n-r-1)}, s )^{m/2} \bigr) .$$
 (In order to facilitate the simplification process, we consistently begin our
indices at~$1$, even when a different starting point would yield less complicated
formulas in the subscripts.) As soon as the simpler cycle is understood, it
should be clear that there is an analogous hamiltonian cycle that includes
subscripts. Thus, the subscripts are essentially a formality, so the correctness
of an overall proof can usually be verified without checking that the authors
have calculated the subscripts correctly.
 \end{sugg}

\begin{defn} \label{FlowDefn} \ 
 \begin{itemize}
 \item A \emph{flow} on the Cayley graph $X = \Cay(G;S)$ is a
function $f \colon G \times S \to \integer$, such that 
 \begin{itemize}
 \item $f(v,s) = - f(v s,s^{-1})$, for all $v \in G$ and $s \in S$, and
 \item $ \sum_{s \in S} f(v,s) = 0$, for all $v \in G$.
 \end{itemize}
 We usually refer to $f(v,s)$ as the \emph{edge-flow} of~$f$ on the oriented edge
$[v](s)$. (Although, for simplicity, we consider only Cayley graphs, flows can be
defined for any graph.)
 \item A \emph{weighting} of~$X$ is a function $\phi \colon G \times S \to
\integer$, such that $\phi(v,s) = - \phi(v s,s^{-1})$, for all $v \in G$ and $s
\in S$. We usually refer to $\phi(v,s)$ as the \emph{weight} of the oriented edge
$[v](s)$.
 \item Given a flow~$f$ on~$X$ and a weighting~$\phi$, the \emph{weighted sum of
the edge-flows of~$f$} is $\sum_{\alpha \in A} \phi(a) f(a)$, where $A$~is any
subset of $G \times S$, such that for each $v \in G$ and $s \in S$, the set
$A$~contains either $(v,s)$ or $(v s,s^{-1})$, but not both. It is independent of
the choice of the set~$A$.
 \end{itemize}
 \end{defn}

\begin{rmk}
 In later sections of the paper, it will sometimes be necessary to define a
particular weighting of~$X$. For convenience, whenever we specify that some
oriented edge $[v](s)$ has a certain weight~$w$, it is implicitly understood that
the oppositely oriented edge $[vs](s^{-1})$ has weight~$-w$.
 \end{rmk}

\begin{notation}
 Suppose $v \in G$, $Y$~is a subgraph of~$X$, and $f \in \flow$. We use $[v]Y$ to
denote the \emph{translate of~$Y$ by~$v$}, and $[v]f$ to denote the
\emph{translate of~$f$ by~$v$}. Namely:
 \begin{itemize}
 \item $[v]Y$ is the subgraph of~$X$ defined by:
 \begin{itemize}
 \item the vertices of $[v]Y$ are the elements of~$G$ of the form $vy$ with $y
\in Y$, and
 \item there is an edge from $v y_1$ to~$v y_2$ in~$[v]Y$ if and only if there is
an edge from~$y_1$ to~$y_2$ in~$Y$.
 \end{itemize}
 \item The edge-flow of $[v]f$ on an oriented edge $[vw](s)$ is defined to be the
same as the edge-flow of~$f$ on the oriented edge $[w](s)$, for $w \in G$ and $s
\in S$.
 \end{itemize}
 \end{notation}

\begin{defn}
 For any graphs $X$ and~$Y$, the \emph{Cartesian product $X \boxprod Y$ of $X$
and~$Y$} is the graph defined as follows:
 \begin{itemize} 
 \item the vertices of $X \boxprod Y$ are the ordered pairs $(x,y)$ with $x \in
X$ and $y \in Y$, and
 \item there is an edge from $(x_1,y_1)$ to $(x_2,y_2)$ if and only if either
 \begin{itemize}
 \item $x_1 = x_2$, and there is an edge from~$y_1$ to~$y_2$ in~$Y$, or
 \item $y_1 = y_2$, and there is an edge from~$x_1$ to~$x_2$ in~$X$.
 \end{itemize}
 \end{itemize}
 We use $X^p$ to denote the Cartesian product $X \boxprod X \boxprod \cdots
\boxprod X$ of $p$~copies of~$X$.
 \end{defn}

\begin{obs}
 If $S_i$ is a symmetric generating set for~$G_i$, for $i = 1,2$, then
 $$ \Cay(G_1;S_1) \boxprod \Cay(G_2;S_2) \iso 
 \Cay \Bigl( G_1 \times G_2; \bigl( S_1 \times \{e\} \bigr) \cup \bigl( \{e\}
\times S_2 \bigr) \Bigr) .$$
 \end{obs}

\begin{defn} \ 
 \begin{itemize}
 \item A \emph{basic $4$-cycle} in~$X$ is any $4$-cycle of the form
$[v](s,t,s^{-1},t^{-1})$ with $v \in G$ and $s,t \in S$.
 \item An element~$s$ of~$G$ is an \emph{involution} if $s$~is of order~$2$; that
is, if $s^2 = e$ and $s \neq e$.
 \item An element~$s$ of~$S$ is a \emph{redundant} generator if $\langle S'
\rangle = G$.
 \item The generating set~$S$ is \emph{irredundant} if none of its elements are
redundant.
 \item For a fixed element~$s$ of~$S$, an edge of~$X$ is an \emph{$s$-edge} if it
is of the form $[v](s)$ or $[v](s^{-1})$, for some $v \in G$.
 \end{itemize}
 \end{defn}

\begin{notation} \ 
 \begin{itemize}
 \item $\integer_n$ denotes the additive group of integers modulo~$n$.
 \item $K_n$ denotes the complete graph on $n$~vertices.
 \item $C_n$ denotes the cycle of length~$n$. (With the expectation that it will
not cause confusion, $C_1$, $C_2$ and~$C_3$ are used to denote certain more
general cycles in the proof of Corollary~\ref{OddHeight4Cycles}.)
 \item $|g|$ denotes the order of the element~$g$ of~$G$.
 \item $|H|$ denotes the order of the subgroup~$H$ of~$G$  (that is, the number
of elements of~$H$).
 \item $\langle g \rangle$ denotes the subgroup of~$G$ generated by the
element~$g$.
 \item $\langle A \rangle$ denotes the subgroup of~$G$ generated by the subset~$A$
of~$G$.
 \end{itemize}
 \end{notation}

\begin{defn} \ 
 \begin{itemize}
 \item For any graph~$Y$, we call $K_2 \boxprod Y$ the \emph{prism over $Y$}.
 \item A \emph{M\"obius ladder} is a graph isomorphic to $\Cay \bigl(
\integer_{2n}; \{\pm1, n\} \bigr)$, for some natural number~$n$.
 \item $X$ is \emph{the square of an even cycle} if there exist $s$ and~$t$
in~$S$, such that
 \begin{itemize}
 \item $S = \{s^{\pm1},t^{\pm1}\}$,
 \item $t = s^2$, and
 \item $t^2 \neq e$.
 \end{itemize}
 (The final condition is a convention: we do not consider the cubic graph $K_4$
to be the square of an even cycle.)
 It is not difficult to show that if $X$~is isomorphic to the square of an even
cycle, then $X$~itself is the square of an even cycle.
 \end{itemize}
 \end{defn}

\begin{warn}
 If we write $S = \{s^{\pm1},t^{\pm1}\}$, then it is obvious that $|S| \le 4$.
However, it need not be the case that $|S|$ is exactly~$4$, unless additional
restrictions are explicitly imposed. For example, it could be the case that $s =
s^{-1}$, or that $s = t$.
 \end{warn}

\begin{warn}
 We use $p$~and~$q$ to denote arbitrary integers; they are \emph{not} assumed
to be prime numbers.
 \end{warn}

\section{Remarks on involutions in~$S$} \label{InvolSect}

\begin{obs} \label{noirredinvols}
 If $S$ is an irredundant generating set, or, more generally, if no involution
in~$S$ is a redundant generator, then we may assume that $S$~contains no more
than one involution. To see this, let
 \begin{itemize}
 \item $S_1$ be the set of involutions in~$S$, 
 \item $G_1$ be the subgroup generated by~$S_1$, and 
 \item $G_2$ be the subgroup generated by $S \setminus S_1$. 
 \end{itemize}
 Then $G_1 \cap G_2$ is trivial (because none of the elements of~$S_1$ are
redundant), so $G = G_1 \times G_2$. Hence 
 $$ X \iso \Cay( G_1 ; S_1) \boxprod \Cay( G_2; S \setminus S_1 ) .$$
 Now, the desired conclusion follows by noting that
 $ \Cay( G_1 ; S_1)$ is isomorphic to either $(C_4)^p$ (if $|S_1| = 2p$) or $K_2
\boxprod (C_4)^p$ (if $|S_1| = 2p+1$), for some natural number~$p$.
 \end{obs}

\begin{obs} \label{noinvols}
 If $|S| = 4$, and 
 $X$~is not the prism over a M\"obius ladder, 
 then we may assume that $S$~does not contain any involutions. Specifically:
 \begin{enumerate}
 \item \label{noinvols-K2xprism}
 If $X$ is isomorphic to the cartesian product of~$K_2$ with a prism over a
cycle~$C_n$, then \pref{noirredinvols} applies.
 \item If $X$ is obtained from the prism over a cycle of length~$2n$ by adding
the diagonals, that is, if
 $$X \iso \Cay \bigl( \integer_2 \times \integer_{2n}; \{ (1,0), (1,n),
(0,1)^{\pm1} \} \bigr) ,$$
 then
 $$X \iso 
 \begin{cases}
 \Cay \bigl( \integer_2 \times \integer_{2n}; \{(1,n/2)^{\pm1}, (0,1)^{\pm1}\}
\bigr) 
 & \mbox{if $n$~is even}, \\
 \Cay \bigl( \integer_4 \times \integer_n; \{(1,0)^{\pm1}, (2,1)^{\pm1}\}
\bigr) 
 & \mbox{if $n$~is odd}. \\
 \end{cases}
 $$
 \end{enumerate}
 It is not difficult to see that these cases are exhaustive, given the list of
cubic graphs in Observation~\ref{CubicCayley}.
 \end{obs}

\begin{lem} \label{notallinvols}
 If there exists $s \in S$, such that
 \begin{itemize}
 \item $s$ is a redundant involution in~$S$,
 and
 \item $S'$ is irredundant,
 \end{itemize}
 then there is a generating set~$T'$ for an abelian group~$H$, and an
involution~$t$ in~$H$, such that
 \begin{enumerate}
 \item $X \iso \Cay \bigl( H; T' \cup \{t\} \bigr)$,
 \item $X' \iso \Cay(H;T')$, 
 \item $t$ is an involution, and
 \item not every element of~$T'$ is an involution.
 \end{enumerate}
 \end{lem}

\begin{proof}
 We may assume that every element of~$S'$ is an involution. (Otherwise, take $H =
G$, $T' = S'$, and $t = s$.) Write $s = t_1 + t_2 + \cdots + t_p$, where
$t_1,\ldots,t_p$ are distinct elements of~$S'$, with $p \ge 2$. Let 
 \begin{align*}
 H &= \integer_4 \times \bigl\langle S' \setminus \{t_1,t_2\} \bigr\rangle 
 \subseteq \integer_4 \times G'
 ,\\
 T' &= \{(1,e)^{\pm1}\} \cup \Bigl( \{0\} \times \bigl( S' \setminus
\{t_1,t_2\} \bigr) \Bigr)
 \subseteq H
 , \\
 t &= (2, t_3 + \cdots + t_p) \in H
 .
 \end{align*}
 Then it is not difficult to verify the desired conclusions.
 \end{proof}

\section{Graphs that are cubic or of odd order} \label{CubicSect}

Let us recall the observations of \cite{LockeWitte} that describe exactly which
flows are in~$\hamilton'$, for the cases where Theorem~\ref{oddorder}
or~\ref{cubic} gives an imprecise answer. For completeness, we include all of
these results, even though the proofs in later sections require only
\pref{H+C=F(K3xK3)}, \pref{BipMobius}, and~\pref{EvenPrism}. 

We state the following result (and its corollary) only for~$G'$, because of our
standing assumption \pref{|G|even} that $|G|$ is even.

\begin{lem}[{\cite[pf.~of Prop.~3.1]{LockeWitte}}] \label{K3xK3}
 Suppose
 \begin{itemize}
 \item $G' \iso \integer_3 \times \integer_3$, and
 \item $S' = \{t^{\pm1}, u^{\pm1}\}$.
 \end{itemize}
 Give 
 \begin{itemize}
 \item weight~$1$ to the oriented $t$-edge $[v](t)$, for each $v \in \langle t
\rangle$,
 \item weight~$-1$ to the oriented $t$-edge $[v](t)$, for each $v \in u \langle t
\rangle$,
 and
 \item weight~$0$ to each of the other edges of~$X'$.
 \end{itemize}
 Then a flow is in~$\hamilton'$ if and only if the weighted sum of its edge-flows
is divisible by~$3$.
 \end{lem}

In the situation of Lemma~\ref{K3xK3}, it is easy to see that the weighted sum of
the edge-flows of any basic $4$-cycle is nonzero, so the following observation is
an easy consequence (cf.\ proof of \pref{deg4+4cyc}).

\begin{cor} \label{H+C=F(K3xK3)}
  If
 \begin{itemize}
 \item $G' \iso \integer_3 \times \integer_3$,
 \item $|S'| = 4$,
 and
 \item $\hamilton$ contains some basic $4$-cycle of~$X'$,
 \end{itemize}
 then $\flow' \subseteq \hamilton + \hamilton'$.
 \end{cor}

\begin{lem}[{\cite[pf.~of Prop.~3.3(1a)]{LockeWitte}}]
 \label{BipMobius}
 Suppose $X$ is a M\"obius ladder, and $X$ is bipartite, so we may write
 \begin{itemize}
 \item $S = \{t^{\pm1},u\}$,
 \item $G = \langle t \rangle$,
 \item $|t| = |G| = 2n$, where $n$~is odd, and
 \item $u = t^n$ is an involution.
 \end{itemize}
 Give
 \begin{itemize}
 \item weight $(-1)^i$ to the oriented $u$-edge $[t^i](u)$, for each~$i$, and
 \item weight $0$ to each $t$-edge.
 \end{itemize}
 Then a flow is in~$\hamilton$ if and only if the weighted sum of its edge-flows
is divisible by~$n$.
 \end{lem}

\begin{lem}[{\cite[pf.~of Prop.~3.3(2a)]{LockeWitte}}]
 \label{EvenPrism}
 Suppose $X$ is the prism over a cycle, and $X$ is bipartite, so we may write
 \begin{itemize}
 \item $S = \{t^{\pm1},u\}$,
 \item $G = \langle t \rangle \times \langle u \rangle$, and
 \item $|G| = 2|t| = 2n$, where $n$~is even.
 \end{itemize}
 Give
 \begin{itemize}
 \item weight $(-1)^j$ to the oriented $t$-edge $[t^i u^j](t)$, for each~$i$
and~$j$, and
 \item weight $0$ to each $u$-edge.
 \end{itemize}
 Then a flow is in~$\hamilton$ if and only if the weighted sum of its edge-flows
is divisible by~$n-1$.
 \end{lem}

\begin{lem}[{\cite[pf.~of Prop.~3.3(2b)]{LockeWitte}}]
 Suppose $X$ is the prism over a cycle, and $X$ is \emph{not} bipartite, so
we may write
 \begin{itemize}
 \item $S = \{t^{\pm1},u\}$,
 \item $G = \langle t \rangle \times \langle u \rangle$,
 and
 \item $|G| = 2|t| = 2n$, where $n$~is odd.
 \end{itemize}
 Give
 \begin{itemize}
 \item weight $1$ to the oriented $t$-edge $[v](t)$, for each $v \in \langle t
\rangle$, and
 \item weight $0$ to all of the other edges of~$X$.
 \end{itemize}
 Then a flow is in~$\hamilton$ if and only if 
 \begin{enumerate}
 \item the weighted sum of its edge-flows is divisible by~$n-1$, and
 \item the flow on the oriented edge $(t)$ is the negative of the flow on the
oriented edge $[u](t)$.
 \end{enumerate}
 \end{lem}

\section{Relations among $\even$, $\even'$, $\hamilton$, $\hamilton'$, and $2
\flow$} \label{HFSect}

 In most cases, our goal in this paper is to show $\even \subseteq \hamilton$,  and
our proof proceeds by induction on $|S|$. Thus, we usually know that $\even'
\subseteq \hamilton'$, and we wish to show that $\even \subseteq \hamilton$. This
section presents some of our main tools to accomplish this. They are of three
general types:
 \begin{enumerate}
 \item It suffices to show $2 \flow \subseteq \hamilton$: \pref{ALWthm}.
 \item \label{HFSect-EinH+E'}
 Results that show $\even \subseteq \hamilton + \even'$: \pref{E=H+E'(G'=G)},
\pref{E=H+2F'}, \pref{E=H+E'(G'<>G)}.
 \item \label{HFSect-H'inH}
 Results that show $\hamilton' \subseteq \hamilton$ (or, in some cases, show only
that $2\hamilton' \subseteq \hamilton$): \pref{G=G'->H'inH}, \pref{H'(G'even)},
\pref{2H'inH(G'odd)notinvol}, \pref{2H'inH(G'odd)deg5}.
 \end{enumerate}
 Note that if $\even' \subseteq \hamilton'$, then combining \pref{HFSect-EinH+E'}
with the strong form of~\pref{HFSect-H'inH} yields the desired conclusion $\even
\subseteq \hamilton$.

 Most of the results in this section assume that $\hamilton$ contains certain
basic $4$-cycles; results in Sections~\ref{4cycleSect} and~\ref{GridSect} show
that $\hamilton$ often contains every basic $4$-cycle.

\begin{thm}[{Alspach-Locke-Witte \cite[Thm.~2.1]{AlspachLockeWitte}}]
\label{ALWthm}
 Let $X$~be a connected Cayley graph on a finite abelian group.
 If $X$ is not a prism over an odd cycle, then $\even \subseteq \hamilton + 2\flow$.
 \end{thm}

We state the following result for $X'$, rather than~$X$, because it applies to
all groups, including those of odd order.

\begin{thm}[Chen-Quimpo \cite{ChenQuimpo}] \label{ChenQuimpoThm}
 Suppose $|S'| \ge 3$, and let $v$~and~$w$ be any two distinct vertices of~$X'$.
 \begin{enumerate}
 \item If $X'$ is not bipartite, then there is a hamiltonian path from~$v$ to~$w$.
 \item If $X'$ is bipartite, then either 
 \begin{itemize}
 \item there is a hamiltonian path from~$v$ to~$w$, or
 \item there is a path of even length from~$v$ to~$w$.
 \end{itemize}
 \end{enumerate}
 \end{thm}

\begin{cor}[{cf.~pf.~of \cite[Cor.~3.2]{AlspachLockeWitte}}] \label{E=H+E'(G'=G)}
  Suppose $s \in S$, such that
 \begin{itemize}
 \item $G' = G$, 
 \item $|S'|\ge 3$, and 
 \item either $X$~is bipartite or $X'$~is not bipartite.
 \end{itemize}
 Then $\even \subseteq \hamilton + \even'$.
 \end{cor}

\begin{lem} \label{weight1}
 Let $s \in S$, and assume
 \begin{itemize}
 \item $|S| \ge 4$,
 \item $G' \neq G$, and 
 \item $(s,t,s^{-1},t^{-1}) \in \hamilton$, for every $t \in S'$.
 \end{itemize}
 Give 
 \begin{itemize}
 \item weight~$1$ to the oriented $s$-edge $[v](s)$ if $v \in s^{-1} G'$, and 
 \item weight~$0$ to all of the other edges of~$X$.
 \end{itemize}
 If $k$~is the weighted sum of the edge-flows of some element~$H$ of~$\hamilton$,
then
 $k\flow \subseteq \hamilton + \flow'$.
 \end{lem}

\begin{proof}
 Let $f \in \flow$. We wish to show $k f \in \hamilton + \flow'$. By adding an
appropriate multiple of~$H$ to~$k f$, we obtain a flow~$f_1$, such that the
weighted edge-sum of~$f_1$ is~$0$. 

For each $i$, with $-1 \le i \le |G/G'| - 2$, let
 $$ E_i = \bigl\{\, [v](s) \mid v \in s^i G' \bigr\} .$$
 Then, for $0 \le i \le |G/G'| - 2$, the union $-E_{-1} \cup E_i$ is the set of
oriented edges that start in $\cup_{j=0}^i s^j G'$, and end in the complement. So
the net flow of~$f_1$ through the edges in~$E_i$ must equal the net flow through
the edges in $E_{-1}$, which is~$0$. 

Therefore, by adding appropriate multiples of basic $4$-cycles of the form
 $$ [v](s,t,s,t^{-1}) ,$$
 with $v \in s^m G'$ and $t \in S'$, to~$f_1$, we obtain a flow that does not use
any edges of~$E_m$. Repeating this for all~$m$ (including $m = -1$), we obtain a
flow~$f_2$ that does not use any $s$-edges.

So $f_2$ is a sum of flows on various cosets of~$G'$. The following claim shows
that $f_2 \in \hamilton + \flow'$, so we conclude that $k f \in \hamilton + f_2
\subseteq \hamilton + \flow'$, as desired.

\begin{proofclaim}
 For any $f' \in \flow'$, and any $v \in G$, we have $[v]f' \in \hamilton +
\flow'$.
 \end{proofclaim}
 We may assume 
 \begin{itemize}
 \item $f'$ is a cycle $[w](t_1,\ldots,t_m)$, with each $t_j \in S'$, and $w \in
G'$, and
 \item $v = s^r$, for some $r > 0$. 
 \end{itemize}
 Then
 $$[v]f' - f' = \sum_{i=1}^{r} \sum_{j=1}^m [s^{i-1} w t_1 t_2 \ldots t_{j-1}]
(s,t_j, s^{-1}, t_j^{-1}) ,$$
 is a sum of basic $4$-cycles, so $[v]f' - f' \in \hamilton$.
 \end{proof}

\begin{cor}[{cf.~\cite[Lem.~3.8]{AlspachLockeWitte}}] \label{E=H+2F'}
 Suppose $s \in S$, and we have
 \begin{itemize}
 \item $|S| \ge 4$,
 \item $G' \neq G$,
 \item $(s,t,s^{-1}, t^{-1}) \in \hamilton$, for every $t \in S'$,
 and
 \item either
 \begin{enumerate}
 \renewcommand{\theenumi}{\alph{enumi}}
 \item \label{E=H+2F'-Xbip}
 $X$ is bipartite, or
 \item \label{E=H+2F'-X'notbip}
 $X'$ is not bipartite.
 \end{enumerate}
 \end{itemize}
 Then $\even \subseteq \hamilton + 2\flow'$.
 \end{cor}

\begin{proof}
 It suffices to show $\flow \subseteq \hamilton + \flow'$, for then multiplying
by~$2$ yields $2\flow \subseteq \hamilton + 2\flow'$, and then the desired
conclusion follows from Theorem~\ref{ALWthm}.

 Let $m = |G/G'|$. Note that, because $|S| \ge 4$, we have $|G'| \ge 3$.

\setcounter{proofcase}{0}

\begin{proofcase}
 Assume $X$ is bipartite.
 \end{proofcase}
 Let $H' = (t_1,t_2,\ldots,t_n)$ be a hamiltonian cycle in~$X'$. There is some
$r$ with $s^{-m} = t_1 t_2 \cdots t_r$. (Note that, because $X$ is bipartite, we
know that $r+m$ is even.) If $m$ is even, define 
 \begin{equation} \label{E=H+2F'pf-H}
 H = \Bigl( \bigl( s^{m-1}, t_{2i-1}, s^{-(m-1)}, t_{2i} \bigr)_{i=1}^{r/2},
 \bigl( (t_{r+i})_{i=1}^{n-r-1}, s, (t_{n-i}^{-1})_{i=1}^{n-r-1}, s \bigr)^{m/2}
\Bigr);
 \end{equation}
  whereas, if $m$ is odd, let
 \begin{align*}
  H &= \Bigl( \bigl( s^{m-2}, t_{2i-1}, s^{-(m-2)}, t_{2i} \bigr)_{i=1}^{(r+1)/2},
 \\
 & \hskip1in
 \bigl( (t_{r+1+i})_{i=1}^{n-r-2}, s, (t_{n-i}^{-1})_{i=1}^{n-r-2}, s
\bigr)^{(m-1)/2},
 \\
 & \hskip1in
 (t_{r+1+i})_{i=1}^{n-r-1}, (t_i)_{i=1}^{r}, s \Bigr) .
 \end{align*}
 Then Lemma~\ref{weight1} (with $k = 1$) implies $\flow \subseteq \hamilton +
\flow'$, as desired.

\begin{proofcase}
 Assume $X'$ is not bipartite, and $|S'| \ge 3$.
 \end{proofcase}
 Choose nonidentity elements $g_1,g_2,\ldots,g_m$ of~$G'$, such that $g_1 g_2
\cdots g_m = s^{-m}$.  Theorem~\ref{ChenQuimpoThm} implies, for each~$i$, that
there is a hamiltonian path $(t_{i,j})_{j=1}^{n-1}$ in~$X'$ from~$e$ to~$g_i$.
Define the hamiltonian cycle
 $H = \bigl( (t_{i,j})_{j=1}^{n-1}, s \bigr)_{i=1}^{m}$.
 Then Lemma~\ref{weight1} (with $k = 1$) implies $\flow \subseteq \hamilton +
\flow'$, as desired.

\begin{proofcase}
 Assume $X'$ is not bipartite, and $|S'| = 2$.
 \end{proofcase}
 Let $t \in S'$, so $S = \{s^{\pm1}, t^{\pm1}\}$. Because $X'$ is not bipartite,
we know $|t|$ is odd. We have $s^{-m} = t^r$ for some~$r$ with $0 \le r \le
|t|-1$. We may assume $r$~is even, by replacing $t$~with its inverse if
necessary. Define $H$ as in \eqref{E=H+2F'pf-H}, with $t_i = t$ for all~$i$.
Then Lemma~\ref{weight1} (with $k = 1$) implies $\flow \subseteq \hamilton +
\flow'$, as desired.
 \end{proof}

\begin{cor}[{cf.~\cite[Lem.~3.8]{AlspachLockeWitte}}] \label{E=H+E'(G'<>G)}
 Suppose $s \in S$, and we have
 \begin{itemize}
 \item $|S| \ge 4$,
 \item $G' \neq G$,
 and
 \item $(s,t,s^{-1}, t^{-1}) \in \hamilton$, for every $t \in S'$.
 \end{itemize}
 Then $\even \subseteq \hamilton + \even'$.
 \end{cor}

\begin{proof}
 It suffices to show $2\flow \subseteq \hamilton + \even'$, for then the desired
conclusion follows from Theorem~\ref{ALWthm}. We may assume that $X'$~is
bipartite, but $X$~is not bipartite, for otherwise Corollary~\ref{E=H+2F'}
applies (and yields a stronger conclusion).

 Let $m = |G/G'|$. Because $|S| \ge 4$, we have $|G'| \ge 3$.
 Let $(t_1,t_2,\ldots,t_n)$ be a hamiltonian cycle in~$X'$. 

We will construct a hamiltonian cycle~$H$ in~$X$, such that $H$ contains the
oriented edges $[s^{m-1}](s)$ and $[s^{m-1}t_1](s)$, but no other oriented edges
of the form $\pm[v](s)$ with $v \in s^{-1} G'$. Then Lemma~\ref{weight1}
(with $k = 2$) implies $2 \flow \subseteq \hamilton + \flow'$. Since $X'$ is
bipartite, this means $2\flow \subseteq \hamilton + \even'$, as desired.

There are $p$~and~$q$
with $0 \le p,q < n$, such that
 $s^m = t_1 t_2 \ldots t_p$ and $s^m t_1 = t_1 t_2 \ldots t_q$. 
 Because $X'$~is bipartite, but $X$~is~not bipartite, it is not difficult to see
that the cycle $(t_1,t_2,\ldots,t_p, s^{-m})$ must be odd. Similarly, the cycle 
$(t_2,t_3,\ldots,t_q, s^{-m})$ must also be odd. Therefore
 \begin{itemize}
 \item $p$~and~$m$ have opposite parity, and
 \item $q$~and~$m$ have the same parity.
 \end{itemize}
 Furthermore, we may assume the hamiltonian cycle $(t_1,t_2,\ldots,t_n)$ has been
chosen so that
 \begin{equation} \label{minimalp}
 \mbox{$p$ is as small as possible.}
 \end{equation}
 Note that if $p>q$ and $q \neq 0$, then replacing the hamiltonian cycle
$(t_1,t_2,\ldots,t_n)$ with $(t_2,t_3,\ldots,t_n, t_1)$ replaces~$p$ with~$q-1$.
Because $q-1 < q < p$, this contradicts \pref{minimalp}. Hence
 \begin{itemize}
 \item either $p < q$ or $q = 0$.
 \end{itemize}

 If $m$~is odd, then $q \neq 0$ (because $q$ must be odd), so we must have $p <
q$. Define
 \begin{align*}
  H &= \Bigl( s, (t_i)_{i=1}^{p},
(t_{p+2i-1},s^{-1},t_{p+2i},s)_{i=1}^{(q-p-1)/2},
 \\
 &\hskip1in
 (t_{q-1+i})_{i=1}^{n-q-1}, 
 \bigl( s, (t_{n-1-i}^{-1})_{i=1}^{n-2}, s, (t_i)_{i=1}^{n-2} \bigr)^{(m-3)/2},
 \\
 &\hskip1in
 s, (t_{n-1-i}^{-1})_{i=1}^{n-3}, s, (t_{q+i})_{i=1}^{n-q-1}, s^{m-1}, t_n, s,
(t_{p+1-i}^{-1})_{i=1}^{p} \Bigr) .
 \end{align*}
 If $m$ is even and $p < q$, define
 \begin{align*}
 H &= \Bigl( s^m, \bigl( s, t_{p+2-2i}^{-1}, s^{-1}, t_{p+1-2i}^{-1}
\bigr)_{i=1}^{(p-1)/2}, s,
 \\
 &\hskip1in
 \bigl( s, (t_{i+1})_{i=1}^{n-2}, s, (t_{n-i}^{-1})_{i=1}^{n-2} \bigr)^{(m-2)/2},
 s, (t_{q+1-i}^{-1})_{i=1}^{q-p-1}, s,
 \\
 &\hskip1in
 (t_{p+1+i})_{i=1}^{q-p}, s^{-1},
 \bigl( t_{q+2i}, s, t_{q+2i+1}, s^{-1} \bigr)_{i=1}^{(n-q-2)/2}, t_n \Bigr).
 \end{align*}
 If $m$ is even and $p > q = 0$, then, by considering the hamiltonian cycle
$(t_1^{-1}, t_n^{-1}, t_{n-1}^{-1}, \ldots, t_2^{-1})$, we see that $p = 1$.
Define
 \begin{equation*}
  H = \Bigl(s^m, \bigl((t_{i+1})_{i=1}^{n-2}, s, 
    (t_{n-i}^{-1})_{i=1}^{n-2}, s\bigr)^{m/2}\Bigr) .
 \end{equation*}
 This completes the proof.
 \end{proof} 

At one point (namely, in Subcase~\ref{degree4pf-bip-notcirc} of the proof of
Proposition~\ref{degree4}), we will need the following more general (but weaker)
version of \pref{E=H+E'(G'<>G)}. Note that this result assumes only that
$\hamilton$ contains the doubles of the basic $4$-cycles, not that it contains the
$4$-cycles themselves.

\begin{cor} \label{2E=H+2E'(G'<>G)}
 Suppose $s \in S$, and we have
 \begin{itemize}
 \item $|S| \ge 4$,
 \item $G' \neq G$,
 and
 \item $2(s,t,s^{-1}, t^{-1}) \in \hamilton$, for every $t \in S'$.
 \end{itemize}
 Then $2\even \subseteq \hamilton + 2\even'$.
 \end{cor}

\begin{proof}
 Let $\basics$ be the subgroup of~$\flow$ that is generated by the basic
$4$-cycles of the form $(s,t,s^{-1},t^{-1})$ with $t \in S'$.
 The proof of Corollary~\ref{E=H+E'(G'<>G)} shows that $\even \subseteq \hamilton +
\basics + \even'$. Therefore $2 \even \subseteq 2 \hamilton + 2 \basics + 2
\even'$. By assumption, we have $2 \basics \subseteq \hamilton$ (and it is obvious
that $2 \hamilton \subseteq \hamilton$), so we conclude that $2\even \subseteq
\hamilton + 2\even'$, as desired.
 \end{proof}

\begin{obs} \label{G=G'->H'inH}
 If $G' = G$, then $X'$ is a spanning subgraph of~$X$, so $\hamilton' \subseteq
\hamilton$. (In particular, if $G' = G$ and $\even \subseteq \hamilton +
\hamilton'$, then $\even \subseteq \hamilton$.)
 \end{obs}

\begin{lem}[{cf.~\cite[Lem.~3.6]{AlspachLockeWitte}}] \label{H'(G'even)}
  Suppose $s \in S$, such that
 $$\mbox{$|G'| \ge 4$ is even.}$$
 \begin{enumerate}
 \item \label{H'(G'even)-1}
 If $(s,t,s^{-1},t^{-1}) \in \hamilton$, for every $t \in S'$, then $\hamilton'
\subseteq \hamilton$.
 \item \label{H'(G'even)-2}
 If $2 (s,t,s^{-1},t^{-1}) \in \hamilton$, for every $t \in S'$, then
$2\hamilton' \subseteq \hamilton$.
 \end{enumerate}
 \end{lem}

\begin{proof}
 We prove only~\pref{H'(G'even)-1}, because \pref{H'(G'even)-2} is very similar.
Let $m = |G/G'|$. Given a hamiltonian cycle $H' = (t_1,t_2,\ldots,t_n)$ in~$X'$,
we have a corresponding hamiltonian cycle
 $$H = \bigl( s^{m-1},t_{2i-1},s^{-(m-1)},t_{2i} \bigr)_{i=1}^{n/2}$$
 in~$X$.
 Because 
 $$H - H' = \sum_{i=0}^{m-2} \sum_{j=1}^{n/2}[s^i t_1 t_2 \ldots t_{2j-2}] \bigl(
s,t_{2j-1},s^{-1},t_{2j-1}^{-1} \bigr) $$
  is a sum of basic $4$-cycles, we conclude that $H' \in \hamilton$.
 \end{proof}

\begin{lem} \label{2H'inH(G'odd)notinvol}
 Suppose $s \in S$, such that
 \begin{itemize}
 \item $|s| > 2$,
 \item $|G'| \ge 3$ is odd, and 
 \item $(s,t,s^{-1},t^{-1}) \in \hamilton$, for every $t \in S'$. 
 \end{itemize}
 Then $2 \hamilton' \subseteq \hamilton$. 
 \end{lem}

\begin{proof}
 Let $H' = (t_1,t_2,\ldots,t_n)$ be a hamiltonian cycle in~$X'$. Let $m =
|G/G'|$. There is some $r$ with $s^{-m} = t_1 t_2 \cdots t_r$. 
 Because $n$~is odd, we may assume $r$ is even (by replacing $H'$ with $-H' =
(t_n^{-1}, t_{n-1}^{-1}, \ldots, t_{1}^{-1})$ if necessary). Note that, because
$|s| > 2$, we have $m > 2$ if $r = 0$.
 Define 
 \begin{align*}
 \epsilon &= 
 \begin{cases}
 0  &\mbox{if $r=0$;}\\
 1 &\mbox{otherwise}
 ,
 \end{cases}
 \\
 H_1 &=
 \Bigl( \bigl( t_{n-i+1}^{-1} \bigr)_{i=1}^{n-r-1+2\epsilon},
 \bigl( s^{m-1},t_{r-2i+1}^{-1}, s^{-(m-1)}, t_{r-2i}^{-1}
\bigr)_{i=1}^{\epsilon(r-2)/2},
 \\
 &\hskip1in
 s, \bigl( s^{m-2} \bigr)^\epsilon, t_1^{-1}, \bigl( s^{-(m-2)} \bigr)^\epsilon,
 \\
 &\hskip1in
 \bigl( t_{n-2i+2}^{-1},s^{m-2},t_{n-2i+1}^{-1},s^{-(m-2)}
\bigr)_{i=1}^{(n-r-3+2\epsilon)/2},
 \\
 &\hskip1in
 t_{r+3-2\epsilon}^{-1}, s^{m-2}, \bigl( t_2^{-1},s^{-(m-3)},t_1^{-1},s^{m-3}
\bigr)^{1-\epsilon}, s \Bigr)
 \\
 \intertext{and}
 H_2 
 &= \Bigl( \bigl( s^{m-1},t_{2i-1},s^{-(m-1)}, t_{2i} \bigr)_{i=1}^{r/2},
 \\
 &\hskip1in
  \bigl( (t_{r+i})_{i=1}^{n-1-r}, s, (t_{n-i}^{-1})_{i=1}^{n-1-r}, s
\bigr)^{m/2} \Bigr) 
 .
 \end{align*}
 Then $H_1 - H_2 + 2 H'$ is a sum of basic $4$-cycles, so $2H'$ belongs
to~$\hamilton$.

Perhaps we should elaborate further. Let $Y$ be the spanning subgraph of~$X$ with
 \begin{align*}
  E(Y) &= \bigset{
 [s^i t_1 t_2 \cdots t_j] (t_{j+1})
 }{
 \begin{matrix} 
 0 \le i \le m-1, \\
 0 \le j \le n-2
 \end{matrix}
 }
 \\ &\hskip1in
 \cup
 \bigset{
 [s^i t_1 t_2 \cdots t_j] (s)
 }{
 \begin{matrix} 
 0 \le i \le m-2, \\
 0 \le j \le n-1
 \end{matrix}
 }
 .
 \end{align*}
 Then 
 $$ H_1 - H_2 + H' + [s] H'
 \in \flow(Y)
 ,$$
 and $Y$~is naturally isomorphic to the cartesian product of the two paths
$(s)^{m-1}$ and $(t_1,t_2,\ldots,t_{n-1})$,
 so any flow on~$Y$ is a sum of basic $4$-cycles. Furthermore,
 $$ [s] H' - H'
 = \sum_{i=1}^n
 [t_1 t_2 \cdots t_{i-1}] (s,t_i, s^{-1}, t_i^{-1} ) $$
 is a sum of basic $4$-cycles. Therefore
 $$ H_1 - H_2 + 2 H'
  = \bigl( H_1 - H_2 + H' + [s] H' \bigr)
 - \bigl( [s] H' - H' \bigr)
 $$
 is a sum of basic $4$-cycles, as claimed.
 \end{proof}

\begin{lem} \label{2H'inH(G'odd)deg5}
 Suppose $S=\{s,t^{\pm 1},u^{\pm 1}\}$, such that
 \begin{itemize}
 \item $|S| = 5$,
 \item $|G'|$ is odd, 
 \item $s^2 = e$, 
 \item $G' \not= \langle t \rangle$,
 \item $G' \not= \langle u \rangle$,
 and 
 \item every basic $4$-cycle is in~$\hamilton$.
 \end{itemize}
 Then $2 \hamilton' \subseteq \hamilton$. 
 \end{lem}

\begin{proof}
 Let $p = |G'|/|t|$ and $n = |t|$. Note that $p$~and~$n$ must be odd, since
$|G'|$~is odd. There is some $r$ with $u^{-p} = t^r$.  Because $n$~is odd, we may
assume $r$~is even (by replacing $u$ with~$u^{-1}$ if necessary).

Define
 \begin{align*}
 H_1 &= \Bigl( t^ {n-1},u^{p-2},t^{-(n-3)},u,t^{n-1},s,t^{n-1},u^{-(p-1)},t,
 \\ & \qquad
 \bigl( t^{n-2},u,t^{-(n-2)},u \bigr)^{(p-1)/2}\sharp,
 s,t^{-1},
 \\ & \qquad
 \bigl( u^{-1},t^{p-2},u^{-1},t^{-(p-2)} \bigr)^{(p-3)/2},
 u^{-1} \Bigr).
 \end{align*}
 Then $H_1 - 2(t^n)$ is a sum of basic 4-cycles. (For example, this follows from
the observation that  
 $$H_1 - [u^{p-1}](t^n) - [s u^{p-1}](t^n)$$
 belongs to $\flow(Y)$, where $Y$~is a spanning subgraph of~$X$ that is naturally
isomorphic to the cartesian product of the paths $(s)$, $(t^{m-1})$,
and~$(u^{p-1})$; cf.\ proof of \pref{2H'inH(G'odd)notinvol}). Therefore
 \begin{equation} \label{2H'inH(G'odd)deg5pf-2tn}
 2(t^n) \in \hamilton
 .
 \end{equation}

 Now define
 \begin{align*}
 H_2 &= \Bigl( s,u^{-1},(t^{-r},u^{-1},t^{r},u^ {-1})^{(p-1)/2}\sharp,
 \\ & \qquad
 \bigl( t,u^{p-2},t,u^{2-p} \bigr)^{(n-r-1)/2}, u^{-1},t^{2-n}, s,
 \\ & \qquad
 \bigl( t^{n-2},u,t^{2-n},u \bigr)^{(p-1)/2},t^{n-1},u^{1-p} \Bigr).
 \end{align*}
 By adding certain basic $4$-cycles involving~$s$ to~$H_2$, we can obtain an
element of~$\hamilton$ that is in $X'$ and uses only one edge of the form
$[t^i]u^{-1}$. We can now apply the proof of Lemma~\ref{weight1} to~$X'$, with
$u$~taking the role of~$s$, to obtain the conclusion that $\flow' \subseteq
\hamilton + (t^n)$, since the basic $4$-cycles are in~$\hamilton$ and not just
in~$\hamilton'$. In particular, this tells us that $2\hamilton' \subseteq
\hamilton + 2(t^n)$. Combining this with \pref{2H'inH(G'odd)deg5pf-2tn}, we see
that $2\hamilton' \subseteq \hamilton$, as desired.
 \end{proof}

\section{Some basic $4$-cycles in~$\hamilton$} \label{4cycleSect}

In this section, we show, for $s,t \in S$, that $\hamilton$ often contains the
flow $2(s,t,s^{-1},t^{-1})$ \see{4cyc=} and, if $|G|$ is divisible by~$4$, the
basic $4$-cycle $(s,t,s^{-1},t^{-1})$ \see{4cyc}. Our main tool is the
construction described in Lemma~\ref{4c}, which was already used in
\cite{AlspachLockeWitte,LockeWitte} (and goes back to \cite{Marusic}).

\begin{lem} \label{4c}
 Suppose 
 \begin{itemize}
 \item $x,y,z \in S$,
 \item $v,w \in G$, and 
 \item $H_+$ and~$H_-$  
 are oriented hamiltonian cycles in~$X$.  
 \end{itemize}
 Then:
 \begin{enumerate}
    \item \label{4c-s}
       If $H_+$ contains both the oriented path $[v](x,y,x^{-1})$ and the
       oriented edge $[v x z](y)$, then 
       $$ (x,y,x^{-1},y^{-1}) + [x](z,y,z^{-1},y^{-1}) \in \hamilton .$$
    \item \label{4c-d}
       If $H_-$ contains both the oriented path $[w](x,y,x^{-1})$ and the
       oriented edge $[w x y z](y^{-1})$, then 
      $$ (x,y,x^{-1},y^{-1}) - [x](z,y,z^{-1},y^{-1}) \in \hamilton .$$
    \item \label{4c-2}
       If both \pref{4c-s} and~\pref{4c-d} apply, then
       $2 (x,y,x^{-1},y^{-1}) \in \hamilton$.
    \item \label{4c-change}
       If
       \begin{itemize}
         \item $(x,y,x^{-1},y^{-1}) \in \hamilton$, and 
         \item either \pref{4c-s} or \pref{4c-d} applies, 
       \end{itemize}
       then $(z,y,z^{-1},y^{-1}) \in \hamilton$.
 \end{enumerate}
 \end{lem}

\begin{proof}
 We may assume $v = w = e$.

 \pref{4c-s} Construct a hamiltonian cycle~$H_+'$ by replacing 
 \begin{itemize}
 \item the path $(x,y,x^{-1})$ with the edge~$(y)$ and 
 \item the edge $[x z](y)$ with the path $[x z](z^{-1},y,z)$. 
 \end{itemize}
 Then $H_+ - H_+'$ is the sum of the two given $4$-cycles.

 \pref{4c-d} Construct a hamiltonian cycle~$H_-'$ by replacing 
 \begin{itemize}
 \item the path $(x,y,x^{-1})$ with the edge~$(y)$ and 
 \item the edge $[x y z](y^{-1})$ with the path $[x y z](z^{-1},y^{-1},z)$. 
 \end{itemize}
 Then $H_- - H_-'$ is the difference of the two given $4$-cycles.

 \pref{4c-2} Adding the flows from \pref{4c-s}
and~\pref{4c-d} results in $2 (x,y,x^{-1},y^{-1})$.

 \pref{4c-change}  The difference of $(x,y,x^{-1},y^{-1})$ and the flow that
results from either of~\pref{4c-d} and~\pref{4c-s} is
$\pm[x](z,y,z^{-1},y^{-1})$.
 \end{proof}

\begin{prop} \label{4cyc=}
 If $s,t \in S$, such that
 \begin{itemize}
 \item $|S| = 4$,
 \item $|s| \ge 3$,
 \item $t \notin \langle s \rangle$, and
 \item $t^2 \notin \{e,s^{\pm1}\}$,
 \end{itemize}
 then, letting $Q = (s,t,s^{-1}, t^{-1})$, we have
 $$Q \equiv -Q \equiv [v]Q \pmod{\hamilton} ,$$
 for all $v \in G$. In particular, $2Q \in \hamilton$.
 \end{prop}

\begin{proof}
 It suffices to show
 \begin{enumerate}
 \renewcommand{\theenumi}{\alph{enumi}}
 \item \label{4cyc=pf-sdiff}
 $Q - [s]Q \in \hamilton$, 
 \item \label{4cyc=pf-tdiff}
 $Q - [t]Q \in \hamilton$, and
 \item \label{4cyc=pf-sum}
 either $Q + [s]Q \in \hamilton$, or $Q + [t]Q \in \hamilton$.
 \end{enumerate}

 Let $n = |G|/|s|$, and write $t^n = s^r$, with $0 \le r < |s|$. Because $t
\notin \langle s \rangle$, we have $G \neq \langle s \rangle$, so $n \ge 2$. Note
that $n |s| = |G|$ is even, so $n$ and~$|s|$ cannot both be odd. If $|s|$~is even
and $n \ge 3$, define, for future reference, the hamiltonian cycle
 \begin{equation} \label{4cyc=pf-H*}
 H_* = \Bigl( t^{n-1}, s^{|s|-1}, t^{-1}, s^{-(|s|-2)}, t^{-(n-2)},
 \bigl( s, t^{n-3}, s, t^{-(n-3)} \bigr)^{(|s|-2)/2}, s \Bigr) 
 .
 \end{equation}

Now let us begin by establishing~\pref{4cyc=pf-sdiff}. Let
 $$ H_{\mathrm{a}} =
 \begin{cases}
 \Bigl( t^{n-1}, s, \bigl( s^{|s|-2}, t^{-1}, s^{-(|s|-2)}, t^{-1}
\bigr)^{n/2}\sharp, s^{-1} \Bigr)
 & \mbox{if $|s|$ is odd} 
 ,\\
 \hfil
 \bigl( t, s^{|s|-1}, t^{-1}, s^{-(|s|-1)} \bigr)
 & \mbox{if $n = 2$}
 ,\\
 \hfil \mbox{$H_*$ as in \eqref{4cyc=pf-H*}}
 & \mbox{otherwise} 
 .
 \end{cases}
 $$
 Then $H_{\mathrm{a}}$ contains both the oriented path $[s^{-2} t^{n-1}] (s, t^{-1},
s^{-1})$ and the oriented edge $[t^{n-2}] (t)$, so Lemma~\fullref{4c}{d}
(with $x = s$, $y = t^{-1}$, $z = s$, and $w = s^{-2} t^{n-1}$) establishes
\pref{4cyc=pf-sdiff}. 

All that remains is to establish \pref{4cyc=pf-tdiff} and~\pref{4cyc=pf-sum}.

\setcounter{proofcase}{0}

\begin{proofcase}
 Assume $n \ge 3$.
 \end{proofcase}

\begin{subcase} 
 Assume $|s|$ is odd.
 \end{subcase}
 Note that $n$ must be even (so $n \ge 4$), because $n |s| = |G|$ is even. We may
assume $t^n \neq e$ (so $r \neq 0$), for otherwise, by interchanging
$s$~and~$t$, we may transfer to one of the cases where $|s|$~is even and $n$~is
odd. We may assume $r$ is odd, by replacing $s$ with its inverse if necessary.
Define hamiltonian cycles
 \begin{align*}
 H_- &= \Bigl( \bigl( t^{n-3}, s,t^{-(n-3)}, s \bigr)^{(|s|-r)/2}, 
 s^{r-1}, t,  
 \\
 &\hskip1in
 \bigl( t^{n-2}, s^{-1} ,t^{-(n-2)} s^{-1} \bigr)^{(r-1)/2}, t^{n-3},
s^{-(|s|-r)}, t, s^{|s|-r}, t \Bigr) 
 \end{align*}
 and
 $$ H_+ = \Bigl( \bigl(t^{n-1}, s, t^{-(n-1)}, s \bigr)^{(|s|-r)/2},
 \bigl( s^{r-1}, t, s^{-(r-1)}, t \bigr)^{n/2} \bigr) \Bigr) .$$
 Then
 \begin{itemize}
 \item $H_-$ contains both the oriented path $[t^{n-4}](t,s,t^{-1})$ and the
oriented edge $[t^{n-2}s](s^{-1})$, so Lemma~\fullref{4c}{d} (with $x =
t$, $y = s$, $z = t$, and $v = t^{n-4}$) establishes~\pref{4cyc=pf-tdiff}.
 \item $H_+$ contains the oriented path $[t^{n-2}](t,s,t^{-1})$ and the oriented
edge $[t^n](s)$ (because $r$~is odd), so Lemma~\fullref{4c}{s} (with $x =
t$, $y = s$, $z = t$, and $v = t^{n-2}$) establishes \pref{4cyc=pf-sum}.
 \end{itemize}

\begin{subcase}
 Assume $|s|$ is even and $n \ge 4$.
 \end{subcase}
 Define hamiltonian cycles
 \begin{equation} \label{4cyc=pf-Hminus}
 H_- = \Bigl( \bigl( s, t^{n-2},s, t^{-(n-2)} \bigr)^{(|s|-2)/2},
 s, t^{n-1}, s^{-(|s|-1)}, t^{-(n-1)} \Bigr) 
 \end{equation}
 and
 $$ H_+ = 
 \begin{cases}
 \hfil \mbox{$H_*$ as in \eqref{4cyc=pf-H*}}
 & \mbox{if $r \neq 2$} \\
 \hfil [s^2] H_*
 & \mbox{if $r = 2$} 
 .
 \end{cases}
 $$
 Then:
 \begin{itemize}
 \item $H_-$ contains both the oriented path $[s t^{n-3}](t,s,t^{-1})$ and the
oriented edge $[s^2 t^{n-1}](s^{-1})$, so Lemma~\fullref{4c}{d} (with $x =
t$, $y = s$, $z = t$, and $w = s t^{n-3}$) establishes~\pref{4cyc=pf-tdiff}.
 \item $H_+$ contains both the oriented path $[s t](t^{-1},s,t)$ and the oriented
edge $[s t^{-1}](s)$, so Lemma~\fullref{4c}{d} (with $x = t^{-1}$, $y =
s$, $z = t^{-1}$, and $v = s t$) establishes~\pref{4cyc=pf-sum}.
 \end{itemize}

\begin{subcase}
 Assume $|s|$ is even and $n = 3$.
 \end{subcase}
 We may assume $0 \le r \le |s|/2$, by replacing~$t$ with $t^{-1}$ if necessary.
 Define hamiltonian cycles
 $$ \mbox{$H_-$ as in \eqref{4cyc=pf-Hminus}},$$
 and
 $$ H_+ = \bigl( t^2, s^{|s|-2-r}, t,
 s^{-(|s|-3)}, t, s^{|s|-2}, t, s^{-r}, t, s \bigr) .$$
 Then:
 \begin{itemize}
 \item $H_-$ establishes~\pref{4cyc=pf-tdiff}, exactly as in the previous subcase.
 \item $H_+$ contains both the oriented path $[s^2](s^{-1},t,s)$ and the oriented
edge $(t)$, so Lemma~\fullref{4c}{d} (with $x = s^{-1}$, $y = t$, $z =
s^{-1}$, and $v = s^2$) establishes~\pref{4cyc=pf-sum}.
 \end{itemize}

 \begin{proofcase}
 Assume $n=2$.
 \end{proofcase}
 We have $2 \le r \le |s|-2$, because $t^2 \notin \{e,s^{\pm1}\}$. Notice that
this implies $|s|\ge 4$.

\begin{subcase}
  Assume $|s|$~is even, and $r$~is odd.
 \end{subcase}
 We may assume $r \le |s|/2$, by replacing $s$ with $s^{-1}$ if necessary. 

  If $r = |s|/2$, then $|t| = 4$, so $|G|/|t| = |s|/2 = r \ge 3$ (because $r$ is
odd). Thus, an earlier subcase applies, after interchanging $s$ and~$t$.

 We may now assume $r < |s|/2$. Define
 $$ H_+ = \Bigl( \bigl( t,s,t^{-1},s \bigr)^{(|s|-2r)/2}, \bigl(
s^{r-1},t,s^{-(r-1)},t \bigr)^2 \Bigr) .$$
 Then $H_+$ contains both the oriented path $[s^{-2}](s,t,s^{-1})$ and the
oriented edge $(t)$, so Lemma~\fullref{4c}{s} (with $x = s$, $y = t$, $z =
s$, and $v = s^{-2}$) implies that $Q + [s]Q \in \hamilton$. Therefore, because
 $$ Q = 
 \bigl( t, s^{|s|-1}, t^{-1}, s^{-(|s|-1)} \bigr)
 + \sum_{i=1}^{(|s|-2)/2} [s^{2i-1}] \bigl( Q + [s]Q \bigr) ,$$
 we conclude that $Q \in \hamilton$, which makes
\pref{4cyc=pf-tdiff} and \pref{4cyc=pf-sum} trivial.

\begin{subcase}
 Assume that either
 \begin{itemize}
 \item $|s|$ is odd or
 \item $|s|$ and $r$ are even. 
 \end{itemize}
 \end{subcase}
 We may assume $|s| - r$ is even, by replacing $t$ with its inverse if necessary.
Define 
 $$H_- = \bigl( t,s,t^{-2},s^{r-2},t^{-1},s^{-(|s|-3)},t,s^{|s|-r-2},t \bigr).$$ 
 Then $H_-$ contains both the oriented path $(t,s,t^{-1})$ and the oriented edge
$[st^2](s^{-1})$, so Lemma~\fullref{4c}{d} (with $v=e$, $x=t$, $y=s$ and
$z=t$) tells us that $(-Q) - [t](-Q) \in \hamilton$. This
establishes~\pref{4cyc=pf-tdiff}.

Define the hamiltonian cycle
 \begin{equation*}
 H_+ = \begin{cases}
 \Bigl( t,s,t^{-1},s^{|s|-2},t,s^{-(|s|-3)},t \Bigr)
 &\mbox{if $r=|s|-2$}\\
 \Bigl(t,s,t^{-1},s^{r+1},t^{-1},s^r,\bigl( s,t^{-1},s,t
\bigr)^{(|s|-r-2)/2}\sharp \Bigr) &\mbox{if $r<|s|-2$}
 .
 \end{cases}
 \end{equation*}
 Then $H_+$ contains both the oriented path $(t,s,t^{-1})$ and the oriented edge
$[t^2](s)$, so Lemma~\fullref{4c}{s} (with $v=e$, $x=t$, $y=s$ and $z=t$)
tells us that $(-Q)+[t](-Q) \in \hamilton$. This establishes \pref{4cyc=pf-sum}.
 \end{proof}

\begin{cor} \label{4cyc}
 Suppose $s,t \in S$, such that
 \begin{itemize}
 \item $|G|$ is divisible by~$4$, 
 \item $|s| \ge 3$, 
 \item $t \notin \langle s \rangle$, and
 \item either
 \begin{enumerate}
 \renewcommand{\theenumi}{\alph{enumi}}
 \item \label{4cyc-|G/G'|>2}
 $|G|/|s| \ge 3$, or
 \item \label{4cyc-t2<>s}
 $t^2 \notin \{e,s^{\pm1}\}$.
 \end{enumerate}
 \end{itemize} 
 Then $(s,t,s^{-1},t^{-1}) \in \hamilton$.
 \end{cor}

\begin{proof}

\setcounter{proofcase}{0}

\begin{proofcase} \label{4cycpf-deg4}
 Assume $\langle s, t \rangle = G$ and $t^2 \notin \{e, s^{\pm1}\}$.
 \end{proofcase}

\begin{itemize}
 \item If $|s|$ is even, let $x = s$, $y = t$, $m = |s|$, and $n = |G|/|s|$.
 \item If $|s|$ is odd, let $x = t$, $y = s$, $m = |G|/|s|$, and $n = |s|$.
 \end{itemize}
 In either case,
 $$H = \Bigl( x^{m-1}, y, \bigl( y^{n-2}, x^{-1}, y^{-(n-2)}, x^{-1} \bigr)^{m/2}
\sharp, y^{-1} \Bigr) $$
 is a hamiltonian cycle in~$X$. Letting 
 $$Q = (x,y,x^{-1},y^{-1}) = \pm (s,t,s^{-1},t^{-1}) ,$$
 we have
 $$H = \left( \sum_{i=1}^{m-1} [x^{i-1}]Q \right)
 + \left( \sum_{i=1}^{m/2} \sum_{j=1}^{n-2} [x^{2i-2} y^j] Q \right) , $$
 so $H$ is the sum of 
 $$ (m-1) + \frac{m}{2} (n-2) = \frac{mn}{2} - 1$$
 translates of~$Q$. Because $mn = |G|$ is divisible by~$4$, we know that $(mn/2)
- 1$ is odd. Therefore, Proposition~\ref{4cyc=} implies 
 $$H \equiv \bigl( (mn/2) - 1 \bigr) Q \equiv Q \pmod{\hamilton} .$$
 Since $H \in \hamilton$, we conclude that $Q \in \hamilton$.
 
\begin{proofcase}  \label{4cycpf-highdeg}
 Assume that either $\langle s,t \rangle \neq G$ or $t^2 \in \{e, s^{\pm1}\}$.
 \end{proofcase}
 We show how to reduce to the previous case.

First, let us show 
 $|G|/|s| \ge 3$.
 By hypothesis, if this fails to hold, then $t^2 \notin \{e, s^{\pm1}\}$, so we
may assume that the latter holds. Hence, the assumption of this case implies
$\langle s,t \rangle \neq G$. Since $t \notin \langle s \rangle$, we conclude that
 $$ \frac{|G|}{|s|} =  \frac{|G|}{|\langle s,t \rangle|} \cdot \frac{|\langle
s,t \rangle|}{|s|} \ge 2 \cdot 2 > 3 ,$$
 as claimed.

 Let $(t_1,\ldots,t_n)$ be a hamiltonian cycle in 
 $$\Cay \bigl( \quot; S \setminus \{s,s^{-1}\} \bigr) ,$$
 with $t_1 = t$ \cf{ChenQuimpoThm}. Define permutations $\sigma$~and~$\tau$
of~$G$ by
 \begin{itemize}
 \item $x^\sigma = xs$ and
 \item $x^\tau = x t_i$, where $i = i(x)$ satisfies $1 \le i \le n$ and $x \in
\langle s \rangle t_1 t_2 \cdots t_{i-1}$. 
 \end{itemize}
 Let $Y$ be the spanning subgraph of~$X$ whose edge set is
 $$ E(Y) = \bigl\{\, ( x, x^\sigma ) \mid x \in V(G) \,\bigr\} \cup
 \bigl\{\, ( x, x^\tau ) \mid x \in V(G) \,\bigl\} .$$
 Then $Y$ is a connected, spanning subgraph of~$X$. It is not difficult to
see that $\sigma$~and~$\tau$ generate a transitive, abelian group~$\Gamma$ of
automorphisms of~$Y$ (and any transitive, abelian permutation group is regular),
so $Y$ is isomorphic to the Cayley graph 
 $$Y^* = \Cay\bigl( \Gamma; \{\sigma^{\pm1},\tau^{\pm1}\} \bigr) .$$
 Furthermore, the natural isomorphism carries the $4$-cycle $(s,t,s^{-1},t^{-1})$
to $(\sigma,\tau, \sigma^{-1},\tau^{-1})$. 

Note that $|\sigma| = |s|$ and $|\Gamma|/|\sigma| = |G|/|s| \ge 3$, so $\tau^2
\notin \langle \sigma \rangle$.  From Case~\ref{4cycpf-deg4}, we know that
 $(\sigma,\tau,\sigma^{-1}, \tau^{-1}) \in \hamilton(Y^*)$. Hence, via the
isomorphism, we see that
 $(s,t,s^{-1},t^{-1}) \in \hamilton(Y) \subseteq \hamilton(X)$.
 \end{proof}

\begin{rmk}
 The assumption that $t \notin \{e, s^{\pm1}\}$ is necessary in
Proposition~\ref{4cyc=} and Corollary~\ref{4cyc}, as is seen from
Theorem~\fullref{cubic}{prism} and Proposition~\ref{SquareCycle}. If $|S| = 4$
and $X$~is not bipartite, then the assumption that $|G|$ is divisible by~$4$ is
necessary in Corollary~\ref{4cyc}, as is seen from Proposition~\ref{WeirdCase}.
 \end{rmk}

\section{The graphs of degree 4 with $\hamilton = \even$} \label{Degree4Sect}

In this section, we show that $\hamilton = \even$ in many cases where $|S| = 4$
\see{degree4}. In Section~\ref{Except4Sect}, we will calculate $\hamilton$ in the
cases not covered by this result.

\begin{prop} \label{degree4}
 If 
 \begin{itemize}
 \item $|S| = 4$, and 
 \item $X$ is \emph{not} the square of an even cycle, and
 \item either
 \begin{enumerate}
 \renewcommand{\theenumi}{\alph{enumi}}
 \item \label{degree4-bip}
 $X$~is bipartite, or
 \item \label{degree4-4}
 $|G|$ is divisible by~$4$,
 \end{enumerate} 
 \end{itemize}
 then $\hamilton = \even$.
 \end{prop}

We preface the proof with an observation on bipartite graphs and with the
treatment of a special case.

\begin{lem} \label{Z2bipartite}
 If
 \begin{itemize}
 \item $|G|$ is not divisible by~$4$, and
 \item every element of~$S$ has even order,
 \end{itemize}
 then $X$ is bipartite.
 \end{lem}

\begin{proof}
 Let $H = \{\,g^2 \mid g \in G\,\}$. Because $|G| \equiv 2 \pmod{4}$, we see that
$|H| = |G|/2$, and $|H|$~is odd.  We know that no element of~$S$ belongs to~$H$
(because the elements of~$H$ have odd order), so the subgraph of~$X$ induced by
each of the two cosets of~$H$ has no edges. Therefore, the coset decomposition $G
= H \cup Hg$ is a bipartition of~$G$, so $G$~is bipartite.
 \end{proof}

 \begin{lem} \label{Mobius2layers}
 If $X \iso K_2 \boxprod Y$, where $Y$~is a M\"obius ladder, then $\hamilton =
\even$.
  \end{lem}
 
 \begin{proof}
  We may assume 
 $$X = \Cay(\integer _2 \times \integer _{2n}; \{ s, t^{\pm1}, u \}) ,$$
 where $s =(0,n)$, $t = (0,1)$, and $u = (1,0)$.

\setcounter{step}{0}

\begin{step} \label{Mobius2layersPf-utut}
 We have $(u,t,u,t^{-1}) \in \hamilton$.
 \end{step}
 Define the hamiltonian cycle
  $$ H = \bigl( (s,t)^{n} \sharp, u,(s,t^{-1})^{n} \sharp,u \bigr). $$
 Then the sum $H + [s]H$ has edge-flow~$0$ on each $s$-edge, so $H + [s]H \in
\flow'$. Under the weighting of~$X'$ specified in Lemma~\ref{EvenPrism}, with~$2n$
in the role of~$n$, the weighted sum of the edge-flows of $H + [s]H$ is $4(n-1)$,
which is relatively prime to~$2n-1$. Thus, we conclude from Lemma~\ref{EvenPrism}
that $H + [s]H$ generates $\flow'/\hamilton'$, so $\flow' \subseteq \hamilton +
\hamilton'$. Because $X'$ is a spanning subgraph of~$X$, we have $\hamilton'
\subseteq \hamilton$, so this implies $\flow' \subseteq \hamilton$. Therefore
$(u,t,u,t^{-1}) \in \hamilton$, as desired.

\begin{step} \label{Mobius2layersPf-(susu)inH}
 We have $(s,u,s,u) \in \hamilton$.
 \end{step}
  The hamiltonian cycle
  $$ \bigl( u, t^{-(n-1)}, s, t^{n-1} \bigr)^2 $$
 contains both the oriented path $[t^{-1}](t,u,t^{-1})$ and the oriented edge
$[su](u)$, so Lemma~\fullref{4c}{d} (with $x = t$, $y = u$, $z = s$,
and $w = t^{-1}$) implies
 $$ (t,u,t^{-1},u) - [t](s,u,s,u) \in \hamilton .$$
 From Step~\ref{Mobius2layersPf-utut}, we know $(t,u,t^{-1},u) \in \hamilton$, so
we conclude that $(s,u,s,u)$ also belongs to~$\hamilton$.

\begin{step} \label{Mobius2layersPf-E(Y)inH}
 For
 $$Y = \Cay \bigl( \{0\} \times \integer _{2n}; \{ s,t^{\pm1} \} \bigr) \subseteq
X ,$$
 we have $\even(Y) \subseteq \hamilton$.
 \end{step}
 From Step~\ref{Mobius2layersPf-(susu)inH} and Lemma~\fullref{H'(G'even)}{1}
(with $u$ in the role of~$s$), we see that $\hamilton(Y) \subseteq \hamilton$.
Thus, 
 $$ \mbox{it suffices to show $\even(Y) \subseteq \hamilton + \hamilton(Y)$.} $$
 We may assume $n$ is odd, for otherwise Theorem~\fullref{cubic}{Mobius-nonbip}
implies $\even(Y) = \hamilton(Y) \subseteq \hamilton + \hamilton(Y)$. Consider the
hamiltonian cycle
 $$H' = \bigl(s,t^{n-2},u,t^{-(2n-3)},u,t^{n-2},s,u,t,u\bigr).$$
 We have
 $$ (s,t^{-1},s,t) 
 = H' - (t^{-1},u,t,u) - \sum_{i=1}^{2n-3} [t^i](t,u,t^{-1},u) \in \hamilton
.$$
 Under the weighting of~$Y$ specified in  Lemma~\ref{BipMobius}, with $s$ in the
role of~$u$, the weighted sum of the edge-flows of $(s,t^{-1},s,t)$ is~$\pm2$,
which is relatively prime to~$n$. Thus, Lemma~\ref{BipMobius} implies that
$(s,t^{-1},s,t)$ generates $\flow(Y)/\hamilton(Y)$, so we conclude that $\even(Y)
\subseteq \hamilton + \hamilton(Y)$.

\begin{step}
 Completion of the proof.
 \end{step}
 Given any even flow $f \in \even$, we wish to show $f \in \hamilton$. Adding
appropriate $4$-cycles of the forms  $[v](u,s,u,s)$ and $[v](u,t,u,t^{-1})$
eliminates all edges of $[u]Y$ from~$f$, and, hence, all $u$-edges as well,
leaving an even flow $f_1 \in \even(Y)$. From Step~\ref{Mobius2layersPf-E(Y)inH},
we know that $\even(Y) \subseteq \hamilton$, so we have $f_1 \in \hamilton$.
Hence, $f \in \hamilton$, as desired.
 \end{proof}

\begin{proof}[{\bf Proof of Proposition~\ref{degree4}}]
 By Remark~\ref{noinvols} and Lemma~\ref{Mobius2layers}, we may assume that $S$
has no involutions. Let $S= \{s^{\pm 1},t^{\pm 1} \}$. 
 
\setcounter{proofcase}{0}

 \begin{proofcase} \label{degree4pf-bip}
 Assume $X$ is bipartite.
 \end{proofcase}
 In this case, we know $|t|$ is even, and $t^2 \neq s^{\pm 1}$.  

\begin{subcase}
 Assume $t \in \langle s \rangle$.
 \end{subcase}
 Write $t = s^r$. We may assume $2 \le r < |s|/2$, by replacing~$t$ with its
inverse if necessary. Because $X$ is bipartite, we know that $r$~is odd. Give 
 \begin{itemize} 
 \item weight $0$ to each $s$-edge, and
 \item weight $(-1)^i$ to each oriented $t$-edge $[s^i](t)$.
 \end{itemize}
 Then the weighted sum of the edge-flows of the hamiltonian cycle
 $$ H_1 = \Bigl( \bigl( t, s, t^{-1}, s \bigr)^{(r-1)/2} , t ,  s^{|s|-2r+1}
\Bigr) $$
 is~$r$, and the weighted sum of the edge-flows of the hamiltonian cycle
 $$ H_2 = \Bigl( t, \bigl( t, s^{-1}, t^{-1}, s^{-1} \bigr)^{(r-1)/2},  t^2, 
s^{|s|-2r-1} \Bigr) $$
 is~$2-r$. 

Given any flow~$f$ on~$X$, we wish to show $f \in \hamilton$. Because $r$ is
relatively prime to~$2-r$, some integral linear combination of~$H_1$ and~$H_2$
has the same weighted edge-sum as $f$; thus, by subtracting this linear
combination, we may assume the weighted edge-sum of~$f$ is~$0$. Then, by
subtracting a linear combination of hamiltonian cycles of the form
 $$ [v] \bigl( t, s^{-(r-1)}, t, s^{|s|-r-1} \bigr) , $$
 we may assume that $f$ does not use any $t$-edges. Then $f$ is a multiple of the
hamiltonian cycle $(s)^{|s|}$, so $f \in \hamilton$.

\begin{subcase} \label{degree4pf-bip-notcirc}
 Assume $t \not\in \langle s \rangle$.
 \end{subcase}
 By Lemma~\ref{4cyc=}, we know $2(s,t,s^{-1},t^{-1}) \in \hamilton$. 
 We have
 \begin{align*}
 \even &\subseteq \hamilton + 2 \flow && \mbox{\see{ALWthm}} \\
 &\subseteq \hamilton + 2 \even && \mbox{($X$ is bipartite, so $\flow = \even$)} \\
 &\subseteq \hamilton + 2 \even' && \mbox{(see~\ref{2E=H+2E'(G'<>G)}, with the
roles of $s$ and~$t$ interchanged)} \\
 &\subseteq \hamilton + 2 \hamilton' && \mbox{($X'$ is an even cycle, so $\even' =
\hamilton'$)} \\
 &\subseteq \hamilton && \mbox{\fullsee{H'(G'even)}{2}}
 ,
 \end{align*}
 as desired.

\begin{proofcase}
 Assume $X$ is not bipartite, so $|G|$ is divisible by~$4$. 
 \end{proofcase}
 Let $m = |G|/|t|$, and write $s^m = t^r$, for some~$r$, with $0 \le r < |t|$.

\begin{subcase}
 Assume $\langle t \rangle \neq G$.
 \end{subcase}
 Because $X$ is not the square of an even cycle, we know $s^2 \notin
\{t^{\pm1}\}$. Therefore Corollary~\fullref{4cyc}{t2<>s} (with the roles of $s$
and~$t$ interchanged) implies that the $4$-cycle $(s,t,s^{-1},t^{-1})$ is
in~$\hamilton$.

If $|t|$ is even, then
 \begin{align*}
 \even &\subseteq \hamilton + \even' && \mbox{\see{E=H+E'(G'<>G)}} \\
 &\subseteq \hamilton + \hamilton' && \mbox{($X'$ is a cycle, so $\flow' =
\hamilton'$)} \\
 &\subseteq \hamilton && \mbox{\fullsee{H'(G'even)}{1}} 
 ,
 \end{align*}
 as desired.

If $|t|$ is odd, then
 \begin{align*}
 \even &\subseteq \hamilton + 2\flow' && \mbox{\fullsee{E=H+2F'}{X'notbip}} \\
 &\subseteq \hamilton + 2\hamilton' && \mbox{($X'$ is a cycle, so $\flow' =
\hamilton'$)} \\
 &\subseteq \hamilton && \mbox{\see{2H'inH(G'odd)notinvol}} 
 , 
 \end{align*}
 as desired.

\begin{subcase}
 Assume $\langle t \rangle = G$.
 \end{subcase}
 Since $X$ is not bipartite (and $|t|=|G|$ is even), $r$~must be even.  So
$\langle s \rangle \subseteq \langle t^2 \rangle \neq G$. Thus, by interchanging
$s$ and~$t$, we can move out of this subcase.
 \end{proof}

For future reference, let us record the following special case of the
proposition. (Note that no bipartite graph is the square of an even cycle.)

\begin{cor} \label{deg4bip}
 If $|S| = 4$ and $X$ is bipartite, then $\hamilton = \even$.
 \end{cor}

\section{The graphs of degree 4 with $\hamilton \neq \even$} \label{Except4Sect}

In this section, we provide an explicit description of~$\hamilton$ for the graphs
of degree~$4$ that are not covered by Proposition~\ref{degree4}
\seeand{SquareCycle}{WeirdCase}.  We also establish two corollaries that will be
used in the study of graphs of higher degree \seeand{deg4+4cyc}{weird+2tn}.

\begin{prop} \label{SquareCycle}
 Suppose $X$ is the square of an even cycle, so 
 $$ \mbox{$S = \{s^{\pm1},t^{\pm1}\}$ with $t = s^2$}.$$
 Give
 \begin{itemize}
 \item weight~$0$ to each $s$-edge, and 
 \item weight $(-1)^i$ to each $t$-edge $[s^i](t)$. 
 \end{itemize}
 Then a flow is in~$\hamilton$ if and only if the weighted sum of its edge-flows
is divisible by $|G|-2$.
 \end{prop}

\begin{proof}
 Let $n = |G|/2$.
 All hamiltonian cycles are of one of the following two forms:
 $$H_1 = \pm [v] \bigl( t^{n-1}, s, t^{-(n-1)}, s^{-1} \bigr) $$
 or
 \begin{align*}
 H_2 &= \pm [v] \bigl (s^{n_0}, (t,s^{-1},t), s^{n_1}, (t,s^{-1},t), s^{n_2},
 \\ & \hskip1in
 (t,s^{-1},t),  \ldots,
  s^{n_{k-1}}, (t,s^{-1},t), s^{n_k} \bigr)
 \end{align*}
 (for some $k \ge 0$ and $n_0,\ldots,n_k \ge 0$ with $3k + \sum n_i = 2n$). In
both cases, it is easy to see that the weighted sum of the edge-flows is
divisible by $2n-2$.

Conversely, given any flow~$f$ such that  the weighted sum of its edge-flows
is~$m(2n-2)$, for some integer~$m$, we wish to show $f \in \hamilton$.  The
weighted sum of the edge-flows of $f - m H_1$ is~$0$, so, by adding appropriate
multiples of hamiltonian cycles of the form~$H_2$ (with $k = 1$ and $n_0 = 0$), we
obtain a flow~$f'$ that does not use any $t$-edges. Then $f'$ is a multiple of the
hamiltonian cycle $(s^{2n})$, so $f' \in \hamilton$. Therefore $f \in \hamilton$.
 \end{proof}

\begin{prop} \label{WeirdCase}
 Suppose 
 \begin{itemize}
 \item $|S| = 4$,
 \item $X$ is not bipartite, 
 \item $|G|$ is not divisible by~$4$, and
 \item $X$ is not the square of an even cycle,
 \end{itemize}
 so, by Lemma~\ref{Z2bipartite},
 \begin{itemize}
 \item $S = \{t^{\pm1},u^{\pm1}\}$, where
 \item $t$~has odd order, and 
 \item $t \neq u^{\pm 2}$.
 \end{itemize}
 Give
 \begin{itemize}
 \item weight~$(-1)^j$ to each oriented $t$-edge $[t^i u^j](t)$, and
 \item weight~$0$ to each~$u$-edge.
 \end{itemize}
 Then a flow belongs to~$\hamilton$ if and only if the weighted sum of its
edge-flows is divisible by~$4$. 
 \end{prop}

This proposition is obtained by combining Lemmas~\ref{WeirdVanish}
and~\ref{WeirdGetAll}.

\begin{obs} \label{WeirdEven<>Weight}
 In the situation of Proposition~\ref{WeirdCase}, we know that $|G|/|t|$ is even,
so it is not difficult to see that $\sum_{v \in G} f\bigl( [v](u) \bigr)$ is
even, for all $f \in \flow$. Therefore, a flow on~$X$ is even if and only if the
weighted sum of its edge-flows is even.
 \end{obs}

\begin{rmk}
 Additionally, some even flows (such as any
basic $4$-cycle) have weight~$2$, so the result implies that $\hamilton \neq
\even$. In fact, $\even/\hamilton \iso \integer_2$.
 \end{rmk}

\begin{lem} \label{WeirdVanish}
  If 
 \begin{itemize}
 \item $S$ and $G$ are as described as in Proposition~\ref{WeirdCase}, and
 \item $H$ is any hamiltonian cycle in~$X$,
 \end{itemize}
 then the weighted sum of the edge-flows of~$H$ is divisible by~$4$.
 \end{lem}

\begin{proof}
 Because it is rather lengthy, and involves arguments of a different sort from
those in the  rest of the paper, this proof has been postponed to a section of its
own \seeSect{WeirdVanishSect}.  The reader can easily verify that this proof does
not rely on any of the subsequent results in the present section.
 \end{proof}

\begin{lem} \label{WeirdGetAll}
 If 
 \begin{itemize}
 \item $S$ and $G$ are as described as in Proposition~\ref{WeirdCase}, and
 \item $f$ is a flow, such that the weighted sum of the edge-flows of~$f$ is
divisible by~$4$, 
 \end{itemize}
 then $f \in \hamilton$.
 \end{lem}

\begin{proof}
 Let $\basics$ be the subgroup of~$\flow$ generated by the basic $4$-cycles, and
let $s = u$.

\setcounter{step}{0}

\begin{step} \label{WeirdGetAllpf-FinH+Q+F}
 We have $\even \subseteq \hamilton + \basics + 2\flow'$. 
 \end{step}
  Obviously, $\hamilton + \basics$ contains every basic $4$-cycle, so this
follows from (the proof of) Corollary~\fullref{E=H+2F'}{X'notbip}. (Since we used
only the fact that $\hamilton$ contained all hamiltonian cycles and all basic
$4$-cycles, we may replace $\hamilton$ there by $\hamilton + \basics$.)

\begin{step} \label{WeirdCasePf-f=4cyc}
 We may assume $f \in \basics$.
 \end{step}
 We have
 \begin{align*}
 \even &\subseteq \hamilton + \basics + 2\flow' && \mbox{(from
Step~\ref{WeirdGetAllpf-FinH+Q+F})} \\
 &= \hamilton + \basics + 2\hamilton' && \mbox{($X'$ is a cycle, so $\flow' =
\hamilton'$)} \\
 &\subseteq \hamilton + \basics && \mbox{(see proof of
\pref{2H'inH(G'odd)notinvol})} .
 \end{align*}
 Thus, because $f \in \even$ \see{WeirdEven<>Weight}, we may write $f = H + Q$,
with $H \in \hamilton$ and $Q \in \basics$. 
 By assumption, the weighted sum of the edge-flows of~$f$ is divisible by~$4$. By
Lemma~\ref{WeirdVanish}, the weighted sum of the edge-flows of~$H$ is also
divisible by~$4$. Therefore the weighted sum of the edge-flows of~$Q$ must also
be divisible by~$4$. So there is no harm in replacing~$f$ with~$Q$.

\begin{step} 
 Completion of the proof.
 \end{step}
 From Step~\ref{WeirdCasePf-f=4cyc}, we may assume that $f$ is a sum of some
number of basic $4$-cycles. The weighted sum of the edge-flows of any basic
$4$-cycle is $\pm2$, so we conclude that the number of $4$-cycles in the sum is
even. Thus, Proposition~\ref{4cyc=}, with the roles of $s$ and $t$ interchanged,
implies $f \in \hamilton$, as desired.
 \end{proof}

\begin{cor} \label{deg4+4cyc}
 If $s \in S$, such that
 \begin{itemize}
 \item $|S'| = 4$,
 \item $|G'|$ is even,
 and
 \item $\hamilton$ contains some basic $4$-cycle~$C$ of~$X'$,
 \end{itemize}
 then $\even' \subseteq \hamilton + \hamilton'$.
 \end{cor}

\begin{proof}
 We may assume that either
 \begin{itemize}
 \item $X'$ is the square of an even cycle, or
 \item $X'$ is not bipartite, and $|G|$ is not divisible by~$4$,
 \end{itemize}
 for otherwise Proposition~\ref{degree4} implies $\even' \subseteq \hamilton'$.

 Given any even flow $f \in \even'$, we wish to show that $f \in \hamilton +
\hamilton'$.   Under the weighting of~$X'$ specified in
Proposition~\ref{SquareCycle} or~\ref{WeirdCase} (as appropriate), the weighted
sum of the edge-flows of~$C$ is~$\pm2$, and, because $f$~is an even flow, it is
not difficult to see that the weighted sum of the edge-flows of~$f$ is even.
Therefore, there is an integer~$m$, such that the weighted sum of the edge-flows
of $f - mC$ is~$0$. Therefore, Proposition~\ref{SquareCycle} or~\ref{WeirdCase}
(as appropriate) asserts that $f - mC \in \hamilton'$. Because $mC \in
\hamilton$, we conclude that $f \in \hamilton + \hamilton'$, as desired.
 \end{proof}

\begin{cor} \label{weird+2tn}
 If $s,t \in S$, such that
 \begin{itemize}
 \item $|S'| = 4$,
 \item $|G'|$ is even, but not divisible by~$4$,
 \item $X'$ is not the square of an even cycle,
 \item $t \in S'$,
 \item $|t|$~is odd, and
 \item $\hamilton$ contains the flow $2 (t^{|t|})$,
 \end{itemize}
 then $\even' \subseteq \hamilton + \hamilton'$.
 \end{cor}

\begin{proof}
 Let $C = 2 (t^{|t|})$. Under the weighting of~$X'$ specified in
Proposition~\ref{WeirdCase}, the sum of the edge-flows of~$C$ is $2|t| \equiv 2
\pmod{4}$. Thus, for any $f \in \even'$, the weighted sum of the edge-flows of
either~$f$ or $f - C$ is is divisible by~$4$. Therefore,
Proposition~\ref{WeirdCase} asserts that either $f \in \hamilton'$ or $f - C \in
\hamilton'$. Because $C \in \hamilton$, we conclude that $f \in \hamilton +
\hamilton'$, as desired.
 \end{proof}

\section{The proof of Lemma~\ref{WeirdVanish}}
 \label{WeirdVanishSect}

This entire section is devoted to the proof of Lemma~\ref{WeirdVanish}. (None of
the definitions, notation, or intermediate results are utilized in other sections
of the paper.) After embedding $X$ on the torus \see{EmbedOnT2}, we assign an
integer modulo~$4$, called the ``imbalance" \fullsee{ImbDefn}{imb}  to
certain cycles (namely, those that are ``essential" and have even length). Then
we show that this geometrically-defined invariant can be used to calculate the
weighted sum of the edge-flows of the cycle \see{wt=len+imb}.
Lemma~\ref{WeirdVanish} follows easily from this formula.

\begin{assump}
 Throughout this section, $S$, $G$, $t$, $u$, and the weighting of~$X$ are as
described in Proposition~\ref{WeirdCase}.
 \end{assump}

\begin{notation}
 We use 
 \begin{itemize}
 \item $\wt(f)$ to denote the weighted sum of the edge-flows of a flow~$f$, and
 \item $\len(P)$ to denote the length of a path~$P$.
 \end{itemize}
 \end{notation}

\begin{defn} \label{EmbedOnT2} \ 
 \begin{enumerate}
 \item Define
 \begin{itemize}
 \item $m = |t|$,
 \item $n = |G : \langle t \rangle|$,
 and 
 \item choose an \emph{even} integer~$r$, such that $u^n = t^r$.
 \end{itemize}
 \item Embed $X$ on the torus $\torus^2 = \real^2/\integer^2$, by identifying the
vertex $t^a u^b$ of~$X$ with the point
 $$\left( \frac{a}{m} + \frac{rb}{mn}, \frac{b}{n} \right) $$
 of~$\torus^2$, and embedding the edges in the natural way (as line segments).
 \end{enumerate}
 \end{defn}

\begin{notation}
 Suppose $P$ is any path in~$X$, and $C$ is any cycle in~$X$, such that neither
the initial vertex nor the terminal vertex of~$P$ lies on~$C$. Intuitively, we
would like to
 $$ \begin{matrix}
 \mbox{define $\chi_C(P)$ to be the parity of} \\
 \mbox{the number of times that~$P$ crosses~$C$.}
 \end{matrix} $$
 (Note that if $P$ coincides with~$C$ on some subpath, then this is counted as a
crossing if and only if $P$~exits $C$~on the opposite side from the one it
entered on.)  It would be possible to formalize the definition in purely
combinatorial terms, but we find it convenient to use a topological approach.

 We may think of~$P$ as a continuous curve on the torus, and $C$ as a knot (or
loop) on the torus. By perturbing~$P$ slightly, we can obtain a curve~$P'$ on the
torus, with the same endpoints as~$P$, such that $P'$~is homotopic to~$P$, and
every intersection of~$P'$ with~$C$ is transverse (and is not a double point
of~$P'$). Let 
 $$ \mbox{$\chi_C(P) = |P' \cap C| \pmod{2}$.} $$
 This is well-defined (modulo~$2$) because $P'$ is homotopic to~$P$ (cf.\
\cite[\S\S73--74]{SeifertThrelfall}).
 \end{notation}

\begin{defn}
 A cycle $C$ in~$X$ is \emph{essential} if the corresponding knot on the torus is
\emph{not} homotopic to a point.

More concretely, a cycle $[v](s_1,\ldots,s_n)$ is essential if and only if either
 $$ |\{\, i \mid s_i = t \,\} | \neq |\{\, i \mid s_i = t^{-1} \,\} | $$
 or
 $$ |\{\, i \mid s_i = u \,\} | \neq |\{\, i \mid s_i = u^{-1} \,\} | .$$
 \end{defn}

\begin{defn} \label{ImbDefn} 
 Let $C$ be any essential, even cycle in~$X$.
 \begin{enumerate}
 \item \label{ImbDefn-color} 
 For two vertices $v$ and~$w$ in $X \setminus C$, we say that \emph{$v$
and~$w$ have the same color {\rm(}with respect to~$C${\rm)},} if 
 $$ \mbox{$\len(P) + \chi_C(P)$ is even,} $$
 where $P$~is any path in~$X$ from~$v$ to~$w$. (This is independent of the choice
of the path~$P$ (see Lemma~\fullref{esseven}{len+chi} below).) This is an
equivalence relation on $V(X \setminus C)$, and has (no more than) two
equivalence classes. 

We may refer to the vertices in one equivalence class as being ``black," and the
vertices in the other equivalence class as being ``white." This is a $2$-coloring
of $X\setminus C$.

 \item \label{ImbDefn-imb} 
 If  $K$ and~$W$ are the number of vertices of $X \setminus C$ that are
black, and the number that are white, respectively, we define the
\emph{imbalance} $\imb(C)$ to be
 $$ \imb(C) = K - W \pmod{4}.$$
 Because 
 $ K + W = |G| - \len(C) $
 is even, we know that $ K - W $ is either $0$ or~$2$ $\pmod{4}$; therefore, $ K
- W \equiv W-K \pmod{4} $, so $\imb(C)$ is well-defined (modulo~$ 4 $),
independent of the choice of which equivalence class is colored black and which
is colored white. 
 \end{enumerate}
 \end{defn}

\begin{obs} \label{altimb} 
 Because this concept is the foundation of this entire section, we describe an
alternate approach to the $2$-coloring that determines $\imb(C)$. The graph $X$
has a natural double cover~$\Xdouble$ that is bipartite. Specifically, 
 $$\Xdouble = \Cay \bigl( G_2; \{t_2, u_2\} \bigr),$$
 where 
 $$ G_2 = \langle \, t_2, u_2 \mid 
 \mbox{ $t_2^{2m}=e$, $u_2^n = t_2^r$, $t_2 u_2 = u_2 t_2$} \, \rangle .$$
 The inverse image of~$C$ in~$\Xdouble$ consists of two disjoint cycles $\Cdouble_1$
and~$\Cdouble_2$, with $\len(\Cdouble_1) = \len(\Cdouble_2) = \len(C)$. (This would
be false if $C$ were not an even cycle.)
 There is a natural embedding of~$\Xdouble$ on the torus~$\torus^2$, and $\torus^2
\setminus ( \Cdouble_1 \cup \Cdouble_2 ) $ has exactly two connected
components. (This would be false if $C$ were not essential.) Choose one connected
component~$\Xdouble^\circ$. The vertices in~$\Xdouble^\circ$ are in one-to-one
correspondence with the vertices in~$X \setminus C$. Because $\Xdouble$ is
bipartite, the vertices in~$\Xdouble^\circ$ have a natural $2$-coloring. Under
the natural correspondence with $V(X \setminus C)$, this is precisely the
$2$-coloring defined above, up to the arbitrary choice of which equivalence class
will be black and which will be white.

 To see that this is the same $2$-coloring, note that:
 \begin{itemize}
 \item Each vertex~$v$ of~$X$ has two inverse images $\vdouble_1$ and~$\vdouble_2$
in~$\Xdouble$, one in each component of $\torus^2 \setminus ( \Cdouble_1 \cup
\Cdouble_2 )$. 
 \item Any path from~$\vdouble_1$ to~$\vdouble_2$ in~$\Xdouble$ has odd length, so
$\vdouble_1$ and~$\vdouble_2$ are of opposite color under the $2$-coloring
of~$\Xdouble$. 
 \item A continuous curve in~$\torus^2$ crosses $\Cdouble_1 \cup \Cdouble_2$ an odd
number of times if and only if its two endpoints are in different components of
$\torus^2 \setminus (\Cdouble_1 \cup \Cdouble_2)$ (unless the curve has an
endpoint on $\Cdouble_1 \cup \Cdouble_2$).
 \end{itemize}
 Therefore, two vertices $v$ and~$w$ of $X \setminus C$ have the same color if
and only if $\len(\Pdouble) + \chi_{\Cdouble_1}(\Pdouble) +
\chi_{\Cdouble_2}(\Pdouble)$ is even, where $\Pdouble$~is any path in~$\Xdouble$ from
$\vdouble$ to~$\wdouble$, and $\vdouble$ and $\wdouble$ are inverse images of $v$
and~$w$, respectively (since, by the second and third bullets, we may assume that
$\vdouble$ and~$\wdouble$ are in the same component of $\Xdouble \setminus (\Cdouble_1
\cup \Cdouble_2)$). This establishes that this alternate approach is indeed
consistent with Definition~\ref{ImbDefn}. Furthermore, we see that replacing
$\Xdouble^\circ$ with the other component of $\torus^2 \setminus ( \Cdouble_1 \cup
\Cdouble_2 )$ reverses the color of each vertex of $X \setminus C$.
 \end{obs}

\begin{rem}
 As an aid to the reader's intuition, we informally describe two other ways of
thinking about the $2$-coloring. This also proves that the definition in
\fullref{ImbDefn}{color} is independent of the choice of~$P$, which we will
prove in another way in Lemma~\fullref{esseven}{len+chi}.
 \begin{enumerate}
 \item We give an orientation and a $2$-coloring to the cycle~$C$. We
then $2$-color the rest of~$X$ consistently with $C$ on, say, the left, and
inconsistently with $C$ on the right.
 \item We ``refine" or ``subdivide" the grid~$X$ by a factor of, say,~$3$ in each
direction. Then the complement of~$C$ in the new graph is connected and bipartite,
so it has a natural $2$-coloring, and this restricts to our coloring of $X
\setminus C$.
 \end{enumerate}
  \end{rem}

\begin{notation}\label{OmegaDefn}
 For any essential, even cycle~$C$ in~$X$, let
 $$ \wli (C) \equiv \wt(C) + \len(C) + \imb(C) \pmod{4}
 .$$
 \end{notation}

Our main task in this section is to show that if $C$ is any essential, even cycle
in~$X$, then $\wli(C) \equiv 2 \pmod{4}$ \cf{wt=len+imb}. (This is accomplished
by reducing to the case in which $C$ is ``monotonic"
\seeand{invariance}{makemonotonic}, and then calculating the imbalance in this
special case \see{Imb(monotonic)}.) Once this formula has been established, it
will be easy to prove Lemma~\ref{WeirdVanish}.

\begin{notation}
 Let 
 $$ \mbox{$\widetilde{t} = ( 1,0 )$ and $\widetilde{u} = ( 0,1 )$} ,$$
 and define
 $$ \Xtilde = \Cay \bigl( \integer \times \integer ; \{  \widetilde{t}^{\pm1},
\widetilde{u}^{\pm1} \} \bigr) .$$
 Note that:
 \begin{enumerate}
 \item There is a natural covering map~$\pi$ from~$\Xtilde$ to~$X$, defined by
$\pi( \widetilde{t}^i \widetilde{u}^j ) = t^i u^j$,
 and
 \item $\Xtilde$ has a natural embedding in the plane~$\real^2$.
 \end{enumerate}
 \end{notation}

\begin{rem}
 For convenience, we relax the usual definition of \emph{path} to allow the
initial vertex to be equal to the terminal vertex; that is, cycles are
considered to be paths.
 \end{rem}

\begin{defn}
 Suppose $P$ is a path in~$X$ and $\Ptilde$ is a path in~$\Xtilde$.
 We say that $\Ptilde$ is a \emph{lift} of~$P$ if $\pi(\Ptilde) = P$.

Note that if $P$ is a path in~$X$, and $\widetilde{v}$ is any vertex
of~$\Xtilde$, such that $\pi(\widetilde{v}$) is the initial vertex of~$P$, then
there is a unique lift~$\Ptilde$ of~$P$, such that the initial vertex of~$\Ptilde$
is~$\widetilde{v}$.  Namely, if  $P=[\pi(\widetilde{v})](s_1,\ldots,s_{\ell})$,
then $\Ptilde= [\widetilde{v}](\widetilde{s_1},\ldots,\widetilde{s_{\ell}})$.
Notice that $\len(\Ptilde) = \len(P)$.
 \end{defn}

\begin{defn}
 Suppose $C$ is a cycle in~$X$, and $\Ptilde$ is any lift of~$C$ to a path
in~$\Xtilde$.
 If $( x_1, y_1 )$ and $(x_2, y_2)$ are the initial vertex of~$\Ptilde$ and the
terminal vertex of~$\Ptilde$, respectively, then, because $C$ is a cycle, there
exist $p,q \in \integer$ with
 $$ ( x_2, y_2 ) - (x_1,y_1) = p\,( m,0 ) + q\, ( -r,n ) .$$
 The \emph{knot class} $\knot(C)$ of~$C$ is the ordered pair $(p,q)$.
 (Because we have not specified an orientation of~$C$, this is well-defined only
up to a sign; that is, we do not distinguish the knot class $(p,q)$ from the knot
class $(-p,-q)$. In our applications, it is only the parities of $p$ and~$q$ that
are relevant, and the parities are not affected by a change of sign. Note also
that $p$ depends on the choice of~$r$, but its parity does not.)

 In topological terms, $p$~is the number of times that the knot corresponding
to~$C$ wraps around the torus longitudinally, and $q$~is the number of times that
the knot wraps around the torus meridionally \cite[pp.~17--19]{Rolfsen}.
 \end{defn}

\begin{obs} \label{Ess<>Lift}
 Let $C$ be a cycle in~$X$. Then the following are equivalent:
 \begin{enumerate} \renewcommand{\theenumi}{\alph{enumi}}
 \item $C$ is essential.
 \item $\knot(C) \neq (0,0)$ \cite[Exer.~14 on p.~25]{Rolfsen}.
  \item If $\Ptilde$ is a lift of~$C$ to a path in~$\Xtilde$, then $\Ptilde$ is
\emph{not} a cycle in~$\Xtilde$; that is, the terminal vertex of~$\Ptilde$ is not
equal to the initial vertex of~$\Ptilde$.
 \end{enumerate}
 \end{obs}

\begin{lem} \label{esseven}
 Suppose $C$ is any essential, even cycle in~$X$.
 \begin{enumerate}
 \item \label{esseven-oddvertwrap}
 Let $\Ptilde$ be a lift of~$C$ to a path in~$\Xtilde$. If $( x_1, y_1 )$ and
$(x_2, y_2)$ are the initial vertex of~$\Ptilde$ and the terminal vertex
of~$\Ptilde$, respectively, then $y_2 - y_1 \equiv 2 \pmod{4}$.
 \item \label{esseven-len+chi}
 Let
 \begin{itemize}
 \item $v$ and $w$ be vertices of~$X$ that do not lie on~$C$,
 and  
 \item $P$ and~$Q$ be paths in~$X$ from~$v$ to~$w$.
 \end{itemize}
 Then $\len(P) + \chi_C(P) \equiv \len(Q) + \chi_C(Q) \pmod{2}$.
 \end{enumerate}
 \end{lem}

\begin{proof}
 \pref{esseven-oddvertwrap} We may assume $( x_1, y_1 ) = ( 0,0 )$. Let $\knot(C)
= (p,q)$. Then
 $$ ( x_2, y_2 ) = p\,( m,0 ) + q\, ( -r,n )
 = ( pm - rq , qn ) .$$
 Note that, because $m$~is odd, but $r$ and~$n$ are even, we have
 \begin{equation} \label{len(C)=p}
 \len(C) \equiv x_2 + y_2 = pm - rq + qn \equiv p \pmod{2} .
 \end{equation}
 Since $C$ is an even cycle, we conclude that $p$~is even. 

Because $(p,q)$ is the knot class of an essential knot (namely,~$C$), a theorem of
topology asserts that $\gcd(p,q) = 1$ \cite[p.~19]{Rolfsen}. (This can also be
proved combinatorially.) Because $p$~is even, this implies that $q$~is odd.
Combining this with the fact that $n \equiv 2 \pmod{4}$, we conclude that
 $$ y_2 - y_1 = y_2 = qn \equiv 2 \pmod{4} ,$$
 as desired.

 \pref{esseven-len+chi} It suffices to show that if $C'$ is any cycle in~$X$,
then $\len(C') + \chi_C(C')$ is even.  Let $\knot(C') = (p',q')$. A theorem of
topology \cite[Exer.~7 on p.~28]{Rolfsen} asserts that $\chi_C(C') \equiv p q' -
p' q \pmod{2}$. From the proof of~\pref{esseven-oddvertwrap}, we know that $p$~is
even and $q$~is odd. Thus, 
 $\chi_C(C') \equiv p' \pmod{2}$. Furthermore, \eqref{len(C)=p} (which did not
use the fact that~$C$ is even), with $C'$ in the role of~$C$, shows that $\len(C')
\equiv p' \pmod{2}$. Thus, 
 $\len(C') + \chi_C(C') \equiv p' + p' \pmod{2}$ is even.
 \end{proof}

\begin{defn}
 If $C$ is any cycle in~$X$ that is \emph{not} essential, then $C$ may be
lifted to a cycle~$\Ctilde$ in~$\Xtilde$ \see{Ess<>Lift}. The cycle $\Ctilde$ has
a well-defined interior and exterior in~$\real^2$. Let $\Itilde$ be the set of
vertices of~$\Xtilde$ in the interior of~$\Ctilde$. Then the restriction of~$\pi$
to~$\Itilde$ is one-to-one and maps onto a subset of $X \setminus C$ that
is independent of the choice of the lift~$\Ctilde$. We say that the vertices
in~$\pi(\Itilde)$ are \emph{the vertices in the interior of the region bounded
by~$C$}. 
 \end{defn}

\begin{lem}[Pick's Theorem] \label{PickThm}
 If $C$ is any cycle in~$X$ that is \emph{not} essential, and $N$~is the number
of vertices in the region enclosed by~$C$, then
 $$ \wt(C) \equiv \len(C) + 2N - 2 \pmod{4} .$$
 \end{lem}

\begin{proof}
  Lift $C$ to a cycle $\Ctilde$ in~$\Xtilde$, and let $A$ be the area of the
region bounded by~$\Ctilde$. Then $C$ is the sum of~$A$ basic $4$-cycles (and the
weight of any basic $4$-cycle is $\pm2$), so 
 $$ \wt(C) \equiv 2A \pmod{4} .$$
 From Pick's Theorem \cite[pp.~27--31]{Honsberger}, we know that 
 $$ 2A = \len(\Ctilde) + 2N - 2.$$
 Because $\len(\Ctilde) = \len(C)$, the desired conclusion is immediate.
 \end{proof}

It will be helpful to know that making certain changes to an essential, even
cycle~$C$ does not affect $\wli(C) \pmod{4}$.

\begin{eg} \label{wli4cyc}
 Let $C$ be an essential, even cycle in~$X$.
 \begin{enumerate}
 \item Suppose $C$ contains the subpath $[v](x,y,x^{-1})$, for some $v \in G$ and
some $x,y \in S$. Let $C'$ be the cycle obtained from~$C$ by replacing this
subpath with the single edge $[v](y)$. Then
 \begin{itemize}
 \item $\wt(C) - \wt(C') = \wt \bigl( [v](x,y,x^{-1},y^{-1} ) \bigr) = \pm 2$,
 \item $\len(C) - \len(C') = 2$,
 and
 \item $\imb(C) - \imb(C') \equiv 0 \pmod{4}$, because $X \setminus C' = (X
\setminus C) \cup \{ vx, vxy\}$ and the two additional vertices are of opposite
color (because they are adjacent in~$X$).
 \end{itemize}
 So $\wli(C) \equiv \wli(C') \pmod{4}$.
 \item Suppose $C$ contains the subpath $[v](x,y)$, for some $v \in G$ and
some $x,y \in S$. If the vertex~$[v](y)$ is not on~$C$, let $C'$ be the cycle
obtained from~$C$ by replacing the subpath $[v](x,y)$ with the subpath
$[v](y,x)$. Then
 \begin{itemize}
 \item $\wt(C) - \wt(C') = \wt \bigl( [v](x,y,x^{-1},y^{-1} ) \bigr) = \pm 2$,
 \item $\len(C) - \len(C') = 0$,
 and
 \item $\imb(C) - \imb(C') \equiv 2 \pmod{4}$, because the symmetric
difference of $X \setminus C$ and $X \setminus C'$ is $\{ vx, vy\}$, and, using
consistent colorings, the color of~$vy$ in $X \setminus C$ is the opposite of the
color of~$vx$ in $X \setminus C'$ (because they are an even distance apart, but on
opposite sides of the two cycles).
 \end{itemize}
 So $\wli(C) \equiv \wli(C')  \pmod{4}$.
 \end{enumerate}
 \end{eg}

The following result is a weak form of the assertion that changing $C$ by a
homotopy does not change $\wli(C)$. On an intuitive level, it can be justified by
claiming that a sequence of the two types of replacements described in
Example~\ref{wli4cyc} will transform~$C$ into $C' = C - P + Q$. Our formal proof
takes a different approach.

\begin{prop} \label{invariance}
 Suppose
 \begin{itemize}
 \item $C$ is an essential, even cycle in~$X$,
 \item $P$ is a subpath of~$C$,
 \item $Q$ is a path in~$X$,
 \item the initial vertex and terminal vertex of~$Q$ are the same as those of~$P$,
 \item $Q$ does not intersect~$C$ {\rm(}except at the  endpoints of~$Q${\rm)},
 and
 \item the cycle $P - Q$ is not essential.
 \end{itemize}
 Then $C - P + Q$ is an essential, even cycle, and
 $\wli( C - P + Q ) \equiv \wli(C)\pmod{4}$.
 \end{prop}

\begin{proof}
 Because $P - Q$ is not essential, and $\Xtilde$ is bipartite, it is clear that
$P - Q$ is an even cycle \cf{Ess<>Lift}. So $C - P + Q$ is an even cycle. It is
also easy to see that $C - P + Q$, like~$C$, must be essential. For example, one
may note that
 \begin{align*}
  \knot(C - P + Q) &= \knot(C) - \knot(P - Q) \\
 &= \knot(C) - (0,0) \\
 &\neq (0,0) .
 \end{align*}

Letting $N$ be the number of vertices in the region enclosed by $P - Q$, and
applying Pick's Theorem \pref{PickThm}, we have
 \begin{align}
 \wt(C) - \wt( C - P + Q )
 &= \wt( P - Q ) \notag\\ 
 &\equiv \len( P - Q ) + 2N - 2 \pmod{4} \label{chgwt} \\
 &= \len( P)  + \len(Q ) + 2N - 2
 .\notag
 \end{align}
 Obviously,
 \begin{equation}\label{chglen}
 \len( C) - \len( C - P + Q )
 = \len( P) - \len(Q) .
 \end{equation}
 Also, we know, from Lemma~\ref{chgimblem} below, that
 \begin{equation}\label{chgimb} 
 \imb( C) - \imb( C - P + Q )
 \equiv 2N + 2 \len(P) - 2 \pmod{4}.
 \end{equation}

Combining \pref{chgwt}, \pref{chglen}, and \pref{chgimb}, we obtain
 \begin{align*}
 \wli(C) - \wli(C - P + Q )
 &= \bigl( \wt (C) - \wt (C - P + Q ) \bigr) \\
 & \qquad +  \bigl( \len (C) - \len (C - P + Q ) \bigr) \\
 & \qquad\qquad
 + \bigl( \imb (C) - \imb (C - P + Q ) \bigr) \\
 &\equiv \bigl( \len(P) + \len(Q) + 2N - 2 \bigr) \\
 & \qquad
 +  \bigl( \len (P) - \len (Q) \bigr) \\
 & \qquad\qquad
 + \bigl( 2N + 2\len(P) - 2 \bigr) \\
 &= 4 \len(P) + 4 N - 4  \\
 &\equiv 0 \pmod{4}
 ,
 \end{align*}
 so $ \wli(C) \equiv \wli(C - P + Q ) \pmod{4}$, as desired.
 \end{proof}

\begin{lem} \label{chgimblem}
 Under the assumptions of Proposition~\ref{invariance}, and letting $N$ be the
number of vertices in the region enclosed by $P-Q$, we have
 $\imb( C - P + Q ) - \imb(C) \equiv 2N + 2 \len(P) - 2 \pmod{4}$.
 \end{lem}

\begin{proof}
 We use the description of~$\imb(C)$ given in Observation~\ref{altimb}.
 \begin{itemize}
 \item Let $\Cdouble_1$ and $\Cdouble_2$ be the two lifts of~$C$ to cycles in~$\Xdouble$.
 \item Let $\Pdouble_1$ and $\Pdouble_2$ be the two lifts of~$P$ to paths in~$\Xdouble$,
with $\Pdouble_1 \subseteq \Cdouble_1$ and $\Pdouble_2 \subseteq \Cdouble_2$.
 \item Let $\Qdouble_1$ and $\Qdouble_2$ be the two lifts of~$Q$ to paths in~$\Xdouble$,
such that $\Qdouble_1$ has the same endpoints as~$\Pdouble_1$, and $\Qdouble_2$~has
the same endpoints as~$\Pdouble_2$.
 \item For a curve~$R$ on~$\torus^2$, we use $R^\circ$
to denote the \emph{interior} of~$R$; that is $R^\circ = R \setminus
\{v_1,v_2\}$, where $v_1$ and~$v_2$ are the endpoints of~$R$.
 \item Let $\Xdouble^\circ$ be the component of $\torus^2 \setminus (
\Cdouble_1 \cup \Cdouble_2 )$ that contains~$\Qdouble_1^\circ$.
 \item Let $\ldouble_1$ and~$\ldouble_2$ be the interior regions in~$\torus^2$ that
are bounded by the cycles $\Pdouble_1 - \Qdouble_1$ and $\Pdouble_2 - \Qdouble_2$,
respectively.
 \item For any subset~$A$ of~$\torus^2$, we let 
 $K_A$ and~$W_A$, respectively, be the number of black vertices of~$\Xdouble$ that are
contained in~$A$, and the number of white vertices of~$\Xdouble$ that are contained
in~$A$. We define 
 $$\kw(A) = K_A - W_A .$$
 \end{itemize}
 Note that:
 \begin{itemize}
 \item From the definition of $\imb(C)$, we have
 $$\imb(C) \equiv \kw(\Xdouble^\circ) \pmod{4}.$$
 \item The inverse image of~$C-P+Q$ in~$\Xdouble$ is 
 $$ (\Cdouble_1 - \Pdouble_1 + \Qdouble_1) \cup (\Cdouble_2 - \Pdouble_2 + \Qdouble_2)
,$$ 
 and one of the components of the complement of this inverse image is
 $$ \bigl( \Xdouble^\circ \setminus (\ldouble_1 \cup \Qdouble_1^\circ) \bigr) \cup
(\ldouble_2 \cup \Pdouble_2^\circ) .$$
 Therefore,
 \begin{align*}
  \imb(C-P+Q) \equiv \kw(\Xdouble^\circ) - \kw(\ldouble_1) -
\kw(\Qdouble_1^\circ)
 \qquad& \\
 {} + \kw(\ldouble_2) + \kw(\Pdouble_2^\circ)  & \pmod{4}.
 \end{align*}
 \item For any vertex $\vdouble_1$ of~$\ldouble_1$, the image~$v$ of~$\vdouble_1$
in~$X$ is a vertex in the interior of the region enclosed by $P - Q$. Thus, there
is a vertex~$\vdouble_2$ in~$\ldouble_2$ whose image is also~$v$. It was noted in
Observation~\ref{altimb} that the color of~$\vdouble_1$ must be the opposite of
the color of~$\vdouble_2$. From this (and the same argument with $\ldouble_1$
and~$\ldouble_2$ interchanged), we conclude that $K_{\ldouble_1} = W_{\ldouble_2}$
and $W_{\ldouble_1} = K_{\ldouble_2}$. Therefore
 \begin{align*}
  -\kw(\ldouble_1) + \kw(\ldouble_2)
 &= - (K_{\ldouble_1} - W_{\ldouble_1}) + (K_{\ldouble_2} - W_{\ldouble_2}) \\
 &= - (W_{\ldouble_2} - K_{\ldouble_2}) + (K_{\ldouble_2} - W_{\ldouble_2}) \\
 &= 2 (K_{\ldouble_2} - W_{\ldouble_2}) \\
 &\equiv 2 (K_{\ldouble_2} + W_{\ldouble_2}) \pmod{4} \\
 &= 2N
 . 
 \end{align*}
 \item If $P$ and~$Q$ are of odd length, then $\Pdouble_2$ and $\Qdouble_1$ both
have an even number of vertices, half black and half white (and the same
is true of their interiors), so $\kw(\Qdouble_1^\circ) = 0$ and
$\kw(\Pdouble_2^\circ) = 0$. Because
 $2 \len(P) - 2$ is divisible by~$4$,
 we have 
 \begin{align} \label{Imb(Q)vsImb(P)odd(alt)}
 -\kw(\Qdouble_1^\circ) + \kw(\Pdouble_2^\circ)
  &= 0 \\
 &\equiv 2 \len(P) - 2 \pmod{4} . \notag
 \end{align}
 On the other hand, if $P$ and~$Q$ are of even length, then $\kw(\Qdouble_1^\circ)
= \pm1$ and $\kw(\Pdouble_2^\circ) = \pm1$. Assume, without loss of generality,
that $\kw(\Qdouble_1^\circ) = -1$. Then the interior of $\Qdouble_1$ has an extra
white vertex, so the endpoints of~$\Qdouble_1$ must both be black. These are the
same as the endpoints of~$\Pdouble_1$, so the endpoints of~$\Pdouble_1$ must both
be black. The endpoints of~$\Pdouble_2$ are of the opposite color (because they
project to the same vertices of~$X$); they must both be white. Hence, the interior
of~$\Pdouble_2$ has an extra black vertex, so $\kw(\Pdouble_2^\circ) = 1$. Noting
that $2 \len(P) - 2 \equiv 2 \pmod{4}$ in the current case, we conclude that
  \begin{align*}
 -\kw(\Qdouble_1^\circ) + \kw(\Pdouble_2^\circ)
  &= -(-1) + 1 \\
 &\equiv 2 \len(P) - 2 \pmod{4} ,
 \end{align*}
 which has the same form as the conclusion of \pref{Imb(Q)vsImb(P)odd(alt)}.
 \end{itemize}
 The desired conclusion is obtained by combining these calculations.
 \end{proof}

\begin{defn} \ 
 \begin{itemize}
 \item A path~$[v](s_1,\ldots,s_n)$ in~$X$ is \emph{increasing} if 
 $\{ s_1,\ldots,s_n \} \subseteq \{t,u\}$;
 that is, if no~$s_i$ is equal to $t^{-1}$ or~$u^{-1}$.
 \item A path~$[v](s_1,\ldots,s_n)$ in~$X$ is \emph{decreasing} if 
 $\{ s_1,\ldots,s_n \} \subseteq \{t,u^{-1}\}$;
 that is, if no~$s_i$ is equal to $t^{-1}$ or~$u$.
 \item A path in~$X$ is \emph{monotonic} if it is either increasing or decreasing.
 \end{itemize}
 \end{defn}

\begin{cor} \label{makemonotonic}
 If $C$ is any essential, even cycle in~$X$, then there is an essential, even
cycle~$C'$ in~$X$, such that
 \begin{enumerate}
 \item $\wli(C') \equiv \wli(C) \pmod{4}$,
 and
 \item $C'$ is monotonic.
 \end{enumerate}
 \end{cor}

\begin{proof}
 We may assume, without loss of generality, that 
 \begin{equation} \label{noshorter}
 \begin{matrix}
 \mbox{there does not exist an essential, even cycle~$C''$ in~$X$, such that} \\
 \mbox{$ \wli(C'') \equiv \wli(C) \pmod{4}$
 and 
 $\len(C'') < \len ( C ) $.}
 \end{matrix}
 \end{equation}

Lift $C$ to a path $\Ctilde$ in~$\Xtilde$.  Assume, for simplicity, that the
initial vertex of~$\Ctilde$ is $ ( 0,0 ) $, and let $ ( a,b ) $ be the terminal
vertex of~$\Ctilde$. We may assume $a \ge 0$, by reversing the orientation
of~$C$, if necessary.
 
Let us show that we may also assume $b \ge 0$, by interchanging $u$ with $u^{-1}$
(and negating~$r$), if necessary. The interchange negates~$b$ and preserves
monotonic paths. It does not affect $\len(C)$, $\imb(C)$, or $\wt(C)$, so it does
not affect $\wli(C)$. 

Suppose $C$ is not increasing. (This will lead to a contradiction.) Then either
$t^{-1}$ or~$u^{-1}$ appears in~$C$. The argument in the two cases is similar, so
let us assume, for definiteness, that an edge of the form $[v](t^{-1})$ appears
in~$C$. Because $b \ge 0$, we know that $C$ also contains an edge of the
form~$[v](t)$. Let $P'$ be a shortest subpath of~$C$ that contains at least one
edge of each form. Then $P'$ must 
 begin with an edge of one of these two forms, 
 end with an edge of the other form,
 and all other edges must be of the form $[v](u^{\pm1})$. 
 Now $P'$ cannot contain both an edge of the form  $[v](u)$ and an edge of the
form  $[v](u^{-1})$, so we may assume
 $$ P' = (t^{-1}, u^k, t)$$
 for some $k > 0$ (since other cases are similar).
 Let 
 \begin{itemize}
 \item $d_1,d_2 \ge 0$ be maximal, such that
 $C$ contains the paths
 $ [u^k] ( u^{-d_1} ) $
 and
 $ [u^{d_2} ]( u^{-d_2} ) $,
 \item $Q$ be the path
 $ [ u^{d_2} ] ( u^{k - d_1 - d_2} )$,
 and
 \item $P$ be the subpath of~$C$ from~$ u^{d_2} $ to~$u^{k-d_1}$
(containing~$P'$).
 \end{itemize}
 Note that:
 \begin{enumerate}
 \item The path $Q$ is disjoint from~$C$, except at its endpoints. (This follows
from the minimality of $k = \len(P')$ and the maximality of $d_1$ and~$d_2$.)
 \item The cycle $P - Q$ is not essential (because it lifts to a cycle in
$\Xtilde$, namely, the subpath of~$\Ctilde$ from~$(0,0 )$ to $( 0,k )$, plus the
vertical path from $( 0,k )$ to $(0,0 )$).
 \item The cycle $C - P + Q$ is strictly shorter than~$C$. (The length of~$Q$ is
strictly less than the length of~$P$, because $Q$ is a straight (vertical) path,
but $P$, which includes the horizontal edge $( t^{-1} )$, is inefficient.)
 \end{enumerate}
 From Proposition~\ref{invariance}, we see that
 $ \wli(C - P + Q) \equiv \wli(C) \pmod{4}$. Because $C - P + Q$ is strictly
shorter than~$C$, this contradicts~\pref{noshorter}.
 \end{proof}

In order to calculate the imbalance of a monotonic cycle~$C$, we formulate a
more convenient description of the $2$-coloring of $X \setminus C$ in this
case.

\begin{notation}
  Let $C$ be any monotonic, essential, even cycle in~$X$.
 \begin{enumerate}
 \item For vertices $u$ and~$v$ of~$C$, we let $\dist_C(u,v)$ be the length of
a shortest subpath of~$C$ with endpoints $u$ and~$v$.
 \item For any vertex $v \in X \setminus C$, it is easy to see that there are a
unique vertex $c(v)$ of~$C$ and a unique $h(v) \in \integer^+$, such that
 $$ c(v) \, t^{h(v)} = v  $$
 and 
 $$ \mbox{$ \{ c(v) t, c(v) t^2, \ldots, c(v) t^{h(v)-1} \}$ is disjoint
from~$C$.} $$
 \item Similarly, for any vertex $v \in X \setminus C$, there are a unique
vertex $c'(v)$ of~$C$ and a unique $h'(v) \in \integer^+$, such that
 $$ c'(v) \, t^{-h'(v)} = v $$
 and
 $$ \mbox{$ \{ v t, v t^2, \ldots, v t^{h'(v)-1} \}$ is disjoint from~$C$.} $$
 \end{enumerate}
 \end{notation}

\begin{lem} \label{ImbIncDefn}
 Let $C$ be any monotonic, essential, even cycle in~$X$, and fix a vertex $c_0$
of~$C$. Then the black vertices can be distinguished from the white vertices by
the parity of 
 $\dist_C \bigl( c_0, c(v) \bigr) + h(v)$.

In other words, two vertices $v$ and~$w$ of $X \setminus C$ have the same color
if and only if 
 $$\dist_C \bigl( c_0, c(v) \bigr) + h(v) \equiv \dist_C \bigl( c_0, c(w) \bigr)
+ h(w) \pmod{2} . $$
 \end{lem}

\begin{proof}
 Let $P_0$ be a path of length $\dist_C \bigl( c(v), c(w) \bigr)$ from~$c(v)$
to~$c(w)$ in~$C$, and let 
 $P = [v](t^{-h(v)}, P_0, t^{h(w)})$.
 Then $P$ is a path from~$v$ to~$w$, such that 
 $$ \len(P) = \dist_C \bigl( c(v), c(w) \bigr) + h(v) + h(w) ,$$
 and it is easy to see that $\chi_C(P) = 0$. (For example, if $C$ is increasing,
then translating $P_0$ by the vector 
  $$\left( \frac{1}{2m} - \frac{r}{2mn}, -\frac{1}{2n} \right) $$
 in~$\torus^2$ results in a path that is disjoint from~$C$, so it is easy to find
a path homotopic to~$P$ that is disjoint from~$C$.) Therefore, $v$ and~$w$ have
the same color if and only if $\len(P)$ is even. The desired conclusion follows.
 \end{proof}

\begin{rmk}
 In Lemma~\ref{ImbIncDefn}, one could use $c'(v)$ and $h'(v)$, in place of $c(v)$
and $h(v)$: vertices $v$ and~$w$ of $X\setminus C$ have the same color if
and only if
 $$\dist_C \bigl( c_0, c'(v) \bigr) + h'(v) \equiv \dist_C \bigl( c_0, c'(w)
\bigr) + h'(w) \pmod{2} . $$

Note that, because $X$~is \emph{not} bipartite, there is some $v \in X
\setminus C$, such that 
 \begin{equation} \label{OppParity}
 \dist_C \bigl( c_0, c(v) \bigr) + h(v)
 \not\equiv \dist_C \bigl( c_0, c'(v) \bigr) + h'(v)
 \pmod{2} .
 \end{equation}
 Then it is easy to see that we have inequality for \emph{all} $v \in X
\setminus C$ (and all $c_0 \in C$).
 \end{rmk}

\begin{defn}
 We will call a $t$-edge $[t^i u^j](t)$ of~$X$
 $$ \begin{cases}
 \mbox{blue} & \mbox{if $j$ is even}; \\ \hfil \mbox{red} & \mbox{if $j$ is odd}.
 \end{cases} $$
 \end{defn}

\begin{prop} \label{Imb(monotonic)}
 If $C$ is any monotonic, essential, even cycle in~$X$, then $\imb(C) \equiv 2B
\pmod{4}$, where $B$~is the number of blue edges in~$C$.
 \end{prop}

\begin{proof}
 Assume, for definiteness, that $C$ is increasing. Then, because $C$ is not of the
form $[c_0](t^m)$, there are a sequence $(\ell_1, \ldots, \ell_z)$ of natural
numbers and some $c_0 \in C$, such that
 \begin{equation}\label{Cinc}
 C = [c_0] ( u, t^{\ell_1}, u, t^{\ell_2}, \ldots, u, t^{\ell_z} ) .
 \end{equation}
 Because $z$ is a multiple of~$n$, we know that
 $$ \mbox{$z$ is even}.$$
 Thus, as $C$ is even, there is no harm in assuming that $c_0\in \langle
t\rangle$, so
 \begin{equation}\label{B=Sum(l)}
 B=\sum_{i=1}^{z/2}\ell_{2i}.
 \end{equation}
 For $1 \le j \le z$, let 
 \begin{itemize}
 \item $c_j=c_0 (u t^{\ell_1})(u t^{\ell_2}) \cdots (u t^{\ell_j}) $,
 \item $h_j$ be the unique positive integer, such that 
 $$ \mbox{$c_j t^{h_j} \in C$ and $\{ c_j t, c_j t^2, \ldots, c_j t^{h_j-1} \}$ is
disjoint from~$C$,} $$
 \item $c'_j = c_j t^{h_j} $,
 \item $d_j = \dist_C \bigl( c_0, c_j \bigr)$, 
 \item $d'_j = \dist_C \bigl( c_0, c'_j \bigr) $,
 \item $K$ be the number of black vertices in~$X \setminus C$,
 and
 \item $W$ be the number of white vertices in~$X \setminus C$.
 \end{itemize}
 Note that $X \setminus C$ is the disjoint union of the paths $ [c_j
t](t^{h_j-2})$, with the convention that such a path is empty if $h_j = 1$. Thus,
we may calculate $K - W$ by determining the excess (which is negative if there is
actually a deficiency) of black vertices in each of these paths, and adding up
the results.

A path in $X\setminus C$
 \begin{itemize}
 \item has an excess of $1$~black vertex if and only if both of its endpoints are
black,
 \item has an excess of $-1$~black vertex if and only if both of its endpoints
are white,
 and
 \item has an excess of $0$  black vertices otherwise.
 \end{itemize}
 By interchanging black with white if necessary, we may assume that a vertex~$v$
of $X \setminus C$ is black if and only if $\dist_C \bigl( c_0, c(v) \bigr)
+ h(v)$ is even \see{ImbIncDefn}. 
 \begin{itemize}
 \item Thus, the initial vertex of the path $[c_j t](t^{h_j-2})$ is black if and
only if $d_j$ is odd.
 \item Using \pref{OppParity}, we see that the terminal vertex of the path $[c_j
t](t^{h_j-2})$ is black if and only if $d'_j$ is even.
 \end{itemize}
 Therefore
 \begin{align*}
 K - W
 &= | \{\, j \mid \mbox{$d_j$ is odd and $d'_j$ is even} \,\} |
 \\
 & {} \qquad - | \{\, j \mid \mbox{$d_j$ is even and $d'_j$ is odd} \,\} |
 \\
 &= \bigl( | \{\, j \mid \mbox{$d_j$ is odd and $d'_j$ is even} \,\} | \\
 & {} \qquad\qquad
 + | \{\, j \mid \mbox{$d_j$ is odd and $d'_j$ is odd} \,\} | \bigr)
 \\
 & {} \qquad - \bigl( | \{\, j \mid \mbox{$d_j$ is even and $d'_j$ is odd} \,\} |
 \\
 & {} \qquad\qquad\qquad
 + | \{\, j \mid \mbox{$d_j$ is odd and $d'_j$ is odd} \,\} | \bigr)
 \\
 &= | \{\, j \mid \mbox{$d_j$ is odd} \,\} |
 - | \{\, j \mid \mbox{$d'_j$ is odd} \,\} |
 .
 \end{align*}
 Note that a vertex~$w$ of $C$ belongs to
 $$
 \begin{cases}
 \{ c'_1, \ldots, c'_z\}
 & \mbox{if and only if the edge $[w u^{-1}] (u)$ is in~$C$}, \\
 \{ c_1, \ldots, c_z\}
 & \mbox{if and only if the edge $[w ] (u)$ is in~$C$}.
 \end{cases} $$
 This implies that $\{ c'_1 u^{-1}, \ldots, c'_z u^{-1} \} = \{ c_1, \ldots,
c_z\}$, so
 $$ | \{\, j \mid \mbox{$d'_j$ is odd} \,\} | 
 = | \{\, j \mid \mbox{$d_j$ is even} \,\} |
 = z - | \{\, j \mid \mbox{$d_j$ is odd} \,\} | .$$
 Therefore
 \begin{align}
 K - W &= 
 | \{\, j \mid \mbox{$d_j$ is odd} \,\} |
 - \bigl( z - | \{\, j \mid \mbox{$d_j$ is odd} \,\} | \bigr) 
 \label{KW=2s} \\
 &= 2 | \{\, j \mid \mbox{$d_j$ is odd} \,\} | - z \notag
 . \end{align}

 From \pref{Cinc} (and the definitions of $d_j$ and~$c_j$), we see that
 \begin{equation} \label{d=sum(ell)}
 d_j \equiv \sum_{i=1}^j (\ell_i + 1) \pmod{2}.
 \end{equation}
 Combining this with the fact that 
 \begin{equation} \label{wheniodd}
 \mbox{$z-i+1$ is odd if and only if $i$~is even,}
 \end{equation}
 we see that
 \begin{align}
 | \{ j \mid \mbox{$d_j$ is odd} \} |
 &\equiv \sum_{j=1}^z d_j \pmod{2} \notag \\
 &\equiv \sum_{i=1}^z (z - i + 1) (\ell_i + 1 ) \pmod{2}
 &&\mbox{\see{d=sum(ell)}} \label{numsodd} \\
 &\equiv \sum_{i=1}^{z/2} (\ell_{2i} + 1 ) \pmod{2} 
 &&\mbox{\see{wheniodd}} . \notag 
 \end{align}
 Hence
 \begin{align*}
 \imb(C)&\equiv
 K - W \pmod{4} &&\mbox{\see{ImbDefn}}\\
 &= 2 | \{ j \mid \mbox{$d_j$ is odd} \} | - z &&\mbox{\see{KW=2s}}\\
 &\equiv 2 \left( \sum_{i=1}^{z/2} (\ell_{2i} + 1 ) \right) - z 
\pmod{4}&&\mbox{\see{numsodd}}
 \\
 & = 2  \sum_{i=1}^{z/2} \ell_{2i} \\
 &= 2B&&\mbox{\see{B=Sum(l)}},
 \end{align*}
 as desired.
 \end{proof}

\begin{cor} \label{wt=len+imb}
 If $C$ is any essential, even cycle in~$X$, then
  $$ \wt(C) \equiv \len(C) + \imb(C) + 2 \pmod{4} .$$
 \end{cor}

\begin{proof}
 It suffices to show $\wli(C) \equiv 2 \pmod{4}$ (see Notation~\ref{OmegaDefn}).
From Corollary~\ref{makemonotonic}, we may assume $C$ is monotonic.

Let
 \begin{itemize}
 \item $B$ be the number of blue $t$-edges in~$C$, 
 \item $R$ be the number of red $t$-edges in~$C$, 
 and
 \item $U$ be the number of $u$-edges in~$C$.
 \end{itemize}
 We have:
 \begin{enumerate}
 \item $\wt(C) = B - R$ (because $C$ is monotonic),
 \item $\len(C) = B + R + U$,
 \item $\imb(C) \equiv 2B \pmod{4}$ \see{Imb(monotonic)}, and
 \item $U \equiv 2 \pmod{4}$ (because $C$ is monotonic;
see~\fullref{esseven}{oddvertwrap}).
 \end{enumerate}
 Therefore, modulo~$ 4 $, we have
 \begin{align*}
 \wli(C)
 &\equiv \wt(C) + \len(C) + \imb(C) \\
 &\equiv ( B - R ) + ( B + R + U ) + 2B \\
 &= 4B + U \\
 &\equiv 0 + 2 \\
 &= 2
 ,
 \end{align*}
 as desired.
 \end{proof}

\begin{proof}[{\bf Proof of Lemma~\ref{WeirdVanish}}]
 Note that  $\len(H) = |G| \equiv 2 \pmod{4}$.

\setcounter{proofcase}{0}

\begin{proofcase}
 Assume $H$ is essential.
 \end{proofcase}
 Because the complement $X \setminus H$ is empty, it is obvious, from the
definition, that 
 $ \imb(H) \equiv 0 \pmod{4}$. Applying Corollary~\ref{wt=len+imb}, we obtain
 $$ \wt(H) 
 \equiv \len(H) + \imb(H) + 2
 \equiv 2 + 0 + 2
 \equiv 0 \pmod{4} .$$

\begin{proofcase}
 Assume $H$ is not essential.
 \end{proofcase}
 From Pick's Theorem \pref{PickThm}, we know that 
 $$ \wt(H) \equiv \len(H) + 2N - 2 \pmod{4},$$
 where $N$~is the number of lattice points in the interior of the region bounded
by~$H$. Because the complement $X \setminus H$ is empty, we have $N = 0$.
Therefore
 $$ \wt(H) 
 \equiv 2 + 2(0) - 2
 = 0 \pmod{4} ,$$
 as desired.
 \end{proof}

\section{$4$-cycles in graphs containing \copy\GRID} \label{GridSect}

Proposition~\ref{P3xP3xP2} is a major ingredient in our study of graphs of degree
at least~$5$. The gist is that if $X$ contains 
 \begin{itemize}
 \item a subgraph~$Y$ that is isomorphic to $P_3 \boxprod P_3 \boxprod P_2$, and
 \item an appropriate hamiltonian cycle~$C$ in the complement $X \setminus
V(Y)$,
 \end{itemize}
 then many hamiltonian cycles in~$Y$ can be extended to hamiltonian cycles
in~$X$. This yields enough hamiltonian cycles to show that $\hamilton$ contains
every basic $4$-cycle.

Later sections do not appeal directly to this result, but only to
Corollaries~\ref{GridsHave4Cycles}, \ref{4DGridsHave4Cycles}, \ref{K2xK3xK3},
and~\ref{OddHeight4Cycles}.

\begin{prop} \label{P3xP3xP2}
 Suppose
 \begin{itemize}
 \item $S = \{s^{\pm1}, t^{\pm1}, u^{\pm1} \}$,
 \item $Y = \bigset{
  s^i t^j u^k }{
 \begin{matrix}
  0 \le i \le 2, \\
 0 \le j \le 2,\\
 0 \le k \le 1
 \end{matrix}
 }$
  consists of $18$~distinct elements of~$G$, and
 \item either
 \begin{enumerate}
 \renewcommand{\theenumi}{\alph{enumi}}
 \item \label{P3xP3xP2-G=Y} $G = Y$, or 
 \item \label{P3xP3xP2-hamcyc}
 there is a hamiltonian cycle~$C$ of $X \setminus Y$, such that $C$~contains
the edge $[t^3](s)$.
 \end{enumerate}
 \end{itemize}
 Then $\hamilton$ contains every basic $4$-cycle.
 \end{prop}

Before proving the proposition, let us establish a lemma that describes the
hamiltonian cycles we will construct.

\begin{lem} \label{lots4cycequiv}
 Suppose 
 \begin{enumerate}
 \renewcommand{\theenumi}{\alph{enumi}}
 \item $S = \{s^{\pm1}, t^{\pm1}, u^{\pm1} \} $,
 \item \label{lots4cycequiv-odd}
 some element of~$\hamilton$ is the sum of an odd number of basic $4$-cycles,
 and 
 \item \label{lots4cycequiv-hamcyc}
 there exist oriented hamiltonian cycles $H_1,\ldots,H_6$ in~$X$, that contain
the following specific oriented paths:
 \begin{enumerate}
 \renewcommand{\theenumii}{{$H_{\arabic{enumii}}$}}
 \renewcommand{\labelenumii}{\theenumii:}
 \item $[v_1](t, u^{-1}, t^{-1})$ and $[v_1 t s^{-1}](u^{-1})$, for some $v_1 \in
X$;
 \item $[v_2](t, u^{-1}, t^{-1})$ and $[v_2 t^2](u^{-1})$, for some $v_2 \in X$;
 \item $[v_3](t, u^{-1}, t^{-1})$ and $[v_3 t s](u^{-1})$, for some $v_3 \in X$;
 \item $[v_4](u^{-1}, t^{-1}, u)$ and $[v_4 u^{-1} s](t^{-1})$, for some $v_4 \in
X$;
 \item $[v_5](t, s, t^{-1})$ and $[v_5 t s u^{-1}](s^{-1})$, for some $v_5 \in X$;
 \item $[v_6](u, t^{-1}, u^{-1})$ and $[v_6  u t^{-1} s](t)$, for some $v_6 \in
X$.
 \end{enumerate}
 \end{enumerate}
 Then $\hamilton$ contains every basic $4$-cycle.
 \end{lem}

\begin{proof}
 Let us begin by establishing that it suffices to show
 \begin{enumerate}
 \renewcommand{\theenumi}{\roman{enumi}}
 \item \label{lots4cycequiv-2Q}
 $2Q \in \hamilton$, for some basic $4$-cycle~$Q$; and
 \item \label{lots4cycequiv-equiv}
 for any two basic $4$-cycles $Q_1$ and~$Q_2$, we have $Q_1 \equiv \pm Q_2
\pmod{\hamilton}$.
 \end{enumerate}
 From \pref{lots4cycequiv-2Q} and~\pref{lots4cycequiv-equiv}, we see that any
even multiple of any basic $4$-cycle belongs to~$\hamilton$.
From~\pref{lots4cycequiv-odd} and~\pref{lots4cycequiv-equiv}, we see that some
odd multiple of any basic $4$-cycle belongs to~$\hamilton$. Subtracting an
appropriate even multiple, we conclude that the $4$-cycle itself belongs
to~$\hamilton$, as desired. Thus, \pref{lots4cycequiv-2Q} and
\pref{lots4cycequiv-equiv} do indeed suffice.

 By applying Lemma~\fullref{4c}{s} to $H_1$, $H_2$, $H_3$, and~$H_4$, we
see that
 \begin{enumerate}
 \item \label{lots4cycequiv-H1}
 $(t, u^{-1}, t^{-1}, u) + [t](s^{-1}, u^{-1}, s, u) \in \hamilton$;
 \item \label{lots4cycequiv-H2}
 $(t, u^{-1}, t^{-1}, u) + [t](t, u^{-1}, t^{-1}, u) \in \hamilton$;
 \item \label{lots4cycequiv-H3}
 $(t, u^{-1}, t^{-1}, u) + [t](s, u^{-1}, s^{-1}, u) \in \hamilton$;
 and
 \item \label{lots4cycequiv-H4}
 $(u^{-1}, t^{-1}, u, t) + [u^{-1}](s, t^{-1}, s^{-1}, t) \in \hamilton$.
 \end{enumerate}
 By applying Lemma~\fullref{4c}{d} to $H_5$ and~$H_6$, we see that
 \begin{enumerate}
 \setcounter{enumi}{4}
 \item \label{lots4cycequiv-H5}
 $(t, s, t^{-1}, s^{-1}) - [t](u^{-1}, s, u, s^{-1}) \in \hamilton$; and
 \item \label{lots4cycequiv-H6}
 $(u, t^{-1}, u^{-1}, t) - [u](s, t^{-1}, s^{-1}, t) \in \hamilton$.
 \end{enumerate}

 To verify \pref{lots4cycequiv-2Q}, note that, because
 \begin{align*}
 [u] (s,t^{-1},s^{-1},t) &= [t^{-1} u](t, s, t^{-1}, s^{-1}) , \\
 [u](u^{-1}, s, u,s^{-1}) &= - [u](s, u^{-1}, s^{-1}, u) , \\
 \intertext{and}
 [t^{-1} u](t, u^{-1}, t^{-1}, u) &= -(u,t^{-1},u^{-1}, t) ,
 \end{align*}
 we have
 \begin{align*}
 2 (u, t^{-1}, u^{-1}, t)
 &= \bigl( (u, t^{-1}, u^{-1}, t) - [u](s, t^{-1}, s^{-1}, t) \bigr)
 && \pref{lots4cycequiv-H6}
 \\ & \qquad
 + \bigl( [t^{-1} u](t, s, t^{-1}, s^{-1}) - [u](u^{-1}, s, u, s^{-1}) \bigr) 
 && \pref{lots4cycequiv-H5}
 \\ & \qquad
 - \bigl( [u](s, u^{-1}, s^{-1}, u) + [t^{-1} u](t, u^{-1}, t^{-1}, u) \bigr)
 && \pref{lots4cycequiv-H3}
 \\
 & \in \hamilton
 .
 \end{align*}

We now establish~\pref{lots4cycequiv-equiv}. Let 
 $$Q_0 = (t, u^{-1}, t^{-1}, u) .$$
 Given any basic $4$-cycle~$Q$ in~$X$, we will show that $Q \equiv \pm Q_0
\pmod{\hamilton}$. 

From \pref{lots4cycequiv-H1} and~\pref{lots4cycequiv-H4}, we see that any basic
$4$-cycle of the form
 $$\mbox{$[v](s, u, s^{-1}, u^{-1})$ 
 \qquad
  or
 \qquad $[v](s, t, s^{-1}, t^{-1})$}$$
 is congruent, modulo~$\hamilton$, to a translate of~$\pm Q_0$. Thus, we may
assume $Q = [w]Q_0$, for some $w \in G$, so it suffices to show
 \begin{enumerate}
 \setcounter{enumi}{18}  
 \renewcommand{\theenumi}{\Alph{enumi}}
 \item \label{lots4cycequiv-S}
 $Q_0 \equiv \pm[s]Q_0 \pmod{\hamilton}$;
 \item \label{lots4cycequiv-T}
 $Q_0 \equiv \pm[t]Q_0 \pmod{\hamilton}$; and
 \item \label{lots4cycequiv-U}
 $Q_0 \equiv \pm[u]Q_0 \pmod{\hamilton}$.
 \end{enumerate}
 First, note that \pref{lots4cycequiv-T} is immediate
from~\pref{lots4cycequiv-H2}. 

For \pref{lots4cycequiv-S}, noting that
 $$ [t](s, u^{-1}, s^{-1}, u) = -[s t](s^{-1}, u^{-1}, s, u) ,$$
 we have
 \begin{align*}
 Q_0 + [s]Q_0
 &=
 \bigl( (t, u^{-1}, t^{-1}, u) + [t](s, u^{-1}, s^{-1}, u) \bigr)
 && \pref{lots4cycequiv-H3} \\
 & \qquad +
 \bigl( [s](t, u^{-1}, t^{-1}, u) + [s t](s^{-1}, u^{-1}, s, u) \bigr)
 && \pref{lots4cycequiv-H1} \\
 &\in \hamilton
 .
 \end{align*}
 For \pref{lots4cycequiv-U}, noting that
 \begin{align*}
 -Q_0 &= [t u^{-1}](u,t^{-1},u^{-1},t) \\
 \intertext{and}
 [u]Q_0 &= [t u](u^{-1}, t^{-1}, u, t),
 \end{align*}
 we have
 \begin{align*}
 -Q_0 + [u]Q_0
 &=
 \bigl( [t u^{-1}](u, t^{-1}, u^{-1}, t) - [t](s, t^{-1}, s^{-1}, t) \bigr)
 && \pref{lots4cycequiv-H6} \\
 & \qquad +
 \bigl( [t u](u^{-1}, t^{-1}, u, t) + [t](s, t^{-1}, s^{-1}, t) \bigr)
 && \pref{lots4cycequiv-H4} \\
 &\in \hamilton
 .
 \end{align*}
 This completes the proof.
 \end{proof}

\begin{proof}[\bf Proof of Proposition~\ref{P3xP3xP2}]
 It suffices to verify hypotheses \pref{lots4cycequiv-odd}
and~\pref{lots4cycequiv-hamcyc} of Lemma~\ref{lots4cycequiv}.

 Let 
 $$ \Ho = (u, s, t, u^{-1}, t^{-1}, s, u, t, u^{-1}, t, u, s^{-1}, u^{-1},
s^{-1}, u, t^{-1}, u^{-1}, t^{-1}) $$
 and
 $$ \He = (s^2, t^2, s^{-2}, u, t^{-1}, u^{-1}, s, u, t, s, t^{-2}, s^{-2},
u^{-1})
 ,$$
 so $\Ho$ and $\He$ are hamiltonian cycles in the subgraph of~$X$ induced by~$Y$.

Let
 $$C' = 
 \begin{cases}
 \hfil 0 & \mbox{if $G = Y$} , \\
 C + [t^2](s, t, s^{-1}, t^{-1}) & \mbox{if $G \neq Y$}
 ,
 \end{cases}
 $$
 and let 
 $$ \mbox{$\Ho' = \Ho + C'$ and $\He' = \He + C'$.}$$
 Because $\Ho$ and~$\He$ each contain the oriented edge $[s t^2](s^{-1})$, and
$C$~contains the oriented edge $[t^3](s)$ (if $G \neq Y$), we see that $\Ho'$
and~$\He'$ are hamiltonian cycles in~$X$.

\pref{lots4cycequiv-odd}
 It is easy to verify that $\Ho$ is the sum of 9 basic $4$-cycles, and that $\He$
is the sum of 8 basic $4$-cycles. Thus, $\Ho - \He$ is the sum of 17 basic
$4$-cycles. Also, we have $\Ho - \He = \Ho' - \He' \in \hamilton$. Therefore,
$\Ho - \He$ is an element of~$\hamilton$ that is the sum of an odd number of
basic $4$-cycles.

\pref{lots4cycequiv-hamcyc} 
 Let
 $$ \mbox{$H_1 = H_2 = H_3 = H_4 = \Ho'$, \qquad $H_5 = H_6 = \He'$,} $$
 $$ v_1 = v_2 = v_3 = s u,
 \qquad
 v_4 = t u, 
 \qquad
 v_5 = s t u,
 \qquad
 v_6 = t^2 
 .
 $$
 Then, because $\Ho'$ and~$\He'$ contain every edge of $\Ho$ or~$\He$, except
(possibly) $[s t^2](s^{-1})$, it is easy to verify that each $H_i$ contains the
oriented paths specified in~\fullref{lots4cycequiv}{hamcyc}.
 \end{proof}

\begin{cor} \label{PaxPbxPcHas4Cycles}
 If
 \begin{itemize}
 \item $S = \{s^{\pm1}, t^{\pm1}, u^{\pm1} \}$, and 
 \item there exist $m,n \ge 3$ and $p \ge 2$, such that every element~$g$ of~$G$
can be written \emph{uniquely} in the form $g = s^i t^j u^k$ with $0 \le i \le
m-1$, $0 \le j \le n-1$, and $0 \le k \le p-1$,
 \end{itemize}
 then $\hamilton$ contains every basic $4$-cycle.
 \end{cor}

\begin{proof}
 By permuting $s$, $t$, and~$u$, we may assume $n = \max\{m,n,p\}$.
 Because $|G|$ is even \see{|G|even}, we cannot have $m = n = p = 3$. Thus,
either $m = n = 3$ and $p = 2$, or $n \ge 4$. In the former case, the desired
conclusion is immediate from Proposition~\fullref{P3xP3xP2}{G=Y}. Thus, we
henceforth assume that $n \ge 4$. We will construct an appropriate hamiltonian
cycle~$C$ as specified in Proposition~\fullref{P3xP3xP2}{hamcyc}.

 For natural numbers 
$a \le m$, $b \le n$, and $c \le p$, and $x \in \{s^{\pm 1}\}$ 
and $y \in \{t^{\pm 1}\}$, let
 \begin{align*}
 A_{a,b}(x,y) &=
 \begin{cases}
 (x^{a-1}, u, x^{-(a-1)}) & \mbox{if $b$ is odd}, \\
 (x^{a-1},y, x^{-(a-1)}, u, x^{a-1}, y^{-1}, x^{-(a-1)}) & \mbox{if $b$ is even},
 \end{cases}
 \\
 C_{a,b}(x,y) &=
 \Bigl( \bigl(x^{a-1}, y, x^{-(a-1)}, y \bigr)^{\lfloor (b-1)/2 \rfloor},
 A_{a,b}(x,y),
 \\ 
 &\hskip1in
 \bigl( y^{-1}, x^{a-1}, y^{-1}, x^{-(a-1)} \bigr)^{\lfloor (b-1)/2 \rfloor},
 u^{-1} \Bigr)
 \\
 \intertext{and}
 X_{a,b,c}(x,y) &= \bigset{
 x^i y^j u^k 
 }{
 \begin{matrix}
  0 \le i \le a-1, \\
  0 \le j \le b-1,\\
  0 \le k \le c-1
 \end{matrix}
 }
 .
 \end{align*}
 Then $C_{a,b}(x,y)$ is a hamiltonian cycle of the subgraph of~$X$ induced by 
 $X_{a,b,2}(x,y)$.

 Let 
 $$ A = 
 \begin{cases}
 -[s^2 t^3]C_{3,n-3}(s^{-1},t)
 & \mbox{if $m = 3$} , \\
 \\
 \begin{matrix}
 -[s^2 t^3]C_{3,n-3}(s^{-1},t) + [s^2 t^3](s,u,s^{-1},u^{-1}) \hfill \\
 \qquad {} + [s^3 t^3] C_{m-3,n-3}(s,t) \hfill \\
 \qquad {} + [s^3 t^3](t^{-1},s,t,s^{-1}) - [s^3 t^2]C_{m-3,3}(s,t^{-1})
 \hfill
 \end{matrix}
 & \mbox{if $m > 3$}
 ,
 \end{cases}
 $$
 so $A$ is a hamiltonian cycle of the subgraph of~$X$ induced by 
 $$ X_{m,n,2}(s,t) \setminus X_{3,3,2}(s,t)
 ,$$
 and $A$ contains the oriented edges $[t^3](s)$ and $[s t^3 u](s^{-1})$.

If $p > 2$, let $B$ be any hamiltonian cycle of the subgraph of~$X$ induced by
 $[u^2] X_{m,n,p-2}(s,t)$,
 such that $B$ contains the oriented edge $[t^3 u^2](s)$. For example, if $p$~is
even, we may let 
 $$ B = -\sum_{i=1}^{(p-2)/2} [u^{2i}] C_{m,n}(s,t)  + \sum_{i=1}^{(p-4)/2} [t^3
u^{2i+1}] (s,u,s^{-1},u^{-1})  .$$
 On the other hand, if $p$~is odd, then $mn$ must be even, so it is easy to
construct a hamiltonian cycle~$H$ of the subgraph of~$X$ induced by
$X_{m,n,1}(s,t)$, such that $H$~contains the oriented edge $[t^3](s)$. If
$[v](x)$ is any other oriented edge of~$H$, then we may let
 \begin{align*}
  B = \sum_{i=2}^{p-1} (-1)^i [u^i] H 
  &+ \sum_{i=1}^{(p-3)/2} [v u^{2i}] (u, x,u^{-1},x^{-1})
 \\ 
  &+ \sum_{i=1}^{(p-3)/2} [t^3 u^{2i+1}] (s,u,s^{-1},u^{-1}) 
 .
 \end{align*}

Let 
 $$C = 
 \begin{cases}
 \hfil A & \mbox{if $p = 2$} , \\
 A + [t^3 u](s,u,s^{-1},u^{-1})
 + B
 & \mbox{if $p > 2$}
 .
 \end{cases}
 $$
 Then $C$ is a hamiltonian cycle of $X \setminus X_{3,3,2}(s,t)$, and $C$
contains the oriented edge $[t^3](s)$, so Proposition~\fullref{P3xP3xP2}{hamcyc}
implies that $\hamilton$ contains every basic $4$-cycle.
 \end{proof}

\begin{cor}  \label{GridsHave4Cycles}
 If 
 \begin{itemize}
 \item $S = \{s^{\pm1}, t^{\pm1}, u^{\pm1} \}$,
 \item $|S| = 5$ or~$6$,
 \item $G' \neq G$,
 \item $\langle t \rangle \neq G'$, and
 \item either
 \begin{enumerate}
 \renewcommand{\theenumi}{\alph{enumi}}
 \item \label{GridsHave4Cycles-not4}
 $|G|$ is \emph{not} divisible by~$4$, or
 \item \label{GridsHave4Cycles-no2}
  no more than one of $|t|$, $|G'|/|t|$ and $|G/G'|$ is equal to~$2$,
 \end{enumerate}
 \end{itemize}
 then $\hamilton$ contains every basic $4$-cycle.
 \end{cor}

\begin{proof}
 Because the product of the three orders in~\pref{GridsHave4Cycles-no2} is~$|G|$,
it is clear that \pref{GridsHave4Cycles-not4} implies
\pref{GridsHave4Cycles-no2}. Thus, we may assume \pref{GridsHave4Cycles-no2}.

 Let $m = |G/G'|$, $n = |t|$, and $p = |G'|/|t|$. Then every element of~$G$ can
be written uniquely in the form
 $s^i t^j u^k$, with $0 \le i \le m-1$, $0 \le j \le n-1$, $0 \le k \le p-1$.
Thus, Corollary~\ref{PaxPbxPcHas4Cycles} applies (with $s$ or~$t$ in the role
of~$u$, if $m = 2$ or $n = 2$).
 \end{proof}

\begin{cor} \label{4DGridsHave4Cycles}
 If
 \begin{itemize}
 \item $|S| \ge 5$,
 \item $S$ is irredundant, and
 \item $S$ contains no more than one involution,
 \end{itemize}
 then $\hamilton$ contains every basic $4$-cycle $(s,t,s^{-1}, t^{-1})$.
 \end{cor}

\begin{proof}
 We may assume $s \neq t^{\pm 1}$. (Otherwise, the $4$-cycle is degenerate, so 
 $(s,t,s^{-1}, t^{-1}) = 0 \in \hamilton$.)
 By an argument similar to the one involving~$\Gamma$ in Case~\ref{4cycpf-highdeg}
of the proof of Corollary~\ref{4cyc}, we may assume $|S|$~is either~$5$ or~$6$.
Also, since $S$ contains no more than one involution, we may assume $|s| \ge 3$ 
(by interchanging $s$ and~$t$ if necessary).
 The desired conclusion follows from either Corollary~\fullref{4cyc}{|G/G'|>2}
(if $|G|$ is divisible by~$4$), or Corollary~\fullref{GridsHave4Cycles}{not4} (if
$|G|$ is not divisible by~$4$).
 \end{proof}

\begin{cor} \label{K2xK3xK3}
 If $X \iso K_3 \boxprod K_3 \boxprod K_2$, then $\hamilton=\even$.
 \end{cor}

\begin{proof}
 Corollary~\ref{4DGridsHave4Cycles} implies that $\hamilton$ contains every basic
$4$-cycle. Let $s$~be the involution in~$S$, so $G' \iso \integer_3 \times
\integer_3$. We have
 \begin{align*}
 \even &\subseteq \hamilton + 2\flow' && \mbox{\fullsee{E=H+2F'}{X'notbip}} \\
 &\subseteq \hamilton + 2\hamilton' && \mbox{\see{H+C=F(K3xK3)}} \\
 &\subseteq \hamilton && \mbox{\see{2H'inH(G'odd)deg5}} 
 ,
 \end{align*}
 as desired.
 \end{proof}

\begin{cor} \label{OddHeight4Cycles}
 If
 \begin{itemize}
 \item $S = \{s^{\pm 1}, t^{\pm 1}, u^{\pm 1} \}$, with $|S| = 6$,
 \item $G' \neq G$, and
 \item $|G/G'|$ is odd,
 \end{itemize}
 then $\hamilton$ contains every basic $4$-cycle.
 \end{cor}

\begin{proof}
 We may assume $\langle t \rangle = \langle u \rangle$, for otherwise
Corollary~\fullref{GridsHave4Cycles}{no2} applies (perhaps after interchanging
$t$ and~$u$), because $|t| \neq 2$, $|u| \neq 2$, and $|G/G'| \neq 2$. Hence
$\langle t \rangle = G'$. In particular, we know $|t| = |G'|$ is even (because
$|G/G'|$ is odd).

We have $u = t^{-q}$ for some~$q$. (Notice the negative sign in the exponent.) We
may assume $2 < q < |t|/2$, by replacing $u$~with its inverse if necessary (and
noting that $q$ is not equal to~$2$ or $|t|/2$, since $|u| = |t|$, and $|t|$~is
even). Then
 $$ Y = \bigset{
 s^i t^j u^k 
}{
 \begin{matrix} 0 \le i \le 2, \\ 0 \le j \le 2, \\ 0 \le k \le 1 \end{matrix} }$$
 consists of 18 distinct elements of~$G$.

Let $m = |G/G'|$ (by assumption, $m$~is odd), and let  $C = C_1 + C_2 + C_3$,
where
 \begin{align*}
 C_1 &= [t^3]\Bigl( s, \bigl( s^{m-2}, t, s^{-(m-2)}, t \bigr)^{(|t|-2q)/2}
\sharp, s^{-1}, t^{-(|t|-2q-1)} \Bigr)
 , \\
 C_2 &=
 \begin{cases}
 \hfil 0 & \mbox{if $m = 3$} , \\
 \begin{matrix}
 [s^{m-1} t^2] \Bigl( t, s^{-1}, t^{-1},
 \hfill \\ \hfill
 \hskip1in
 \bigl( s^{-(m-5)}, t^{-1}, s^{m-5}, t^{-1} \bigr)^q \sharp, 
 s, t^{2q-1} \Bigr)
 \end{matrix}
 & \mbox{if $m \ge 5$}
 ,
 \end{cases}
 \\
 C_3 &=
 \begin{cases}
 \hfil 0 & \mbox{if $q = 3$} , \\
 \begin{matrix}
 [t^{-(2q-3)}] \Bigl( t^{-1}, s, t^{q-3}, s, t^{-(q-4)}, u^{-1}, t^{q-4}, 
 \hfill \qquad \\ \hfill
 s^{-1}, t^{-(q-4)},s^{-1}, t^{q-4}, u, t^{-(q-4)} \Bigr)
 \end{matrix}
 & \mbox{if $q > 3$}
 .
 \end{cases}
 \end{align*}
 Then $C$ is a hamiltonian cycle of $X \setminus Y$, and $C$~contains the
oriented edge $[t^3](s)$, so Proposition~\fullref{P3xP3xP2}{hamcyc} applies.
 \end{proof}

\section{Redundant generators in graphs of large degree} \label{RedundantSect}

In this section, assuming that $|S| \ge 5$, and that $s$~is a redundant
generator, we provide an induction step: if $\hamilton' = \even'$, then
$\hamilton = \even$ \see{G=G'big}. 

\begin{lem} \label{H+H=2s}
 If there exist $s, t \in S$ 
and oriented hamiltonian cycles $H_+$ and $H_-$ of $X$ 
such that
 \begin{itemize}
 \item $G' = G$,
 \item $\even' \subseteq \hamilton + \hamilton'$,
 \item $|S| \ge 4$,
 \item $H_+ = (s) + [t](s) + (\text{\upshape edges in~$X'$})$,
 and
 \item $H_- = (s) + [st](s^{-1}) + (\text{\upshape edges in~$X'$})$,
 \end{itemize}
 then $\hamilton = \even$.
 \end{lem}
 
\begin{proof}
 We have 
 $$ H_+ + H_- = 2(s) + \text{(edges in~$X'$)} ,$$
 so adding or subtracting appropriate translates of $H_+ + H_-$ will eliminate
all the $s$-edges from any flow in $2 \flow$, the result being a flow
in~$\even'$. Therefore
 \begin{equation} \label{H+H=2s-2FinH+E'}
 2\flow \subseteq \hamilton + \even'
 .
 \end{equation}
 We have
 \begin{align*}
 \even &\subseteq \hamilton + 2 \flow && \mbox{\see{ALWthm}} \\
 &\subseteq \hamilton + \even' && \mbox{\see{H+H=2s-2FinH+E'}} \\
 &\subseteq \hamilton + \hamilton' && \mbox{(by assumption)} \\
 &\subseteq \hamilton && \mbox{\see{G=G'->H'inH}} 
 ,
 \end{align*}
 as desired.
 \end{proof}

\begin{lem} \label{lotsofchords}
 If there exist $s$ and~$t$ in~$S$, such that
 \begin{itemize}
 \item $|S| \ge 5$, 
 \item $t \in S'$, 
 \item $G = \langle t \rangle$, and
 \item $\even' \subseteq \hamilton + \hamilton'$,
 \end{itemize}
 then $\hamilton = \even$.
 \end{lem}

\begin{proof}
 Choose $u \in S \setminus \{s^{\pm1},t^{\pm1}\}$. Write $s = t^p$ and $u =
t^q$. We may assume $2 \le p,q \le |t|/2$ (by replacing~$s$ or~$u$ with its
inverse, if necessary).
 We may also assume $X'$ is bipartite, but $X$~is not bipartite, for otherwise we
have
 \begin{align*}
 \even &\subseteq \hamilton + \even' && \mbox{\see{E=H+E'(G'=G)}} \\
 &\subseteq \hamilton + \hamilton' && \mbox{(by assumption)} \\
 &\subseteq \hamilton && \mbox{\see{G=G'->H'inH}}
 ,
 \end{align*}
 as desired. Hence, $p$~is even and $q$~is odd.

 Let 
 $$ H_+ = (s, t^{-(p-1)}, s, t^{|t|-p-1}) ,$$
 so $H_+$ is a hamiltonian cycle in~$X$, and the only oriented $s$-edges in $H_+$
are $(s)$ and $[t](s)$.

 Define $H_-$ from the hamiltonian cycle
 $$ H_-^* = [su]\Bigl( \bigl( u^{-1}, t, u, t \bigr)^{(q-1)/2}, t^{|t|-2q}, u, t
\Bigr) , $$
 by replacing
 \begin{itemize}
 \item the path $[su](u^{-1}, t, u)$ with the oriented edge $[su](t)$ and
 \item the oriented edge $(t)$ with the path $(s,t,s^{-1})$.
 \end{itemize}
 Then $H_-$ is a hamiltonian cycle of~$X$, and the only oriented $s$-edges
in~$H_-$ are $(s)$ and $[st](s^{-1})$.

 Therefore, Lemma~\ref{H+H=2s} applies.
 \end{proof}

\begin{prop} \label{G=G'big}
 If there exists $s \in S$, such that
 \begin{itemize}
 \item $|S| \ge 5$,
 \item $G' = G$,
 \item not every nonidentity element of~$G$ is an involution,
 and
 \item $\even' \subseteq \hamilton + \hamilton'$,
 \end{itemize}
 then $\hamilton = \even$.
 \end{prop}

\begin{proof}
 We may assume that $X'$ is bipartite, but $X$~is not bipartite, for otherwise we
have
 \begin{align*}
 \even &\subseteq \hamilton + \even' && \mbox{\see{E=H+E'(G'=G)}} \\
 &\subseteq \hamilton + \hamilton' && \mbox{(by assumption)} \\
 &\subseteq \hamilton && \mbox{\see{G=G'->H'inH}}
 ,
 \end{align*}
 as desired.

 By assumption, not every nonidentity element of~$G$ is an involution. Because
$S'$ generates~$G$, this implies that we may choose some $t \in S'$ that is not an
involution.  It suffices to find oriented hamiltonian cycles $H_+$ and~$H_-$
in~$X$, such that 
 \begin{itemize}
 \item the only $s$-edges in~$H_+$ are $(s)$ and~$[t](s)$, and
 \item  the only $s$-edges in~$H_-$ are $(s)$ and~$[st](s^{-1})$,
 \end{itemize}
 for then Lemma~\ref{H+H=2s} applies.

Note that, because $X'$~is bipartite, $|t|$~must be even. Also note that we may
assume $\langle t \rangle \neq G$, for otherwise Lemma~\ref{lotsofchords}
applies. 
 Let $m = |G|/|t|$, and let $(u_1,u_2,\ldots,u_{m-1})$ be a hamiltonian path in
 $$ \Cay \bigl( G/ \langle t \rangle ; S' \setminus \{t^{\pm1}\} \bigr) 
 . $$
 We have $s = t^p u_1 u_2 \cdots u_q$, for some $p$ and~$q$.  We may assume $0
\le p \le |t|/2$ and $0 \le q \le m/2$ (by replacing $s$ or $t$ with its
inverse, if necessary). Because $X'$ is bipartite, but $X$~is not bipartite, we
know that $p+q$ is even.

\setcounter{proofcase}{0}

\begin{proofcase}
 Assume $m > 2$, and $p$ and~$q$ are odd.
 \end{proofcase}
 Let
 \begin{align*}
  H_+ &= \Bigl( s,
 \bigl((u_{q+1-i}^{-1})_{i=1}^q, t^{-1}, (u_i)_{i=1}^q, t^{-1} \bigr)^{(p-1)/2},
 (u_{q+1-i}^{-1})_{i=1}^q, \\
 &\hskip1in
 s, \bigl( (u_{q+1-i}^{-1})_{i=1}^q, t, (u_i)_{i=1}^q, t \bigr)^{(|t|-p -
1)/2}\sharp,
 u_{q+1},
 \\ &\hskip1in
 \bigl( (u_{q+1+i})_{i=1}^{m-q-2}, t^{-1}, (u_{m-i}^{-1})_{i=1}^{m-q-2},t^{-1}
\bigr)^{|t|/2}\sharp,
 \\ &
 \hskip3in 
 (u_{q+2-i}^{-1})_{i=1}^{q+1} \Bigr)
 . 
 \end{align*}
 Define $H_-$ from the hamiltonian cycle
 \begin{align*}
 H_-^* &= \Bigl( t,
 \bigl( t^{|t|-2}, u_{2i-1}, t^{-(|t|-2)}, u_{2i} \bigr)_{i=1}^{(q-1)/2},
 t^{|t|-2}, u_{q},
 \\ &\hskip1in
 \bigl( (u_{q+i})_{i=1}^{m-q-1}, t^{-1}, (u_{m-i}^{-1})_{i=1}^{m-q-1}, t^{-1}
\bigr)^{|t|/2}\sharp,
 \\ &\hskip3in 
 (u_{q+1-i}^{-1})_{i=1}^{q} \Bigr),
  \end{align*}
 by replacing
 \begin{itemize}
 \item the path $[stu_{q+1}](u_{q+1}^{-1}, t^{-1}, u_{q+1})$ with the edge
$[stu_{q+1}](t^{-1})$, and 
 \item the edge $(t)$ with the path $(s,t,s^{-1})$.
 \end{itemize}

\begin{proofcase}
 Assume $m > 2$, and $p$ and~$q$ are even.
 \end{proofcase}
 Let $C = 0$ if $m = 4$ and $q = 2$; otherwise, let
 $$
 C = 
 [u_1 u_2 \cdots u_{q+1}]\Bigl( u_{q+2}, 
 \bigl( (u_{q+2+i})_{i=1}^{m-q-3}, t, (u_{m-i}^{-1})_{i=1}^{m-q-3}, t
\bigr)^{|t|/2}\sharp, u_{q+2}^{-1}, t \Bigr).  $$
 Define
 \begin{align*}
  H_+ &=
 \biggl( s, \Bigl( t^{-p}, \bigl( u_{q+2-2i}^{-1}, t^{p+1}, u_{q+1-2i}^{-1},
t^{-(p+1)} \bigr)_{i=1}^{q/2} \Bigr)\sharp,
 s, u_{q+1}, t^{-(p+2)}, 
 \\
 &\hskip1in
 \bigl( t^{-(|t|-p-3)}, u_{q+3-2i}^{-1}, t^{|t|-p-3}, u_{q+2-2i}^{-1}, 
\bigr)_{i=1}^{(q+2)/2} \sharp,
 t
 \biggr)
 + C
 .
 \end{align*}

 If $p \neq 0$, define $H_-$ from the hamiltonian cycle
 \begin{align*}
 H_-^* &= \Bigl( t,
 \bigl( t^{|t|-2}, u_{2i-1}, t^{-(|t|-2)}, u_{2i} \bigr)_{i=1}^{q/2}, \\
 &\hskip1in
 \bigl( (u_{q+i})_{i=1}^{m-q-1}, t, (u_{m-i}^{-1})_{i=1}^{m-q-1}, t
\bigr)^{|t|/2}\sharp, (u_{q+1-i}^{-1})_{i=1}^{q} \Bigr),
  \end{align*}
 by replacing 
 \begin{itemize}
 \item the path $[su_{q+1}](u_{q+1}^{-1}, t, u_{q+1})$ with the oriented edge
$[su_{q+1}](t)$, and 
 \item the edge $(t)$ with the path $(s,t,s^{-1})$.
 \end{itemize}

If $p = 0$, let
 \begin{align*}
  H_- &= \Bigl( s, u_{q+1}, \bigl( (u_{q+1+i})_{i=1}^{m-q-2}, t,
(u_{m-i}^{-1})_{i=1}^{m-q-2}, t \bigr)^{|t|/2}\sharp,
 \\ &\hskip1in
 u_{q+1}^{-1}, t^{-(|t|-3)}, u_q^{-1},
  \bigl( (u_{q-i}^{-1})_{i=1}^{q-1}, t, (u_i)_{i=1}^{q-1}, t \bigr)^{(|t|-2)/2},
 \\ &\hskip1in
 t, u_q, s^{-1}, (u_i)_{i=1}^{q-2}, t^{-1}, (u_{q-1-i}^{-1})_{i=1}^{q-2} \Bigr) .
 \end{align*}

\begin{proofcase}
 Assume $m = 2$, and $p$ and~$q$ are odd.
 \end{proofcase}
 Note that $q = 1$. Let
 $$
  H_+ = \Bigl( s, \bigl(u_1^{-1}, t^{-1}, u_1, t^{-1} \bigr)^{(p-1)/2},
 u_1^{-1}, s,
 \bigl( u_1^{-1}, t, u_1, t \bigr)^{(|t|-p-1)/2}, u_1^{-1} \Bigr)
 $$
 and 
 $$H_- = \Bigl( s,t^{-(|t|-1)},s^{-1},t^{|t|-1} \Bigr) .$$

\begin{proofcase}
 Assume $m = 2$, and $p$ and~$q$ are even.
 \end{proofcase}
 Note that $q = 0$, so $s \in \langle t \rangle$. 
 Let
 \begin{align*}
 H_+ &= \Bigl( s, t^{-(p-1)}, s, t^{|t|-p-2}, u_1, t^{-(|t|-1)}, u_1^{-1} \Bigr)
 ,\\
 H_- &= \Bigl( s, t^{-(p-2)}, u_1, t^{|t|-3}, u_1^{-1}, t^{-(|t|-p-2)}, s^{-1},
u_1, t^{-1}, u_1^{-1} \Bigr)
 .
 \end{align*}
 \end{proof}

To deal with the groups in which every element is an involution, we prove the
following complement to Proposition~\ref{G=G'big}:

\begin{lem} \label{G=G'invols}
 If
 \begin{itemize}
 \item $|S| \ge 5$,
 \item $S$ is redundant,
 \item every nonidentity element of~$G$ is an involution,
 and
 \item $\hamilton' = \even'$, whenever $s$ is a redundant generator in~$S$,
 \end{itemize}
 then $\hamilton = \even$.
 \end{lem}

\begin{proof}
 Let $s$ be a redundant generator in~$S$. We may assume $S'$ is redundant, for
otherwise Lemma~\ref{notallinvols} shows that we may apply
Proposition~\ref{G=G'big}, by replacing $X$ with an isomorphic Cayley graph. Thus,
we may let $t$ be a redundant generator in~$S'$. We may assume (by interchanging
$s$ and $t$, if necessary), that either $X$ is bipartite or $X'$ is not bipartite.
Therefore
  \begin{align*}
 \even &\subseteq \hamilton + \even' && \mbox{\see{E=H+E'(G'=G)}} \\
 &\subseteq \hamilton + \hamilton' && \mbox{(because $\hamilton' = \even'$)} \\
 &\subseteq \hamilton && \mbox{(because $G' = G$)}
 . \end{align*}
 \end{proof}

\section{Two troublesome cases} \label{TroubleSect}

In this section, we treat two special cases
\seeand{EvenHeightOverChords}{cycle+diam+chords}, in order to deal with some
graphs that are not covered by our previous results and do not yield easily to a
proof by induction.

\begin{lem} \label{PrismOverChords4cyc}
 If there exist $s,t \in S$, such that
 \begin{itemize}
 \item $|G/G'| = 2$,
 \item $G' = \langle t \rangle$, and
 \item $|S'| = 4$,
 \end{itemize}
 then $\hamilton$ contains the basic $4$-cycle $(s,t,s^{-1},t^{-1})$.
 \end{lem}

\begin{proof}
 Let $u \in S \setminus \{s^{\pm1}, t^{\pm1}\}$, so $S = \{s^{\pm1},
t^{\pm1}, u^{\pm1}\}$.
 We may write $u=t^q$ for some $q$.  Interchange $u$ with $u^{-1}$, if necessary,
to ensure that 
 \begin{equation*} 
 \mbox{$q$ is even if $|t|$ is odd.}
 \end{equation*}
 Define hamiltonian cycles
 $$
 H_+= \bigl( s, u, t^{-(q-2)}, s^{-1}, t^{|t|-3}, s, t^{-(|t|-q-2)}, u^{-1},
s^{-1}, t^{-1} \bigr)
 $$
 and
 $$H_-= \bigl( t^{|t|-1}, s, t^{-(|t|-1)}, s^{-1} \bigr),$$
 and let 
 $$ Q = (t,s,t^{-1},s^{-1}) ,$$
 so $-Q$ is the basic $4$-cycle specified in the statement of the lemma.
 Then 
 \begin{itemize}
 \item $H_+$ contains both the oriented path $[t^{-2}](t,s,t^{-1})$ and the
oriented edge $[e](s)$, and
 \item $H_-$ contains both the oriented path $[t^{-2}](t,s,t^{-1})$ and the
oriented edge $[s](s^{-1})$,
 \end{itemize}
 so Lemma~\ref{4c}(\ref{4c-s},\ref{4c-d},\ref{4c-2})
 (with $x = t$, $y = s$, $z = t$, and $v = w = t^{-2}$) implies that
 $$ \mbox{$\hamilton$ contains $Q + [t]Q$, $Q - [t]Q$, and $2Q$.}$$

If $|t|$ is even, then 
 $$ Q = H_- - \sum_{i=1}^{|t|-2} [t^i]Q
 = H_- - \sum_{j=1}^{(|t|-2)/2} [t^{2j-1}] \bigl( Q + [t]Q \bigr)
 \in \hamilton ,$$
 as desired.
 
 If $|t|$ is odd, define the hamiltonian cycle
 \begin{equation} \label{PrismOverChords4cyc-H}
 H = \Bigl( t^{q-1},s,t^{-(q-1)},u, \bigl( s^{-1},t,s,t \bigr)^{(|t|-q-1)/2},
s^{-1}, t \Bigr).
 \end{equation}
 We have
 $$H = (t^{|t|}) + [s](u,t^{-q}) - \sum_{i=0}^{(|t|-q-1)/2} [t^{q-1+2i}]Q $$
 and
 $$ H_+ = H_- + [s](u,t^{-q}) - [st](u,t^{-q}) - 2Q - [t]Q ,$$
 so
 \begin{align*}
  H - [t]H &= [s](u,t^{-q}) - [st](u,t^{-q}) - \sum_{i=0}^{(|t|-q-1)/2}
[t^{q-1+2i}] \bigl( Q - [t]Q \bigr) \\
 &\equiv [s](u,t^{-q}) - [st](u,t^{-q}) \pmod{\hamilton} ,
 \end{align*}
 and therefore
 $$ \bigl( H - [t]H \bigr) - H_+ \equiv - H_- + 2Q + [t]Q 
 \equiv Q \pmod{\hamilton} .$$
 Because $\bigl( H - [t]H \bigr) - H_+$ obviously belongs to~$\hamilton$, we
conclude that $Q \in \hamilton$, as desired.
 \end{proof}

\begin{prop} \label{EvenHeightOverChords}
 If there exist $s,t \in S$, such that
 \begin{itemize}
 \item $|G/G'|$ is even,
 \item $G' = \langle t \rangle$, and
 \item $|S'| = 4$,
 \end{itemize}
 then $\hamilton = \even$.
 \end{prop}

\begin{proof}
 Let
 \begin{itemize}
 \item $u \in S \setminus \{s^{\pm1}, t^{\pm1}\}$, so $S = \{s^{\pm1},
t^{\pm1}, u^{\pm1}\}$, 
 \item $m = |G/G'|$, and
 \item $Q = (t,s,t^{-1},s^{-1})$.
 \end{itemize}
 We may write $u=t^q$ for some $q$. 

Let us first establish that $Q \equiv -Q \equiv [v]Q \pmod{\hamilton}$, for all
$v \in G$.
 \begin{itemize}
 \item If $m = 2$, then Lemma~\ref{PrismOverChords4cyc} asserts $Q \in
\hamilton$, which is a stronger statement.
 \item If $m > 2$, then $s^2 \notin \langle t \rangle$, so the desired conclusion
is obtained by applying Proposition~\ref{4cyc=} to the spanning subgraph $\Cay
\bigl( G; \{s^{\pm1}, t^{\pm1}\} \bigr)$, with the roles of $s$ and~$t$
interchanged.
 \end{itemize}

Now let us now show that 
 \begin{equation} \label{EvenHeightOverChordsPf-2tn}
 2(t^{|t|}) \in \hamilton .
 \end{equation}
 We may assume, for the moment, that $q < |t|/2$, by replacing $u$ with its
inverse if necessary.
 (Note that, because $|S'| = 4$, we know $u$~is not an involution, so $q \neq
|t|/2$.)
 Let 
 \begin{equation} \label{EvenHeightOverChordsPf-C}
 C = 
 \begin{cases}
 \hfil 0 & \mbox{if $m = 2$} \\
 {[s]}\Bigl( \bigl( t,  \bigl( s, t^{|t|-2}, s, t^{-(|t|-2)} \bigr)^{(m-2)/2},
 t^{-1}, s^{-(m-2)} \Bigr)
 & \mbox{if $m > 2$}
 ,
 \end{cases}
 \end{equation}
 and define the hamiltonian cycles
 $$ H_1= \bigl( t^{|t|-3},s,t^{-(q-2)},u,t^{|t|-q-1},u,s^{-1},t^2 \bigr) - C
  $$
 and
 $$H_2= \bigl( t^{|t|-1}, s, u^{-1}, t^{q-1}, u^{-1}, t^{-(|t|-q-2)}, s^{-1}
\bigr) + C.$$
 Then
 $$ H_1 + H_2 = 2(t)^{|t|} - [t^{-1}]Q - [t^{-3}]Q \equiv 2(t^{|t|})
\pmod{\hamilton} ,$$
 which proves \pref{EvenHeightOverChordsPf-2tn}.

 Let 
 $$ X^* = \Cay \bigl( G; \{s^{\pm1}, t^{\pm1} \} \bigr) .$$

 \setcounter{proofcase}{0}

\begin{proofcase}
 Assume $|s| > 2$.
 \end{proofcase}
 It suffices to show $\even(X^*) \subseteq \hamilton + \hamilton(X^*)$, for then
Proposition~\ref{G=G'big} applies (with $u$~in the role of~$s$).
 \begin{itemize}
 \item If $X^*$~is not the square of an even cycle and $|G|$ is divisible
by~$4$, then Proposition~\fullref{degree4}{4} implies $\even(X^*) \subseteq
\hamilton(X^*) \subseteq \hamilton + \hamilton(X^*)$.
 \item If $X^*$~is not the square of an even cycle and $|G|$~is not divisible
by~$4$, then Corollary~\ref{weird+2tn} (with $u$ in the role of~$s$) implies
$\even(X^*) \subseteq \hamilton + \hamilton(X^*)$.
 \item If $X^*$ is the square of an even cycle, then $s^2 = t^{\pm1} \in
G'$, so $m = 2$. Therefore, Lemma~\ref{PrismOverChords4cyc} implies that
$\hamilton$ contains a basic $4$-cycle of~$X^*$. Therefore,
Corollary~\ref{deg4+4cyc} (with $u$ in the role of~$s$) implies $\even(X^*)
\subseteq \hamilton + \hamilton(X^*)$.
 \end{itemize}

 \begin{proofcase}
 Assume that $|s| = 2$, and that either $|t|$ or~$q$ is odd.
 \end{proofcase}
  If $|t|$ is odd, we may assume that $q$ is even (by replacing $u$ with its
inverse, if necessary). Thus, $|t|$~and~$q$ are of opposite parity. Define $H$ as
in~\pref{PrismOverChords4cyc-H}. Because $H$ has only a single $u$-edge, we can
eliminate all of the $u$-edges from any flow in~$\flow$, by adding appropriate
translates of the hamiltonian cycle~$H$, leaving us with a flow in
$\flow(X^*)$.
 Therefore $\flow \subseteq \hamilton + \flow(X^*)$, so
 \begin{equation} \label{EvenHeightOverChordsPf-2FinH+2F*}
 2\flow \subseteq \hamilton + 2\flow(X^*)
 .
 \end{equation}
 But any flow in $\flow(X^*)$ is a linear combination of cycles of the forms
$[v]Q$ and $[v](t)^{|t|}$. (Here, we use the assumption that $|s| = 2$.)
Therefore $2 \flow(X^*) \subseteq \hamilton$. Combining this with
\pref{EvenHeightOverChordsPf-2FinH+2F*} and \pref{ALWthm}, we conclude that
$\even \subseteq \hamilton$, as desired.

 \begin{proofcase}
 Assume that $|s| = 2$, and that both $|t|$ and $q$ are even.
 \end{proofcase}
 Let us begin by showing that $\hamilton$ contains every basic $4$-cycle. (Recall
that we already know $\hamilton$ contains the basic $4$-cycle~$Q$.)
 We may assume (by replacing $u$ with its inverse if necessary) that  $2 \le q <
|t|/2$. 
 \begin{itemize}
 \item Define the hamiltonian cycle 
 $$H_+= \Bigl( \bigl( s,t,s^{-1},t \bigr)^{|t|/2} \Bigr).$$
 This contains both the oriented path  $(s,t,s^{-1})$ and the oriented edge $[s
u](t)$ (because $q$~is even), so
Lemma~\ref{4c}(\ref{4c-s},\ref{4c-change}) (with $x = s$, $y = t$,
$z = u$, and $v = e$) implies that $\hamilton$ contains the basic $4$-cycle
$(u,t,u^{-1},t^{-1})$.
 \item Define the hamiltonian cycle 
 $$H_- = \bigl(
t,s,t^{-(|t|-q-1)},u^{-1},t^{q-1},s^{-1},t^{-(q-1)},u,t^{|t|-q-2} \bigr).$$
 This contains both the oriented path $(t,s,t^{-1})$ and the oriented edge $[t s
u](s^{-1})$, so Lemma~\ref{4c}(\ref{4c-d},\ref{4c-change}) (with
$x = t$, $y = s$, $z = u$, and $w = e$) implies that $\hamilton$ contains the
basic $4$-cycle $(u,s,u^{-1},s^{-1})$.
 \end{itemize}
 We have
 \begin{align*}
 \even &\subseteq \hamilton + \even' && \mbox{\see{E=H+E'(G'<>G)}} \\
  &\subseteq \hamilton + \hamilton' && \mbox{\see{deg4+4cyc}} \\
  &\subseteq \hamilton && \mbox{\fullsee{H'(G'even)}{1}} 
 ,
 \end{align*}
 as desired.
 \end{proof}

\begin{prop}\label{cycle+diam+chords}
 If $|S| = 5$ and $\langle t \rangle = G$, for some $t \in S$, then
$\hamilton=\even$.
 \end{prop}

\begin{proof}
 Let 
 \begin{itemize}
 \item $n = |G|/2$, and
 \item $u \in S \setminus \{t^{\pm1}\}$, with $|u| \neq 2$. We have $u =
t^q$, for some $q \neq n$. 
 \end{itemize}
 Note that $S = \{t^{\pm1}, u^{\pm1}, t^n\} = \{t^{\pm1}, t^{\pm q}, t^n\}$.
 
\setcounter{proofcase}{0}

\begin{proofcase}
 Assume $n$ is even.
 \end{proofcase}
 Let $s = u$, so $X'$ is a non-bipartite M\"{o}bius ladder. From
Theorem~\ref{cubic}, we have $\hamilton' = \even'$, so Lemma~\ref{lotsofchords}
applies.

\begin{proofcase}
 Assume $n$ and $q$ are both odd.
 \end{proofcase}
 Note that $X$ is bipartite. Let $s = t^n$ be the involution in~$S$.
 We have 
 \begin{align*}
 \even &\subseteq \hamilton + \even' && \mbox{\see{E=H+E'(G'=G)}} \\
 &\subseteq \hamilton + \hamilton' && \mbox{\see{deg4bip}} \\
 &\subseteq \hamilton && \mbox{\see{G=G'->H'inH}} 
 ,
 \end{align*}
 as desired.

\begin{proofcase}
 Assume $n$ is odd and $q$ is even.
 \end{proofcase}
  We may assume $2 \le q \le n-1$, by replacing $u$ with $u^{-1}$ if necessary.
Let $s = t^n$ be the involution in~$S$.

 Define 
 $$ \mbox{$Q = (t,s,t^{-1},s)$, \qquad $C = (t^n,s)$, \qquad $H = (t,s)^n$} ,$$
 so $Q$~is a basic $4$-cycle, $C$~is a cycle, and $H$~is a hamiltonian cycle
in~$X$.
 Note that
 $$ C - [t]C = Q .$$

Define the hamiltonian cycles 
 \begin{align*}
 H_+ &= (s,t^{n-1},s,t^{-(n-1)}) \\
 \intertext{and} 
 H_- &= 
 \begin{cases}
 (t^{n-3},u^{-1},t,s,t^{-1},u,t^{-(n-3)},s) &\mbox{if $q=n-1$} , \\
(u,t^{-(q-2)},s,t^{n-3},s,t^{-(n-q-2)},u^{-1},s,t^{-1},s) &\mbox{otherwise} 
 .
 \end{cases}
 \end{align*}
 Then
 \begin{itemize}
 \item $H_+$ contains both the oriented path $[t^{-2}](t,s,t^{-1})$ and the
oriented edge $[e](s)$, and
 \item $H_-$ contains both the oriented path $[t^{-2}](t,s,t^{-1})$ and the
oriented edge $[s](s^{-1})$,
 \end{itemize}
 so Lemma~\ref{4c}(\ref{4c-s},\ref{4c-d}) (with $x = t$, $y =
s$, $z = t$, and $v = w = t^{-2}$) implies that
 $$\mbox{$\hamilton$ contains $Q+[t]Q$ and $Q-[t]Q$.} $$
 Therefore
 $$\mbox{$2Q \in \hamilton$ \qquad and \qquad $[v]Q \equiv Q \pmod{\hamilton}$,
for all $v \in G$} .$$

Let us show that 
 \begin{equation} \label{cycle+diam+chords-QinH}
 Q \in \hamilton
 .
 \end{equation}
 We have
 $$ 
 C = H + \sum_{i=1}^{(n-1)/2} [t^{2i-1}]Q
 \equiv \frac{n-1}{2} Q
 \pmod{\hamilton} ,$$
 so
 $$ Q = C - [t]C
 \equiv \frac{n-1}{2} \bigl( Q - [t]Q \bigr)
 \equiv 0
 \pmod{\hamilton} ,$$
 as claimed.

The hamiltonian cycle~$H$ contains both the oriented path $[t](s,t,s)$ and the
oriented edge $[tsu](t)$ (because $q$ is even), so Lemma~\fullref{4c}{s}
(with $x = s$, $y = t$, $z = u$, and $v = t$) implies that 
 $$ (s,t,s^{-1},t^{-1}) + [s](u, t, u^{-1}, t^{-1}) \in \hamilton .$$
 Then, from \pref{cycle+diam+chords-QinH},
we know that $(u, t, u^{-1}, t^{-1}) \in \hamilton$. Hence, Lemma~\ref{deg4+4cyc}
implies $\even' \subseteq \hamilton + \hamilton'$, so Lemma~\ref{lotsofchords}
implies $\hamilton=\even$, as desired.
 \end{proof}

\section{Graphs of degree at least 5} \label{Degree5Sect}

In this section, we show that if $|S| \ge 5$, then $\hamilton = \even$
\see{degree7}. After establishing the special cases with $|S| = 5$ or~$6$
\seeand{degree5}{degree6}, it is very easy to complete the proof by induction
on~$|S|$.

\begin{prop} \label{degree5}
 If $|S| = 5$, then $\hamilton = \even$.
 \end{prop}

\begin{proof}
 
\setcounter{proofcase}{0}

\begin{proofcase}
 Assume some involution in~$S$ is redundant.
 \end{proofcase}
 We may assume there is a redundant involution~$s$ in~$S$, such that $\hamilton'
\neq \even'$, for otherwise Proposition~\ref{G=G'big} or Lemma~\ref{G=G'invols}
applies. (Note that $G' = G$, because $s$ is redundant.) Then
Proposition~\ref{degree4} implies that either
 \begin{itemize}
 \item
 $X'$ is the square of an even cycle, in which case,
Proposition~\ref{cycle+diam+chords} applies, or 
 \item $X'$ is not bipartite, and $|G'|$~is not divisible by~$4$.
 \end{itemize}
 Thus, we may assume that $X'$ is not bipartite, and $|G'|$~is not divisible
by~$4$.

 We have $|G| = |G'| \equiv 2 \pmod{4}$ and, by Lemma~\ref{Z2bipartite}, $S' =
\{t^{\pm1}, u^{\pm1}\}$, where $|t|$~is odd and $|u|$~is even.
 If $G = \langle u \rangle$, then Proposition~\ref{cycle+diam+chords} applies. If
not, then, because $|G|$~is not divisible by~$4$,
Corollary~\fullref{GridsHave4Cycles}{not4} (with the roles of $s$ and~$t$
interchanged) implies that $\hamilton$ contains every basic $4$-cycle~$C$. Hence, 
 \begin{align*}
 \even &\subseteq \hamilton + \even' && \mbox{\see{E=H+E'(G'=G)}} \\
 &\subseteq \hamilton + \hamilton' && \mbox{\see{deg4+4cyc}} \\
 &\subseteq \hamilton && \mbox{\see{G=G'->H'inH}}
 ,
 \end{align*}
 as desired.

\begin{proofcase}
 Assume every involution in~$S$ is irredundant.
 \end{proofcase}
 We may assume that $S$ contains only one involution~$s$ \see{noirredinvols}.
Then $X$~is the prism over~$X'$ (and $S'$ does not contain any involutions).
 Furthermore, we may assume 
 \begin{equation} \label{degree5pf-t<>G'}
 \mbox{$\langle t \rangle \neq G'$, for all $t \in S'$}
 \end{equation} (otherwise, Proposition~\ref{EvenHeightOverChords} applies). In
particular, this implies that $X'$ is not the square of an even cycle.

\begin{subcase}
 Assume $|G'|$ is odd.
 \end{subcase}
 Let $t \in S'$. We know $\langle t
\rangle \neq G'$ and $|G| = 2 |G'|$ is not divisible by~$4$, so
Corollary~\fullref{GridsHave4Cycles}{not4} implies that $\hamilton$ contains
every basic $4$-cycle. We may assume that $X' \not\iso K_3 \boxprod K_3$
(otherwise, Corollary~\ref{K2xK3xK3} applies), so Theorem~\ref{oddorder} asserts
that $\hamilton' = \flow'$.
 We have
 \begin{align*}
 \even &\subseteq \hamilton + 2\flow' && \mbox{\fullsee{E=H+2F'}{X'notbip}} \\
 &\subseteq \hamilton + 2 \hamilton' && \mbox{(because $\hamilton' = \flow'$)} \\
 &\subseteq \hamilton && \mbox{\see{2H'inH(G'odd)deg5}} 
 ,
 \end{align*}
 as desired.

\begin{subcase}
 Assume that $|G'|$ is even and that either $X'$~is bipartite or $|G'|$~is
divisible by~$4$.
 \end{subcase}
 For any $t \in S'$, we have
 $$ \frac{|G|}{|t|}
 = \frac{|G|}{|G'|} \cdot \frac{|G'|}{|t|}
 \ge 2 \cdot 2
 = 4
 ,$$
 so Corollary~\fullref{4cyc}{|G/G'|>2} (with the roles of
$s$~and~$t$ interchanged) implies that the basic $4$-cycle
$(s,t,s^{-1},t^{-1})$ is in~$\hamilton$. Therefore
 \begin{align*}
 \even & \subseteq \hamilton + \even' && \mbox{\see{E=H+E'(G'<>G)}} \\
  & \subseteq \hamilton + \hamilton' && \mbox{\see{degree4}} \\
  & \subseteq \hamilton && \mbox{\fullsee{H'(G'even)}{1}} 
 ,
 \end{align*}
 as desired.

\begin{subcase}
 Assume that $X'$~is not bipartite and that $|G'|$~is even, but not divisible
by~$4$.
 \end{subcase}
 Since $X'$ is not the square of an even
cycle, $X'$~must be as described in Proposition~\ref{WeirdCase}.
 Write $S' = \{t^{\pm1}, u^{\pm1}\}$, with $|t|$~odd and $|u|$~even. 
 \begin{itemize}
 \item We know $\langle u \rangle \neq G'$ \see{degree5pf-t<>G'}.
  \item Because $S'$ does not contain any involutions, we know $|u| \neq 2$.
\item Because $|G'|/|u|$ is odd, we know $|G'|/|u| \neq 2$.
 \end{itemize}
 Thus, Lemma~\fullref{GridsHave4Cycles}{no2} (with $u$ in the role of~$t$)
implies that $\hamilton$ contains every basic $4$-cycle. 
 We have
 \begin{align*}
 \even &\subseteq \hamilton + \even' && \mbox{\see{E=H+E'(G'<>G)}} \\
 &\subseteq \hamilton + \hamilton' && \mbox{\see{deg4+4cyc}} \\
 &\subseteq \hamilton && \mbox{\fullsee{H'(G'even)}{1}} 
 ,
 \end{align*}
 as desired.
 \end{proof}

\begin{prop} \label{degree6}
 If $|S| = 6$, then $\hamilton = \even$.
 \end{prop}

\begin{proof}
 Let us begin by establishing that we may assume 
 \begin{equation} \label{degree6pf-noinvol}
  \mbox{there are no involutions in~$S$.} 
 \end{equation}
 First, note that if $s$ is any redundant involution, then $\hamilton' = \even'$
 (by Proposition~\ref{degree5}), so Proposition~\ref{G=G'big} or
Lemma~\ref{G=G'invols} implies $\hamilton = \even$, as desired. On the other
hand, if all of the involutions in~$S$ are irredundant, then
Observation~\ref{noirredinvols} asserts that $X$ can be realized by a generating
set with at most one involution. Because $X$ has even degree, there must be no
involutions.

\setcounter{proofcase}{0}

\begin{proofcase}
 Assume $S$ is irredundant.
 \end{proofcase}
 Choose $s \in S$, such that 
 $$ \mbox{$|G'|$~is even.} $$
 We know, from Corollary~\ref{4DGridsHave4Cycles}, that 
 \begin{equation} \label{degree6-basic}
 \mbox{$\hamilton$ contains every basic $4$-cycle.} 
 \end{equation}
 We have 
 \begin{align*}
 \even &\subseteq \hamilton + \even' && \mbox{\see{E=H+E'(G'<>G)}} \\
 & \subseteq \hamilton + \hamilton' && \mbox{\seeand{degree6-basic}{deg4+4cyc}} \\
 & \subseteq \hamilton && \mbox{\fullsee{H'(G'even)}{1}}
 , \\
 \end{align*}
 as desired.

\begin{proofcase}
 Assume there is a redundant generator~$s$ in~$S$.
 \end{proofcase}
 It suffices to show $\even' \subseteq \hamilton + \hamilton'$, for then
Proposition~\ref{G=G'big} applies. Note that,
from \pref{degree6pf-noinvol}, we know $|S'| = 4$.

\begin{subcase}
 Assume that $X'$ is not the square of an even cycle, and that either $X'$~is
bipartite, or $|G'|$~is divisible by~$4$.
 \end{subcase}
 Proposition~\ref{degree4} (and the assumption of this subcase) implies that
$\hamilton' = \even'$, so $\even' \subseteq \hamilton + \hamilton'$, as desired.

\begin{subcase}
 Assume $X'$ is the square of an even cycle.
 \end{subcase}
 By choosing~$s$ to be the chord of length~$2$, we may move out of this subcase.

\begin{subcase}
 Assume that $X'$~is not bipartite, and that $|G'|$~is not divisible by~$4$.
 \end{subcase}
 We have  $|G| = |G'| \equiv 2 \pmod{4}$, and, by Lemma~\ref{Z2bipartite}, $S' =
\{t^{\pm1}, u^{\pm1} \}$, where $|t| \equiv |s| \pmod{2}$
 and $|u| \not\equiv |s| \pmod{2}$. 

 It suffices to show that $\hamilton$ contains some basic $4$-cycle of~$X'$, for
then Corollary~\ref{deg4+4cyc} yields
 $\even' \subseteq \hamilton + \hamilton'$, 
 as desired.

\begin{subsubcase}
 Assume $\langle s,t \rangle = G$.
 \end{subsubcase}
 Let $X^* = \Cay \bigl( G; \{s,t\} \bigr)$. Now $X^*$ is bipartite
\see{Z2bipartite}, so $\hamilton(X^*) = \even(X^*)$ \see{deg4bip}. Therefore
Proposition~\ref{G=G'big} (with $u$~in the role of~$s$) implies $\hamilton =
\even$.

\begin{subsubcase}
 Assume  $\langle s,t \rangle \neq G$.
 \end{subsubcase}
 \begin{itemize}
 \item If $\langle s \rangle \neq \langle t \rangle$, then
Corollary~\fullref{GridsHave4Cycles}{not4} applies (with $u$ in the role of~$s$,
and one or the other of $s$ and~$t$ in the role of~$t$), so $\hamilton$ contains
every basic $4$-cycle.
 \item If $|G/\langle s,t \rangle|$ is odd, then Corollary~\ref{OddHeight4Cycles}
(with $u$ in the role of~$s$) implies that $\hamilton$ contains every basic
$4$-cycle.
 \item If $\langle s \rangle = \langle t \rangle$ and $|G/\langle s,t \rangle|$
is even, then Proposition~\ref{EvenHeightOverChords} (with $u$~in the role
of~$s$) implies $\hamilton = \even$.
 \end{itemize}
 This completes the proof.
 \end{proof}

\begin{cor} \label{degree7}
 If $|S| \ge 5$, then $\hamilton = \even$.
 \end{cor}

\begin{proof}
 By Propositions~\ref{degree5} and~\ref{degree6}, we may assume $|S| \ge 7$, so,
by induction on $|S|$, we have
 \begin{equation} \label{degree7-E'inH'}
 \mbox{$\hamilton' = \even'$, for every $s \in S$ with $|G'|$ even.}
 \end{equation}
 
 We may assume 
 $$ \mbox{$S$ is irredundant,} $$
 for otherwise Proposition~\ref{G=G'big} or Lemma~\ref{G=G'invols} implies
$\hamilton = \even$, as desired.

 Choose $s \in S$, such that 
 $$ \mbox{$|G'|$ is even.} $$
 Because $S$ is irredundant, we know $G' \neq G$, and, by
Observation~\ref{noirredinvols}, we may assume $S$ has no more than one
involution. Then, from Corollary~\ref{4DGridsHave4Cycles}, we know that
$\hamilton$ contains every basic $4$-cycle. We have
 \begin{align*}
 \even & \subseteq \hamilton + \even' && \mbox{\see{E=H+E'(G'<>G)}} \\
 &\subseteq \hamilton + \hamilton' && \mbox{\see{degree7-E'inH'}} \\
 &\subseteq \hamilton &&  \mbox{\fullsee{H'(G'even)}{1}}
 ,
 \end{align*}
 as desired.
 \end{proof}

Corollary~\ref{degree7} completes the proof of
Theorem~\ref{mainthm}. The Cayley graphs of degree~$4$ were
considered in \pref{degree4}, \pref{SquareCycle}, and
\pref{WeirdCase}.

\end{document}